\documentclass[10pt]{amsart}

\usepackage{latexsym}
\usepackage{amsmath}
\usepackage{amsfonts}
\usepackage{amssymb}
\usepackage{amscd}
\usepackage{esint}
\usepackage{graphicx}
\usepackage{pstricks}

\AtBeginDocument{%
   \def\MR#1{}
}

\usepackage{pstricks,pstcol,pst-plot,pst-node,pst-tree,amssymb}
\usepackage{pstricks-add}
\renewcommand{\a}{\alpha }
\renewcommand{\b }{\beta }
\renewcommand{\d}{\delta }
\newcommand{\D }{\Delta }

\newcommand{\e }{\varepsilon }
\newcommand{\g }{\gamma}

\renewcommand{\l}{\lambda }
\renewcommand{\L }{\Lambda }

\newcommand{\n }{\nabla }
\newcommand{\var}{\varphi }

\newcommand{\s }{\sigma }

\renewcommand{\t }{\tau }
\renewcommand{\th }{\theta }
\renewcommand{\o }{\omega }

\newcommand{\oz }{{\ov{z}}}

\newcommand{\pa }{\partial}

\newcommand{\ov}{\overline}

\newcommand{\be}{\begin{equation}}
\newcommand{\ee}{\end{equation}}
\newenvironment{pf}{\noindent{\sc Proof}.\enspace}{\rule{2mm}{2mm}\medskip}
\newenvironment{pfn}{\noindent{\sc Proof}}{\rule{2mm}{2mm}\medskip}

\newcommand{\R}{\mathbb{R}}
\newcommand{\C}{\mathbb{C}}
\renewcommand{\H}{\mathbb{H}}

\newcommand{\ou}{\overline{1}}

\begin{document}

\newtheorem{lem}{Lemma}[section]
\newtheorem{pro}[lem]{Proposition}
\newtheorem{thm}[lem]{Theorem}
\newtheorem{rem}[lem]{Remark}
\newtheorem{cor}[lem]{Corollary}
\newtheorem{df}[lem]{Definition}

\bibliographystyle{amsalpha}

\title[] {On the Sobolev quotient of three-dimensional CR manifolds}

\author{Jih-Hsin Cheng$^{(1)}$, Andrea Malchiodi$^{(2)}$,  Paul Yang$^{(3)}$}

\address{}

\thanks{}

\email{}




\maketitle

\begin{center}

{\small }

\end{center}

\footnotetext[1]{E-mail addresses: cheng@math.sinica.edu.tw, andrea.malchiodi@sns.it, yang@math.princeton.edu}

\vspace{-1cm}

\begin{center}

{\small $^{(1)}$ Institute of Mathematics, Academia Sinica and NCTS 

6F, Astronomy-Mathematics Building
No. 1, Sec. 4

Roosevelt Road,
Taipei 10617, TAIWAN

\

$^{(2)}$ Scuola Normale Superiore, 
Piazza dei Cavalieri 7, 
50126 Pisa, 
ITALY

\

$^{(3)}$Princeton University, Department of Mathematics

Fine Hall, Washington Road, 
Princeton NJ 08544-1000 USA

}

\end{center}

\begin{abstract}
We exhibit examples of compact three-dimensional CR manifolds of positive Webster class, {\em Rossi spheres},  
for which the pseudo-hermitian mass as defined in \cite{CMY17} is negative, and for which the infimum of the 
CR-Sobolev quotient is not attained. To our knowledge, this is the first geometric context on smooth closed manifolds
where this phenomenon arises,  in striking contrast to the Riemannian case. 
\end{abstract}

\section{Introduction}

The Yamabe problem  consists in deforming conformally 
the metric of a manifold of dimension $n \geq 3$ so that its scalar curvature becomes a constant. 
Apart from being a natural conformal extension of the Uniformization Problem in two dimensions, 
the question was introduced in \cite{Yam60} for trying to attack Poincar\'e's conjecture. 
Yamabe metrics have also been applied to other contexts, such as the study of 
 degeneration of conformal structures.  For example, in \cite{TV05} it is shown that the set of Yamabe Bach-flat metrics on a 
 four-manifold is compact up to orbifold degeneration.

Writing on $(M,g)$ a conformal metric as $\tilde{g} = u^{\frac{4}{n-2}} g$, the scalar curvature 
transforms as 
\begin{equation}\label{eq:trans-R}
- \frac{4(n-1)}{n-2} \D_g u + S_g u = S_{\tilde{g}} \, u^{\frac{n+2}{n-2}}. 
\end{equation}
Therefore, if one wishes to have $S_{\tilde{g}}$ constant, the 
following elliptic problem must be solved 
\begin{equation}\label{eq:Y} \tag{$Y$}
 - \frac{4(n-1)}{n-2} \D_g u + S_g u = \ov{S} \, u^{\frac{n+2}{n-2}} \quad 
 \hbox{ on } M; \qquad \quad \ov{S} \in \R. 
\end{equation}
Notice that the exponent on the right-hand side of the equation is critical with respect to 
the Sobolev embeddings. In \cite{Yam60} an attempt was made to solve \eqref{eq:Y} by 
lowering the exponent by a small amount, in order to obtain compactness, and then by 
letting it approach the critical one studying the limit of the corresponding solutions. 
The problem with this strategy though is that the weak limit of such solutions might 
be zero. 
Another way to attack \eqref{eq:Y} was to view $\ov{S}$ as a Lagrange multiplier, considering the Sobolev quotient
\begin{equation}\label{eq:Q}
Q_{(M,g)}(u) := \frac{\int_M \left( c_n |\n_g u|^2 + S_g u^2 \right) dV_g}{\left( \int_M |u|^{2^*} dV_g \right)^{\frac{2}{2^*}}} 
= \frac{\int_M S_{\tilde{g}} dV_{\tilde{g}}}{\left( Vol_{\tilde{g}}(M) \right)^{\frac{2}{2^*}}}, 
\end{equation}
where $c_n = \frac{4(n-1)}{(n-2)}$ and $2^* = \frac{2n}{n-2}$. 
If one could realise the minimum of $Q_{(M,g)}(u)$ over all  non-zero $u$'s of class $W^{1,2}(M,g)$, this would give rise to 
a solution of \eqref{eq:Y}: notice that it is sufficient to consider functions 
in $W^{1,2}(M,g)$ that are non-negative, therefore by regularity theory one would obtain a 
positive smooth solution.  Defining then 
$$
  Y(M,g) := \inf_{u \in W^{1,2}(M,g), u \not\equiv 0} Q_{(M,g)}(u), 
$$
it can be proved using \eqref{eq:trans-R} that this quantity is independent of the 
conformal representative of $g$, and will therefore be denoted from now on by 
$Y(M,[g])$. Depending on the sign of the latter quantity, $(M,[g])$ is said to be 
of {\em negative, null} or of {\em positive Yamabe class}. 

It was proved in \cite{Tru68} that there exists a dimensional constant $\e _n > 0$ such that 
$Y(M,[g])$  is attained (and hence \eqref{eq:Y} is  solvable) provided $Y(M,[g]) \leq \e_n$. The 
result applies in particular to all manifolds with conformal classes of metrics of 
negative or null Yamabe class. 

Consider the (normalized) Sobolev quotient in $\R^n$ 
\begin{equation}\label{eq:Sn}
  S_n := \inf_{u \in C^\infty_c(\R^n), u \not\equiv 0} \frac{\int_{\R^n} c_n |\n u|^2 dx}{\left( \int_{\R^n} |u|^{2^*} dx \right)^{\frac{2}{2^*}}}.
\end{equation}
Using the stereographic projection from $S^n$ to $\R^n$ it can be proved that  the 
above quantity coincides with the Yamabe quotient of the round sphere, i.e. for 
all $n \geq 3$ one has 
$$
S_n =  Y(S^n,[g_{S^n}]). 
$$
It was shown in \cite{Aub76} that one always has $Y(M,[g]) \leq S_n$, and that \eqref{eq:Y} 
is solvable provided the strict inequality holds. It was also 
shown in \cite{Aub76} that  $Y(M,[g]) < S_n$ provided $n \geq 6$ and $M$ is 
not locally conformally flat, i.e. when the Weyl tensor of $(M,g)$ is not identically zero. 
It was proved then in  \cite{Scho84} that $Y(M,[g]) < S_n$  in all 
complementary cases (provided $(M,g)$ is not conformally equivalent to the round sphere), 
i.e. when $(M,g)$ has dimension less or equal to $5$ or when it is locally conformally flat. 
While the argument in \cite{Aub76} was based on a local energy expansion, the one 
in \cite{Scho84} relied on the {\em Positive Mass Theorem} in general relativity, see 
\cite{SY-PM}, \cite{SYPM2}, \cite{SY-LCF}, \cite{SY17}, which is in turn related to the 
expansion of the Green's function of the {\em conformal Laplacian} $L_g$ near 
its pole, where  
$$
 L_g u:= 
  - \frac{4(n-1)}{n-2} \D_g u + S_g u. 
$$
In both \cite{Aub76} and \cite{Scho84} the strict inequality was proved by evaluating the 
Yamabe-Sobolev quotient on (suitable perturbations of) highly concentrated 
extremals of \eqref{eq:Sn} (classified in \cite{Aub76}, \cite{Tal76}), suitably glued to 
$(M,g)$. Such extremals, parametrized using the M\"obius group of 
$S^n$, can be chosen arbitrarily peaked near any point: these decay faster 
at infinity in higher dimensions and therefore the correction to the quotient 
due to the geometry of $M$ is more \emph{localized} in space for $n$ large. 
In any case, we always have 
$$
  Y(M,[g]) < S_n \qquad \quad \hbox{ provided } \qquad \quad (M,g) \stackrel{conf.}{\not\simeq} (S^n,g_{S^n}). 
$$

\

We consider in this paper  compact three dimensional pseudo-hermitian manifolds $(M,J,\th)$: 
these are CR manifolds, i.e. endowed with a contact structure $\xi$ and a CR structure $J : \xi \to \xi$
such that $J^2 = - 1$. We assume $(M,J)$ to be {\em pseudo-convex}, namely that it is 
globally defined a 
contact form $\th $ which annihilates $\xi$ and for which $\th \wedge d \th $
is always non-zero (see \cite{BFG83}). We define the \emph{Reeb vector field} as the unique
 $T$ for which $\th (T) \equiv 1$ and $T \lrcorner \; d \th = 0$.
Given $J$ as above, we can define  locally  a vector field $Z_1$ such that 
\begin{equation}  \label{eq:J0}
J Z_1 = i Z_1; \qquad J Z_{\overline{1}} = - i Z_{\overline{1}} \qquad \quad %
\hbox{ where } \quad Z_{\overline{1}} = \overline{(Z_1)}.
\end{equation}
We also let $(\th , \th ^1, \th ^{\overline{1}})$ be the dual triple to $%
(T, Z_1, Z_{\overline{1}})$, so that 
\begin{equation*}
d \th = i h_{1 \overline{1}} \th ^1 \wedge \th ^{\overline{1}} \qquad 
\hbox{
for some } h_{1 \overline{1}} > 0 \quad (\hbox{possibly replacing } \theta %
\hbox{ by } -\theta).
\end{equation*}
In the following we will always assume that $h_{1 \overline{1}} \equiv 1$.

The \emph{connection} 1-form $\o ^1_1$ and the \emph{torsion} $A^1_{\overline{1}}$ are uniquely determined by
the structure equations 
\begin{equation}\label{eq:structure}
\left\{ 
\begin{array}{ll}
d \th ^1 = \th ^1 \wedge \o ^1_1 + A^1_{\overline{1}} \th \wedge \th ^{%
\overline{1}}; &  \\ 
\o ^1_1 + \o ^{\overline{1}}_{\overline{1}} = 0. & 
\end{array}
\right.
\end{equation}
The \emph{Tanaka-Webster curvature} (or {\em Webster curvature}) $R_\th$ (or, simply, $R$) is then defined by the formula 
\begin{equation*}
d \o ^1_1 = R_\th \, \th ^1 \wedge \th ^{\overline{1}} \; \qquad (\hbox{mod } \th ). 
\end{equation*}
A  model with positive curvature is the standard sphere $(S^3, J_{S^3}, \hat{\th})$, 
with $S^3 \subseteq \C^2 = \{(z_1, z_2)\}$, and 
\begin{equation}
\hat{\theta} = \frac 12 i(\bar{\partial}-\partial )(|z^{1}|^{2}+|z^{2}|^{2})
\label{B-1} = \frac 12 i\sum_{k=1}^{2}(z^{k}dz^{\bar{k}}-z^{\bar{k}}dz^{k}); 
\qquad 
Z_{1} = Z_{1}^{S^{3}}= \ov{z}^{2}\frac{\partial }{\partial z^{1}}- \ov{z}^1 \frac{\partial }{%
\partial z^{2}}.
\end{equation}
Similarly to what happens with the classical stereographic projection, the CR three-sphere 
is CR equivalent to the {\em Heisenberg group} $\H^1 = \left\{
(z,t), z \in \C, t \in \R \right\}$, see e.g. \cite{CMY17}.

The Tanaka-Webster curvature enjoys 
conformal properties similar to the  scalar curvature on Riemannian manifolds. 
More precisely, scaling the contact form $\theta$ by a positive function, 
one has the following law for the transformation of the Webster curvature, 
similar to \eqref{eq:trans-R}
 \begin{equation}\label{eq:trans-W}
   L_b u := - 4 \D_b u + R_\th \, u = R_{\tilde{\th}} \, u^{3}; \qquad \qquad 
   \tilde{\th} = u^2 \th . 
 \end{equation}
Here $R_{\tilde{\th}}$ is the Tanaka-Webster curvature corresponding to the
pseudo-hermitian structure $(J, \tilde{\th })$. $\D_b$ stands for the operator defined as follows 
\begin{equation*}
\D_b f = f,_1^{\; \; \; 1} + f,_{\overline{1}}^{\; \; \; \overline{1}} =
f_{, 1\overline{1}} + f_{, \overline{1} 1}, 
\end{equation*}
where we have used $h^{1 \overline{1}} = h_{1 \overline{1}} = 1$ to raise or
lower the indices, and where we  set 
\begin{equation}\label{eq:not-der}
f_1 = f,_{1} := Z_1 f; \qquad \quad f,_{1 \overline{1}} = Z_{\overline{1}%
} Z_1 f - \o ^1_1(Z_{\overline{1}}) Z_1 f; \qquad \quad f,_0 = T f. 
\end{equation}
The CR-invariant sub-Laplacian transforms  covariantly as follows 
\begin{equation*}
\hat{L}_b (\var) = u^{-\frac{Q+2}{Q-2}} L_b(u \var); \qquad \quad \hat{\th }
= u^2 \th ,
\end{equation*}
where $Q = 4$ is the \emph{homogeneous dimension} of the manifold. 
By \eqref{eq:trans-W}, finding $\tilde{\th}$ with constant Webster curvature 
corresponds to solving the following analogous problem to \eqref{eq:Y} 
\begin{equation}\label{eq:W} \tag{$W$}
 L_b u = \ov{R} \, u^{\frac{Q+2}{Q-2}} \quad 
 \hbox{ on } M; \qquad \quad \ov{R} \in \R, \quad u > 0. 
\end{equation}

In \cite{JLCRYam} the counterpart of the result in \cite{Aub76} was obtained, 
i.e. if  the infimum of the CR-Sobolev quotient satisfies 
\begin{equation}  \label{eq:Y(J)}
\mathcal{Y}(M,J) := \inf_{\hat{\th }} \frac{\int_M R_{J,\hat{\th }} \, \hat{\th }
\wedge d \hat{\th }}{\left( \int_M \hat{\th } \wedge d \hat{\th }
\right)^{\frac 12}} < \mathcal{Y}(S^3,J_{S^3}), 
\end{equation}
then it is attained and a solution of 
\eqref{eq:W} exists (indeed, this holds true in any dimension). The same authors verified this condition when the 
dimension is greater or equal to five and  $(M,J)$ is not {\em  spherical}, 
see \cite{JLCR} and \cite{JLExtr}. 

\

However, in the CR setting new phenomena appear, related to the fact that 
most three-dimensional structures are non-embeddable, 
differently from the higher-dimensional case, see \cite{Bou84}, \cite{BE90}. In \cite{CMY17} 
some results in the above directions were proved,  assuming embeddability of the 
structure.

More precisely, a notion of {\em pseudo-hermitian mass} was defined for three-dimensional {\em asymptotically-Heisenberg} 
manifolds  (we refer to the latter paper for precise definitions and details) by setting 
\begin{equation*}
m(J, \th ) := i \oint_{\infty} \o ^1_1 \wedge \th := \lim_{\L \to + \infty}
i \oint_{S_\L } \o ^1_1 \wedge \th ,
\end{equation*}
where  $S_\L = \left\{ \rho = \L \right\}$, $\rho^4 = |z|^4 + t^2$ (with $(z,t)$ coordinates on the Heisenberg group), 
and where $\o ^1_1$ stands for the connection form of the structure. 
The above definition was introduced considering an analogue of the 
\emph{Einstein-Hilbert action}.

As it happens in the Riemannian case, this mass is related to the expansion of the Green's 
function of the conformal sub-Laplacian $L_b$ on a compact manifold $M$. When    $\mathcal{Y}(M,J) > 0$ the latter operator is
invertible, so for any $p \in M$ there exists a Green's function $G_p$ verifying distributionally 
\begin{equation*}
\left( - 4 \D_b + R \right) G_p = 64 \, \pi^2 \, \d_p,
\end{equation*}
where $\d_p$ in the the right-hand side stands for the Dirac delta 
w.r.t. the volume measure $\th \wedge d \th$. In {\em CR normal coordinates} $(z,t)$ (introduced in \cite{JLCR} and discussed in Section \ref{s:CRcoord}) 
$G_p$ writes as  
\begin{equation}\label{eq:Green-A}
G_p = 2 \rho^{-2} + A + O(\rho),
\end{equation}
for some $A \in \R$  and where  $\rho^4(z,t)$ is as above. For the latter expansion, 
we refer to Proposition 5.2 in \cite{CMY17} (here we use an extra factor $4 \pi$ in the definition of $G_p$), and 
to Subsection \ref{ss:notprel} for our notation $O(\rho)$. 
Given $(M,J,\th)$ compact and $p \in M$, 
consider a  blow-up of contact form as follows 
\begin{equation}  \label{eq:bbuu}
N = (M \setminus \{p\}, J,  G_p^2 \th ). 
\end{equation}
As it is shown in \cite{CMY17}, via an inversion of coordinates, 
the manifold $N$ turns out to have 
asymptotically the geometry of the Heisenberg group, and 
its pseudo-hermitian mass satisfies
\begin{equation}\label{eq:mass-A}
  m = 12 \pi A 
\end{equation}
(see Lemma 2.5 there, and recall the difference of $4 \pi$ in our current notation), where $A$ is as above. 
Using crucially a result in \cite{HsYu}, in the same paper it was also proved that the pseudo-hermitian mass 
is non-negative (and zero only when $(M,J,\th)$ is CR equivalent to 
$S ^3$), provided that the {\em CR Paneitz operator} $P$ on $(M,J)$ is non-negative definite. 
The latter operator is 
\begin{equation}
P\varphi :=4(\varphi {_{\bar{1}}}^{\bar{1}}{_{1}}+iA_{11}\varphi ^{1})^{1}, 
\label{Pan}
\end{equation}
it has a relation to the  $\log $-term coefficient in the Szeg\"{o} kernel expansion, 
and it is pseudo-hermitian-covariant,  namely 
$
P_{\hat{\theta}}\varphi =e^{4f} P_{\theta }\varphi  
$
for the conformal change $\theta = e^{2f} \hat{\theta}$ (\cite{Hi}). 
By a result in \cite{CCY1}, manifolds for which $P$ is non-negative and $W > 0$ 
can be embedded into some 
$\C^N$ (see also  \cite{CCY2}). 

The assumption on the positivity on the Paneitz operator is not technical, as 
in \cite{CMY17} some counterexamples for the positivity of the pseudo-hermitian mass 
were also given for structures (arbitrarily) close to the spherical one, 
and hence with positive Webster curvature. 

\

In this paper we are concerned with {\em Rossi spheres}: these  are a one-parameter-family of CR structures on the 3-sphere of the form 
$S^3_s := (S^3,J_{(s)}, \hat{\th})$, where $\hat{\th}$ is as in \eqref{B-1}, and 
where $J_{(s)}$ is characterized by 
\begin{equation}\label{eq:Js}
J_{(s)} Z_{1(s)} = i Z_{1(s)}; \qquad \qquad Z_{1(s)}  = Z_{1}+\frac{s}{\sqrt{1+s^{2}}}Z_{\bar{1}},     \quad 
Z_{\bar{1}(s)}  = Z_{\bar{1}}+\frac{s}{\sqrt{1+s^{2}}}Z_{1}.  
\end{equation}
Rossi spheres are interesting because they are simple examples of CR structures on the three-sphere that cannot be embedded in $\mathbb C^N$. In \cite{Burns} it was shown that all the holomorphic functions on such structures are even functions if $s \neq 0$. On the other hand, there are explicit embeddings in $\mathbb C^3$ of the quotient of the Rossi spheres by the antipodal map, see \cite{C-S-Book}. 
By the above discussion, it follows that the Paneitz operator cannot be non-negative here. In addition, this family of CR structures are homogeneous and if we take the
 standard contact form, it is \emph{pseudo-Einstein}, i.e. $R_{,1} - iA_{11,\bar1}=0$, see \cite{CY13} as well as 
our notation for covariant derivatives in Section \ref{ss:notprel}.

\medskip 

Our first main result in this paper is the following  theorem. 

\medskip

\begin{thm}\label{t:mass}
For $|s|$ small, $s \neq 0$, the pseudo-hermitian mass of the Rossi spheres $S^3_s$ is negative.  
More precisely, one has the expansion 
$$
  m_s =  - 18 \pi s^2 + o(s^2) \qquad \qquad \hbox{ for } s \simeq 0.  
$$
\end{thm}

\

We saw before (in both low-dimensional Riemannian and CR cases) that positivity of the
mass implies attainment of the Sobolev quotient. We also strengthen the relation between 
mass and quotient by means of the following result.  

\medskip

\begin{thm}\label{t:SQ}
For $|s|$ small, $s \neq 0$, the infimum of the CR-Sobolev quotient of $S^3_s$ coincides with 
$\mathcal{Y}(S^3,J_{S^3})$ and  \underline{is not attained}. 
\end{thm}

\medskip

\begin{rem}\label{r:gamara}
(a) The phenomenon in Theorem \ref{t:SQ} is typical of some critical problems in a PDE context, 
like the Yamabe equation on Euclidean domains with Dirichlet boundary conditions or the case of some general elliptic operators 
on manifolds. However, to our knowledge this is the first time this is displayed in a \underline{purely geometric} smooth context.

(b) We recall that in \cite{Ga01} and \cite{GaYa01} the CR-Yamabe problem was  solved for every three dimensional 
CR manifolds, but there solutions were  found via variational arguments and they are not of minimal type. 
Theorem \ref{t:SQ} shows that the use of such methods is in some cases somehow necessary. 
\end{rem}

\medskip

Determining or estimating the {\em mass} of a manifold is in general a hard problem, 
since this is deeply related to the Green's function of the conformal (sub-)Laplacian, which 
is a {\em global} object. After recalling some preliminary facts in Section \ref{s:CRcoord} on 
CR normal coordinates (introduced in \cite{JLCR}) and on Rossi spheres, 
we specialize in Section \ref{s:CR-Rossi} to the latter manifolds, deriving first a suitable conformal 
factor and then expressing pseudo-hermitian coordinates depending on $s$. 
By the special expression of the Green's function in these coordinates, we are able to determine it
quite precisely near the north pole,  up to the constant term $A$ appearing in \eqref{eq:Green-A}. 

By a formal expansion in $s$, worked-out at the beginning of Section \ref{s:Green}, 
it is possible to characterize formally the Green's function for the conformal sub-Laplacian 
on Rossi spheres up to an order $O(s^3)$. However this expansion generates singular terms, 
with a particularly bad behavior near the pole, if expressed with respect to the standard 
complex coordinates of $\C^2$, where $S^3$ embeds. However we verify in the second 
part of the section that the global singular expansion on $S^3$ matches with the 
one done in CR normal coordinates up to an order $O(s^3)$, allowing us to prove Theorem \ref{t:mass}.

In Section \ref{s:min}, arguing by contradiction, we analyse the possible behaviours of 
minimizers for the CR Sobolev quotient. Due to a non-degeneracy result from \cite{MuUg}, 
the analysis of minimizers can be reduced to a finite-dimensional one, and we show that 
the CR-Sobolev quotient of all candidate minimizers is strictly above the spherical one, 
i.e. $\mathcal{Y}(S^3,J_{S^3})$.  With negative mass, this is expected for highly 
concentrated profiles, reversing the expansion in \cite{Scho84}: however such a 
property has to be obtained in \underline{all cases}, i.e. even for non-concentrated profiles.  
In Proposition \ref{p:expansion} this is proved for $s$ small in 
a \underline{fixed} compact set of the CR maps of $S^3$. One needs then to analyze the 
quotient in a regime with loss of compactness, which is  particularly delicate due to the following reason. 
It is known from \cite{Scho84} that the mass of a (given) manifold plays a role in the expansion 
for  Sobolev quotients of highly concentrated functions. In our case this must be 
done \underline{uniformly in $s$}, and the problem could be that the {\em principal term} 
coming from the mass could become negligible as $s \to 0$. To solve this issue 
we exploit a symmetry  $s \to -s$ for Rossi spheres, discussed in Section \ref{s:CRcoord}, 
which implies that all variational expansions are indeed \underline{even} in $s$ and hence 
the mass, which vanishes with $s$, gives still a dominant sign to the asymptotic expansion 
of the CR-Sobolev quotient. Two appendices are devoted to the estimates 
of the latter quantity in two different scaling regimes.

\

\noindent \textbf{Acknowledgements}  J.-H.C. (P.Y., resp.) are grateful to
Scuola Normale Superiore and Princeton University (Academia Sinica in Taiwan, resp.) for the
kind hospitality.  J.-H.C. is supported by the project MOST 107-2115-M-001-011 of Ministry of Science and Technology and NCTS of Taiwan.
A.M. is supported by the project {\em Geometric Variational Problems} from Scuola Normale Superiore and 
by MIUR Bando PRIN 2015 2015KB9WPT$_{001}$.  He is also a member of GNAMPA as part of INdAM.  He also would like to thank
Academia Sinica in Taiwan and Princeton University for the kind hospitality. P. Y. acknowledges support from the NSF for the grant DMS 1509505.

\

\section{Background material}\label{s:CRcoord}

In this section we recall some useful facts about CR manifolds and the 
properties of CR normal coordinates,  constructed in \cite{JLCR}. 
We then describe some general features of Rossi spheres.

\medskip

\subsection{Preliminary facts on CR manifolds}
\label{ss:notprel}

Let us begin by recalling  the following commutation relations on tensors, see Lemma 2.3 in \cite{Lee2} 
(we also refer to this paper for our tensorial notation)
\begin{equation}  \label{eq:comm}
\left\{ 
\begin{array}{ll}
c,_{1 \overline{1}} - c,_{\overline{1} 1} = i c,_0 + k c R; &  \\ 
c,_{01} - c,_{10} = c,_{\overline{1}} A_{11} - k c A_{11},_{\overline{1}}; & 
\\ 
c,_{0\overline{1}} - c,_{\overline{1} 0} = c,_1 A_{\overline{1} \overline{1}%
} + k c A_{\overline{1} \overline{1}},_1. & 
\end{array}
\right.
\end{equation}
Here $c$ is a tensor with $1$ or $\bar{1}$ as sub-indices, $k$ is the number
of $1$-sub-indices of $c$ minus the number of $\overline{1}$-sub-indices of $c$
and where, we recall, we are assuming that $h_{1\overline{1}}=1$ (so $A_{%
\bar{1}\bar{1}}$ $=$ $A_{\bar{1}}^{1}$ and $A_{11}$ is the complex conjugate
of $A_{\bar{1}\bar{1}}).$

\ 

\noindent In the system of coordinates we will describe below, for $(z,t) \in \H^1$ near zero we will set 
\begin{equation}  \label{eq:rho}
\rho^4 = |z|^4 + t^2.
\end{equation}
For $k \in \mathbb{Z}$ we denote by $\tilde{O}(\rho^k)$ a function $f(z,%
\overline{z},t)$ for which $|f| \leq C \rho^k$ for some $C > 0$; we use
instead the symbol $\tilde{O}^{\prime}(\rho^k)$ for a function $f(z,\overline{z}
,t)$ such that 
\begin{equation*}
|f| \leq C \rho^k, \qquad |\partial_z f| \leq C \rho^{k-1} \left| \partial_z
\rho \right|, \qquad |\partial_{\overline{z}} f| \leq C \rho^{k-1} \left|
\partial_{\overline{z}} \rho \right|, \qquad |\partial_t f| \leq C
\rho^{k-1} \left| \partial_t \rho \right|.
\end{equation*}
One can define similarly the symbols $\tilde{O}^{\prime \prime}(\rho^ k)$, $\tilde{O}^{\prime \prime \prime }(\rho^k)$, etc. We will use  $O(\rho^k)$ for a
function which is of the form $\tilde{O}^{(j)}(\rho^k)$ for every integer $j$, 
or for $j$ large enough for our purposes.

Large positive constants are always denoted by $C$, and the value of $C$ is
allowed to vary from one formula to another and also within the same line. When
we want to stress the dependence of the constants on some parameter (or
parameters), we add subscripts to $C$, as $C_\d$, etc.. Also constants with
this kind of subscripts are allowed to vary.

\

Let us recall the notions of \emph{pseudo-hermitian geometry} from 
\cite{Web83} and \cite{Lee1}. 
We would need the following result in \cite{JLCR} on page 313, Proposition 2.5.
For a differential form $\eta ,$ let us denote by $\eta _{(m)}$ the part of
its Taylor series that is homogeneous of degree $m$ in terms of 
parabolic dilations (see \cite{JLCR} for more details). 

\begin{pro}\label{p:phcoord} Let ${\tilde{Z}}_1$ be a special frame dual to
${\tilde{\th}}^1$ (with ${\tilde{h}}_{1 \ou}=2$) and let ${\th^1}={\sqrt{2}}
{\tilde{\th}}^1$ be a unitary coframe ($h_{1 \ou}=1$). Then in
pseudo-hermitian normal coordinates $(z, t)$ with respect to ${\tilde{Z}}_1$,
${\tilde{\th}}^1$, we have

\noindent $(a)$ $\th_{(2)} = \overset{\circ}{\theta}; \qquad \th_{(3)} = 0; \qquad
\th_{(m)} =
 \frac 1 m {\sqrt 2} \left( i  z \th^{\ou} - i \ov{z} \th^1 \right)_{(m)},
\quad m \geq 4$;

\noindent $(b)$ $\th^1_{(1)} = {\sqrt 2}d z; \quad \th^1_{(2)} = 0; \qquad
  \th^1_{(m)} = \frac 1 m \left( {\sqrt 2}z \o^1_1 + 2 t A_{\ou \ou} \th^{\ou}
  - {\sqrt 2}\ov{z} A_{\ou \ou} \th  \right)_{(m)}, \quad m \geq 3$

\noindent $(c)$  $(\o^1_1)_{(1)} = 0; \quad (\o^1_1)_{(m)} = \frac 1
m \left( {\sqrt 2}R ( z \th^{\ou} - \ov{z} \th^1 )  + A_{11,\ou} ( {\sqrt 2}z \th
- 2 t \th^1) - A_{\ou \ou, 1} ( {\sqrt 2}\ov{z} \th - 2 t \th^{\ou})
\right)_{(m)}$,

\noindent $m \geq 2$, where $\overset{\circ}{\theta} = dt + i z d\ov{z} - i \ov{z} dz$.
\end{pro}

\medskip

\begin{df} Given a three dimensional pseudo-hermitian manifold
$(M,\th)$ we define a real symmetric tensor $Q$ as
$$
  Q = Q_{jk} \th^j \otimes \th^k, \qquad j,k \in \left\{ 0, 1, \ou \right\}
$$
with $\th^0:=\th$, whose components with respect to any admissible coframe are given by
$$
  Q_{11} = \overline{Q_{\ou \ou}} = 3 i A_{11}; \qquad \quad Q_{1 \ou} = Q_{\ou 1}
  = h_{1 \ou}R;
$$
$$
  Q_{01} = Q_{10} = \overline{Q_{0 \ou}} = \overline{Q_{\ou 0}} = 4 A_{11,}^{\;\;\;\;\;1}
  + i R_{,1}; \qquad \quad Q_{00} = 16 {\text Im}  A_{11,}^{\;\; \;\;\; 11} - 2 \D_b R.
$$
\end{df}

\ 

\noindent We have then the following result, see page 315 in \cite{JLCR},
Theorem 3.1. 

\begin{pro}\label{p:CRnormcoord} Suppose $M$ is a strictly pseudo-convex pseudo-hermitian manifold of
dimension three, and let $q \in M$. Then for any integer $N \geq 2$ there exists a choice of
contact form $\th$ such that all symmetrized covariant derivatives
of $Q$ with total order less or equal than $N$ vanish at $q$, that
is
\begin{equation}\label{eq:sqv}
    Q_{\langle jk, l \rangle} = 0 \qquad \qquad \hbox{at q if }
   \quad \mathbb{O}(jkl) \leq N.
\end{equation}
\end{pro}

By \emph{CR normal coordinates} of order $N$ we mean the pseudo-hermitian normal coordinates
with $\th $ chosen as in Proposition \ref{p:CRnormcoord}.
We recall (\cite{JLCR}) that for a multi index $l = (l_1, \dots, l_s)$ we count its order as
$$
  \mathbb{O}(l) = \mathbb{O}(l_1) + \cdots + \mathbb{O}(l_s),
$$
where $\mathbb{O}(1) = \mathbb{O}(\ov{1}) = 1$ and where $\mathbb{O}(0) = 2$. 
The symmetrized covariant derivatives are defined by
$$
  Q_{\langle l \rangle} = \frac{1}{s!} \sum_{\s \in \mathbb{S}_s} Q_{\s l}; \qquad
  \s l = \left( l_{\s(1)}, \dots, l_{\s(s)} \right).
$$
In \cite{CMY17}, Proposition A.5, the following result was proved.

\begin{pro}\label{p:CRcoord} In CR normal coordinates of order $N = 4$, we have a
contact form ${\th}$ such that 
$$
  {\th} = \left( 1 + O(\rho^4) \right) \overset{\circ}{\theta} + O(\rho^5) dz
  + O(\rho^5) d\ov{z}; \qquad {\th}^1 = \left( 1 + O(\rho^4)
  \right) {\sqrt 2}dz + O(\rho^4) d\ov{z} + O(\rho^3) \overset{\circ}{\theta};
$$
$$
  {\o}^1_1 = O(\rho^3) d z + O(\rho^3) d\ov{z} + O(\rho^2) \overset{\circ}{\theta};
$$
$$
  {Z}_1 = \left( 1 + O(\rho^4) \right) \overset{\circ}{Z}_{1} + O(\rho^4)
  \overset{\circ}{Z}_{\overline{1}} + O(\rho^5) \frac{\pa}{\pa t}; \qquad \quad {T}
  = \left( 1 + O(\rho^4) \right) \frac{\pa}{\pa t} +
  O(\rho^3) \overset{\circ}{Z}_{1} + O(\rho^3) \overset{\circ}{Z}_{\overline{1}},
$$
where we recall
\begin{equation}\label{eq:circleformulas}
\overset{\circ}{\theta} = dt + i z d\ov{z} - i \ov{z} dz; \qquad \overset{\circ}{Z}_{1} =
 \frac{1}{{\sqrt 2}} \left( \frac{\pa}{\pa z} + i \ov{z} \frac{\pa}{\pa t}
 \right); \qquad \rho^4 = t^2 + |z|^4.
\end{equation}
\end{pro}

\medskip

\subsection{Rossi spheres}

We recall here some properties of {\em Rossi spheres}, introduced in \cite{Rossi} as a 
non-embeddable example of CR manifold (see also \cite{Burns}). These are families of 
CR structures on $S^3$, containing the standard one, obtained in the following way. 

Considering the complex vector field $Z_1$ as in \eqref{B-1} and its conjugate $Z_{\bar 1}$, 
one defines the CR structure $J_{(s)}$ by setting  $J_{(s)}Z_{1(s)}=iZ_{1(s)}$, where 
%
%
%
%
%
%
%
%
%
\begin{equation}
Z_{1(s)} =Z_{1}+\frac{s}{\sqrt{1+s^{2}}}Z_{\bar{1}},  \label{B6-1} \qquad \quad
Z_{\bar{1}(s)} =Z_{\bar{1}}+\frac{s}{\sqrt{1+s^{2}}}Z_{1}.  \notag
\end{equation}
Corresponding to these vector fields, we have the dual forms 
\begin{equation*}
\theta _{(s)}^{1} =(1+s^{2})\theta ^{1}-s\sqrt{1+s^{2}}\theta ^{\bar{1}},
\qquad \quad 
\theta _{(s)}^{\bar{1}} =(1+s^{2})\theta ^{\bar{1}}-s\sqrt{1+s^{2}}\theta
^{1}, 
\end{equation*}
where $\th^1 = z^2 dz^1 - z^1 dz^2$.  Compute%
\begin{equation}
i\theta _{(s)}^{1}\wedge \theta _{(s)}^{\bar{1}} =(1+s^{2})i\theta
^{1}\wedge \theta ^{\bar{1}}  \label{B7} 
=(1+s^{2})d\hat{\theta},  \notag
\end{equation}
where $d\hat{\theta} $ $=$ $i\theta ^{1}\wedge \theta ^{\bar{1}},$ i.e., 
$h_{1\bar{1}}$ $=$ $1.$ Hence, from (\ref{B7}) we get%
\begin{equation*}
h_{1\bar{1}}^{(s)}=\frac{1}{1+s^{2}}\text{ and }h_{(s)}^{1\bar{1}}:=(h_{1%
\bar{1}}^{(s)})^{-1}=1+s^{2}.
\end{equation*}

\noindent By taking 
\begin{equation*}
\tilde{\theta}_{(s)}^{1}=\frac{1}{\sqrt{2(1+s^{2})}}\theta _{(s)}^{1},
\end{equation*}%
\noindent we have $\tilde{h}_{1\bar{1}}^{(s)}$ $=$ $2.$  The Webster
curvature $R$ of $(J,\hat{\theta} )$ is identically equal to $2$. Then we
should take $\omega _{1}^{1}$ $=$ $-2 i \hat{\theta} $ in the structure
equation \eqref{eq:structure}, such that $d\omega _{1}^{1}$ $=$ $2 \theta ^{1}\wedge \theta ^{%
\bar{1}}.$ 
We can then determine, from the structure
equation for $(J_{(s)},\hat{\theta} ),$ that
\begin{equation}
\omega _{1(s)}^{1} = -2 i (1+2s^{2})  \hat{\theta},  \label{BK} \qquad
h_{(s)}^{1\bar{1}}A_{\bar{1}\bar{1}(s)} =4is\sqrt{1+s^{2}},   \qquad
R_{(s)} =2 (1+2s^{2}).  \notag
\end{equation}
%
Dual to $\theta ^{1}$ $=$ $z^{2}dz^{1}$ $-$ $%
z^{1}dz^{2},$ we have%
\begin{equation}
Z_{1} =Z_{1}^{S^{3}}
=z^{\bar{2}}\frac{\partial }{\partial z^{1}}-z^{\bar{1}}\frac{\partial }{%
\partial z^{2}}.
\end{equation}

The sub-Laplacian associated to $(J_{(s)},\hat{\theta})$ reads%
\begin{equation}
\triangle _{b}^{(s)} =h_{(s)}^{1\bar{1}}(Z_{1(s)}Z_{\bar{1}(s)}+Z_{\bar{1}%
(s)}Z_{1(s)}) 
=(1+2s^{2})\triangle _{b}^{(0)}+2s\sqrt{1+s^{2}}(Z_{1}^{2}+Z_{\bar{1}%
}^{2}).
\end{equation}

\noindent It follows that,  at $s=0$, the first-and second-order 
derivatives of $\triangle _{b}^{(s)}$ w.r.t. $s$ are given by 
\begin{equation}\label{eq:dot-L-ddot-L}
- \dot{\D}_{b}=2Z_{\overline{1}}Z_{\overline{1}} +\hbox{conj.}; 
\qquad \qquad - \ddot{\Delta}_b =  - 4 \Delta_b.
\end{equation}
Moreover since $R_s = 2(1 + 2 s^2)$ it follows that, still at $s = 0$
\begin{equation}\label{eq:dot-R-ddot-R}
  \dot{R} = 0; \qquad \qquad \ddot{R} = 8. 
\end{equation}



We next analyze a symmetry property of Rossi spheres, that will 
imply in particular the symmetry of the mass in $s$. Consider the
 diffeomorphism $\iota :$ $S^{3}$ $\rightarrow $ $S^{3}$ defined by%
 \begin{equation}\label{eq:iota}
 \iota (z^{1},z^{2})=(iz^{1},z^{2}),
 \end{equation}
 
 \noindent which fixes the point (0,1). A direct computation shows that $\iota _{\ast
 }Z_{1}^{S^{3}}$ $=$ $iZ_{1}^{S^{3}}$ and hence $\iota _{\ast }Z_{\bar{1}%
 }^{S^{3}}$ $=$ ($-i)Z_{\bar{1}}^{S^{3}}.$ By (\ref{B6-1}), we compute%
 \begin{eqnarray}
 \iota _{\ast }Z_{1(s)} =\iota _{\ast }Z_{1}+\frac{s}{\sqrt{1+s^{2}}}\iota
 _{\ast }Z_{\bar{1}}  \label{B8} 
 =iZ_{1}+\frac{s}{\sqrt{1+s^{2}}}(-i)Z_{\bar{1}}  
 =iZ_{1(-s)}.  
 \end{eqnarray}
 
 \noindent It follows that 
 \begin{eqnarray*}
 (\iota ^{\ast }J_{(-s)})Z_{1(s)} &=&\iota _{\ast }^{-1}J_{(-s)}(\iota _{\ast
 }Z_{1(s)}) 
 =\iota _{\ast }^{-1}J_{(-s)}(iZ_{1(-s)})\text{ (by (\ref{B8}))} 
 =\iota _{\ast }^{-1}(-Z_{1(-s)}) \\
 &=&(-1)(-i)Z_{1(s)}\text{ (by the inverse of }(\ref{B8})) 
 =iZ_{1(s)}=J_{(s)}Z_{1(s)}.
 \end{eqnarray*}
 
 \noindent Hence we have shown 
 \begin{equation}
 J_{(s)}=\iota ^{\ast }J_{(-s)}.  \label{S1}
 \end{equation}
 
 \noindent Let $v_{(s)}$ denote the conformal factor in $\check{\theta}_{(s)}$
 $=$ $e^{2v_{(s)}}\hat{\theta}$, yielding CR normal coordinates with respect to $J_{(s)}$. It then follows that  
 \begin{eqnarray}
 v_{(s)} =\iota ^{\ast }v_{(-s)},  \label{S2} \qquad \quad 
 \check{\theta}_{(s)} =\iota ^{\ast }\check{\theta}_{(-s)},  \notag
 \end{eqnarray}%
 \noindent and hence $\check{G}_{s}=\iota ^{\ast }\check{G}_{-s}$ by observing%
 \begin{equation}
 \iota ^{\ast }\hat{\theta}=\hat{\theta}.  \label{S3}
 \end{equation}%
 \noindent Write 
 \begin{equation*}
 \check{G}_{s}= 2 \rho _{s}^{-2}+A_{s}+O(\rho _{s})
 \end{equation*}%
 \noindent in $s$-CR normal coordinates near $(0,1)$. Then $\rho _{s}$ $=$ $%
 \iota ^{\ast }\rho _{-s}$ $=$ $\rho _{-s}\circ \iota $ and%
 \begin{equation*}
 A_{s}=\iota ^{\ast }A_{-s}=A_{-s}\circ \iota =A_{-s}
 \end{equation*}
 
 \noindent near the point $(0,1)$. So, we have obtained
 \begin{eqnarray*}
 m(J_{(s)},\theta _{(s)}) =12\pi A_{s} 
 =12 \pi  A_{-s} 
 =m(J_{(-s)},\theta _{(-s)}),
 \end{eqnarray*}
 where $\theta_{(s)} = \check{G}_{(s)}^2 \check{\theta}_{(s)}$.  
This property (and other related ones) will be crucial in the last section of the paper.  

\

\section{CR normal coordinates on Rossi spheres}\label{s:CR-Rossi}

In this section we will find the main-order terms of CR normal coordinates on Rossi spheres. 
We first determine the principal term in the required conformal factor, then discuss 
pseudo-hermitian coordinates and finally CR normal coordinates.  This will allow us 
to express with a good precision the Green's function of the conformal sub-Laplacian 
near its pole.

\medskip 

\subsection{Conformal factor in normalized contact form on Rossi spheres.}

Fix $p = (0,1) \in S^3 \subseteq \C^2$ and consider a contact form 
$\check{\theta}'_{(s)}$ $=$ $e^{2v_{(s)}}\hat{\theta}'$, 
where $\hat{\theta}' = 2 \hat{\th} = i(\bar{\partial}-\partial )(|z^{1}|^{2}+|z^{2}|^{2})$ 
yielding CR normal coordinates (see Proposition \ref{p:CRnormcoord}) with 
respect to $J_{(s)}$ for $N = 4$.  We are
going to solve an equation for $v_{(s)}$ as in Lemma 3.11 of Jerison-Lee's
paper (\cite{JLCR}). Write%
\begin{equation}\label{eq:v-vi}
v_{(s)}=v_{2}+v_{3}+... , 
\end{equation}
where $v_{2}$ $\in $ $\mathcal{R}_{2}$ $\subset $ $\mathcal{P}%
_{2}, $ $v_{3}$ $\in $ $\mathcal{P}_{3}$.  Recall that, in the notation 
of \cite{JLCR}, $\mathcal{P}_m$ denotes the vector space of polynomials 
in $(z,t)$ that are homogeneous of degree $m$ in terms of parabolic 
dilations (for which $t$ has homogeneity 2), and $\mathcal{Y}_m \subseteq \mathcal{P}_m$ 
denotes the subspace of polynomials independent of $t$.

 First, write $v_{2}$ $\in $ $%
\mathcal{R}_{2}$ as $v_{2}$ $=$ $az^{2}+bz\bar{z}+c\bar{z}^{2}$ (($z,t)$
being pseudo-hermitian normal coordinates for $\hat{\theta}'$ at $p)$
satisfying%
\begin{equation}
L_{2}v_{2} =-z^{2}Q_{11}-\bar{z}^{2}Q_{\bar{1}\bar{1}}-z\bar{z}Q_{1\bar{1}%
}-\bar{z}zQ_{\bar{1}1}; \label{L2} \qquad \quad
L_{2} =-2|z|^{2}(\partial _{z}\partial _{\bar{z}}+\partial _{\bar{z}%
}\partial _{z})-12,   
\end{equation}
where $Q_{11}=3iA_{11(s)}^{JL}=Q_{\bar{1}\bar{1}}$ and $Q_{1\bar{1}}
=R_{1\bar{1}(s)}^{JL}=Q_{\bar{1}1}$ are w.r.t. the Jerison-Lee coframe
$\theta _{JL}^{1}=\theta
_{(s)}^{1}/\sqrt{1+s^{2}}$ with $h_{1\bar{1}(s)}^{JL}=2$ w.r.t. $\hat{\theta}'$
by the formulas for $Q_{jk}$ on
page 315 in \cite{JLCR} and (\ref{BK}). We compute

\begin{equation}
\tilde{Q}_{11} =3iA_{11(s)}=\frac{12s}{\sqrt{1+s^{2}}}=\tilde{Q}_{\bar{%
1}\bar{1}},  \label{Q11}; \qquad \quad
\tilde{Q}_{1\bar{1}} =R_{1\bar{1}(s)}=h_{1\bar{1}}^{(s)}R_{(s)}=2\frac{%
1+2s^{2}}{1+s^{2}}, 
\end{equation}
\noindent with respect to the  co-frame $\theta _{(s)}^{1}$. A direct computation shows that%
\begin{equation*}
L_{2}v_{2}=-12az^{2}-12c\bar{z}^{2}-16b|z|^{2}, 
\end{equation*}
where $Q_{11}=12a,$ $Q_{\bar{1}\bar{1}}=12c$, $Q_{1\bar{1}}=Q_{%
\bar{1}1}=8b$ and 
\begin{equation}
a =c=s\sqrt{1+s^{2}},  \label{abc} \qquad \qquad
b =\frac{1}{4}(1+2s^{2}).  
\end{equation}
For $v_{3},$ we observe that all $Q_{jk,l}$'s for $j,k,l$ being $1$ or $%
\bar{1}$ vanish since the space derivatives of the constant $R_{0}$ is zero. On
the other hand, $Q_{0k}$ and $Q_{k0}$ for $k$ $=$ $1$ or $\bar{1}$ also
vanish since they involve space derivatives by formulas on page 315 in \cite{JLCR}. Altogether, the right hand side of the equation in Lemma 3.11 in \cite{JLCR} for $m=3$
equals zero, so we have%
\begin{equation*}
L_{3}v_{3}=0.
\end{equation*}
By Lemma 3.9 in \cite{JLCR}, we learn that $L_{3}$ is invertible
on $\mathcal{P}_{3}.$ It follows that 
\begin{equation}
v_{3}=0.  \label{u3}
\end{equation}
Therefore,  from (\ref{abc}) and (\ref{u3}) we get the following result. 

\begin{lem}\label{l:conf-fact-ph}
In pseudo-hermitian coordinates, the conformal factor expands in homogeneous powers as 
\begin{equation}
v_{(s)}=s\sqrt{1+s^{2}}(z^{2}+\bar{z}^{2})+\frac{1}{4}%
(1+2s^{2})|z|^{2}+v_{4}+....  \label{us}
\end{equation}
\end{lem}

\medskip

\subsection{Pseudo-hermitian coordinates on Rossi spheres}

Recall that on Rossi spheres we have 
$$
  \th^1_{(s)} = (1+s^2) \th^1 - s \sqrt{1+s^2} \th^{\ou}; 
  \qquad \qquad \o^1_{1 (s)} = - i (1 + 2 s^2) 2 \th, 
$$
and that pseudo-hermitian coordinates near $(0,1)$ are defined by the equation 
\begin{equation}\label{eq:psi-h}
  \n_{\dot{\s}} \dot{\s} = 2 \, c \, \hat{T}'; \qquad \quad \s(0) = (0,1), 
\end{equation}
where $\hat{T}'$ is the unique vector field such that $\hat{\th}'(\hat{T}') = 1$ and 
$d \hat{\th}'(\hat{T}', \cdot) = 0$.
Recall also that  
$$
  \hat{\th}' = i \sum_{i=1}^2 \left( z^i d \ov{z}^i - \ov{z}^i d z^i  \right); 
  \qquad \qquad \hat{T}' = -  {\text Im}  \left( z^1 \frac{\pa}{\pa z_1} 
  + z^2 \frac{\pa}{\pa z^2} \right) = \frac 12 i 
  \sum_{i=1}^2 \left( z^i \frac{\pa}{\pa z_i} - \ov{z}^i \frac{\pa}{\pa \ov{z}_i} \right). 
$$
Setting 
$$
  \dot{\s} = \a {Z}^{JL}_{1 (s)}  + \b {Z}^{JL}_{\ou (s)}  + \g \hat{T}', 
$$
\eqref{eq:psi-h} becomes 
\begin{eqnarray} \nonumber
 2 c \hat{T}' & = &  \n_{\dot{\s}} \dot{\s}  = \left( \dot{\a} + \a  \o^1_{1 (s)} (\dot{\s})  \right) {Z}^{JL}_{1 (s)} 
 + \left( \dot{\b} + \b  \o^{\ou}_{\ou (s)} (\dot{\s})  \right) {Z}^{JL}_{\ou (s)}  + \dot{\g} \hat{T}' \\ & = & 
 \left( \dot{\a} - i \a (1+2s^2) \g  \right)  {Z}^{JL}_{1 (s)}   + \left( \dot{\b} + i \b (1+2s^2) 
 \g \right) {Z}^{JL}_{\ou (s)}  + \dot{\g} \hat{T}'. 
\end{eqnarray}
If $\tau$ parametrizes the  curve $\s$, the above formulas imply that 
$$
  \g= 2 c \t; \qquad \qquad \frac{\dot{\a}}{\a} = i (1+2s^2) \g; \qquad \qquad
  \frac{\dot{\b}}{\b} = - i (1+2s^2) \g, 
$$
which in turn yields 
$$
  \a(t) = \a(0) e^{i c  (1+2s^2) \t^2}; \qquad \qquad \b(t) = \b(0) 
  e^{- i c  (1+2s^2) \t^2}. 
$$
Therefore we obtained  
$$
  \dot{\s} = \a(0) e^{ i c (1+2s^2)  \t^2 } {Z}^{JL}_{1 (s)}  
  + \b(0) e^{-  i c (1+2s^2)  \t^2 } {Z}^{JL}_{\ou (s)}  + 2 c \t \hat{T}. 
$$
Recall also that ${Z}_{1 (s)}  = {Z}_{1 (0)} + 
\frac{s}{\sqrt{1+s^2}} {Z}_{\ou (0)}$. Hence we need to 
solve for 
\begin{eqnarray*}
\dot{z}_1(\t) & = & dz_1(\dot{\s}(\t)) = \a(0) e^{ i c (1+2s^2)  \t^2 } dz_1 ({Z}^{JL}_{1 (s)} ) 
  + \b(0) e^{-  i c (1+2s^2)  \t^2 } dz_1 ({Z}^{JL}_{\ou (s)})  + 2 c \t dz_1 (\hat{T}') \\ & = & 
  \sqrt{1+s^2} \left( \a(0) e^{i \d \t^2} \ov{z}_2(\t) + \b(0) e^{- i \d \t^2} \frac{s}{\sqrt{1 + s^2}} \ov{z}_2(\t) 
  \right) +  c \t i z_1(\t),  
\end{eqnarray*}
where $\d =  c (1+2s^2) $. Similarly, we obtain 
$$
  \dot{z}_2(\t) = 
  \sqrt{1+s^2} \left( -   \a(0) e^{i \d \t^2} \ov{z}_1(\t) - \b(0) e^{- i \d \t^2} \frac{s}{\sqrt{1 + s^2}} \ov{z}_1(\t) 
  \right) +    c \t i z_2(\t).   
$$
Once we will solve for this system, the pseudo-hermitian coordinates will be given by the map 
\begin{equation}\label{eq:coord-t}
  (z, \ov{z}, t) = (\a (0) \t, \b(0) \t, c \t^2) \quad \longmapsto \quad (0,1) + \int_0^t \dot{\s}(\eta) d \eta. 
\end{equation}
Setting for simplicity 
$$
  A_0 = \a(0); \qquad \qquad B_0 = \b(0) \frac{s}{\sqrt{1 + s^2}}; \qquad \qquad C_0 = 2 c; 
  $$
  $$
   F_0(\t) := \sqrt{1+s^2} (A_0 e^{i \d \t^2} + B_0 e^{- i \d \t^2}) = f_0(\t) + i g_0(\t), 
$$
we have then the system of ODEs
$$
\dot{z}_1(\t) = F_0(\t) \ov{z}_2(\t) + i C_0 \t z_1(\t); \qquad \qquad 
\dot{z}_2(\t) = - F_0(\t) \ov{z}_1(\t) + i C_0 \t z_2(\t), 
$$
which in real form becomes 
$$
  \begin{cases}
  \dot{x}_1(\t) = f_0(\t) x_2(\t) + g_0(\t) y_2(\t) - C_0 \t y_1(\t); \\
  \dot{y}_1(\t) = g_0(\t) x_2(\t) - f_0(\t) y_2(\t) + C_0 \t x_1(\t); \\ 
  \dot{x}_2(\t) = -f_0(\t) x_1(\t) - g_0(\t) y_1(\t) - C_0 \t y_2(\t); \\
  \dot{y}_2(\t) = f_0(\t) y_1(\t) - g_0(\t) x_1(\t) + C_0 \t x_2(\t). 
  \end{cases}
$$
We rewrite this system as 
$$
 \dot{\mathfrak{X}}(\t) = \mathfrak{A}(\t) \mathfrak{X}(\t), 
$$
where 
$$
 \mathfrak{A}(\t) = \left(  \begin{matrix}
  0 & - C_0 \t & f_0(\t) & g_0(\t) \\C_0 \t & 0 & g_0(\t) & - f_0(\t) \\
  - f_0(\t) & - g_0(\t) & 0 & - C_0 \t \\ 
  - g_0(\t) & f_0(\t) & C_0 \t & 0
  \end{matrix} \right).
$$
We can Taylor-expand the solution to an arbitrary order in $\t$. Differentiating the 
above ODE we obtain 
$$
 \ddot{\mathfrak{X}}(\t) = \dot{\mathfrak{A}}(\t) \mathfrak{X}(\t) +   \mathfrak{A}(\t)^2 \mathfrak{X}(\t); 
 \qquad \quad \dddot{\mathfrak{X}}(\t) = \ddot{\mathfrak{A}}(\t) \mathfrak{X}(\t) 
 + (\mathfrak{A}(\t)  \dot{\mathfrak{A}}(\t) + 2 \dot{\mathfrak{A}}(\t) \mathfrak{A}(\t)) \mathfrak{X}(\t) + \mathfrak{A}(\t)^3 \mathfrak{X}(\t). 
$$
We have that 
$$
\mathfrak{A}(0) = \left(
\begin{array}{cccc}
 0 & 0 & \text{Re} A_0+\text{Re} B_0 & \text{Im} A_0+\text{Im} B_0 \\
 0 & 0 & \text{Im} A_0+\text{Im} B_0 & -\text{Re} A_0-\text{Re} B_0 \\
 -\text{Re} A_0-\text{Re} B_0 & -\text{Im} A_0-\text{Im} B_0 & 0 & 0 \\
 -\text{Im} A_0-\text{Im} B_0 & \text{Re} A_0+\text{Re} B_0 & 0 & 0 \\
\end{array}
\right)$$
$$
  \dot{\mathfrak{A}}(0) = \left(  \begin{matrix}
        0 & -C_0  & 0 & 0 \\C_0 & 0 & 0 & 0 \\
       0& 0 & 0 & -C_0 \\ 
        0 & 0 & C_0 & 0
        \end{matrix} \right); 
$$
$$
 \ddot{\mathfrak{A}}(0) = \left(
 \begin{array}{cccc}
  0 & 0 & 2 d (\text{Im} B_0-\text{Im} A_0) & -2 d (\text{Re} B_0-\text{Re} A_0) \\
  0 & 0 & -2 d (\text{Re} B_0-\text{Re} A_0) & 2 d (\text{Im} A_0-\text{Im} B_0) \\
  2 d (\text{Im} A_0-\text{Im} B_0) & 2 d (\text{Re} B_0-\text{Re} A_0) & 0 & 0 \\
  2 d (\text{Re} B_0-\text{Re} A_0) & 2 d (\text{Im} B_0-\text{Im} A_0) & 0 & 0 \\
 \end{array}
 \right).
$$
In conclusion, looking at the first three terms in the Taylor expansion of $\mathfrak{X}(\t)$ near $(0,1)$ we find that 
$$
  \mathfrak{X}(\t) = 
  \left( \begin{matrix}
  \tilde{\t} (\text{Re} A_0+\text{Re} B_0) \\\tilde{\t} (\text{Im} A_0+\text{Im} B_0) \\ 
  1-\frac{1}{2} \tilde{\t}^2 \left((\text{Im} A_0+\text{Im} B_0)^2+(\text{Re} A_0+\text{Re} B_0)^2\right) \\ \frac{C_0 \tilde{\t}^2}{2} 
  \frac{1}{2(1+s^2)}
  \end{matrix} \right) 
  + o(\t^2); \qquad \tilde{\t} = \sqrt{1+s^2} \, \t. 
$$
Recalling \eqref{eq:coord-t}, we then obtain the following result.

\begin{lem}
Pseudo-hermitian coordinates near $(0,1)$ on Rossi spheres w.r.t. $\hat{\th}' = 2 \hat{\th}$ are given by the following map 
$$
 (z, \ov{z}, t) \longmapsto \left(  \begin{matrix}
        \sqrt{ (1+s^2)} \left( z + \frac{s}{\sqrt{1+s^2}} \ov{z} \right) \\ 
            1 - \frac 12  (1+s^2) \left| z + \frac{s}{\sqrt{1+s^2}} \ov{z}  \right|^2 + i \frac{t}{2}
            \end{matrix} \right) + o(\rho^2).
$$
Inverting in the first component, we have in particular that 
\begin{equation}\label{eq:z-from-z1}
z = \frac{1}{\sqrt{1+s^2}}(1+s^2) \left( z_1 - \frac{s}{\sqrt{1+s^2}} \ov{z}_1 \right) + o(\rho^2).
\end{equation}
\end{lem}

\medskip

\subsection{CR normal coordinates}

Recalling \eqref{eq:v-vi}, Lemma \ref{l:conf-fact-ph} and using \eqref{eq:z-from-z1}, we get 
%
\begin{eqnarray}\nonumber
v_2  = \frac{1}{4}  \left(3 (z_1^2 + \ov{z}_1^2) s \sqrt{s^2+1} \left(2 s^2+1\right)- |z_1|^2 \left(12 s^4+12 s^2-1\right) \right) 
  = A_1 (z_1^2 + \ov{z}_1^2 ) + B_1 |z_1|^2, 
\end{eqnarray}
where 
\begin{equation}\label{eq:A1B1}
  A_1 = \frac{1}{8} 2 \left(s^2+1\right)^{1/2} 3 s \left(2 s^2+1\right); 
  \qquad \qquad B_1 = - \frac{1}{8} 2  \left(12 s^4+12 s^2-1\right). 
\end{equation}
Recall that also 
$$
  \hat{\th}^1_{(0)} = z^2 d z^1 - z^1 dz^2; \qquad  \qquad 
  \hat{\th}' = i \sum_{i=1}^2 \left( z^i d \ov{z}^i - \ov{z}^i d z^i  \right), 
$$
and that 
$$
  \hat{\th}^1_{(s)} = (1+s^2) \hat{\th}^1_{(0)}  - s \sqrt{1+s^2} \, \hat{\th}^{\ou}_{(0)}.
$$
Conformally changing the contact form and recalling Appendix 1.1.1 in \cite{CMY17}, we have that 
$\hat{\th}^1_{(s)}$ transforms as  
\begin{equation}\label{eq:th-hat-1-}
  \hat{\th}^1_{(s)} \quad \longmapsto \quad e^v (\hat{\th}^1_{(s)} + 2 i v^1 \hat{\th}'), 
\end{equation}
where  ($d \hat{\th}' = 2 i \hat{\th}^1_{(0)}  \wedge \hat{\th}^{\ou}_{(0)}$)
$$
  v^1  = h^{1 \ou}_{(s)} \hat{Z}_{\ou (s)} v = \frac{(1+s^2)}{2} \hat{Z}_{\ou (s)} v. 
$$
By computing explicitly, it turns out that 
\begin{eqnarray} \nonumber
(v_2)^1 = \left( \frac{\left(s^2+1\right)}{2} \frac{\ov{z}_2 s (2 {A_1} z_1 +{B_1} \ov{z}_1)}{\sqrt{s^2+1}}+\frac{\left(s^2+1\right)}{2}
   z_2 (2 {A_1} \ov{z}_1+{B_1} {z_1})\right), 
\end{eqnarray}
which can be written as 
\begin{equation}\label{eq:new-33}
(v_2)^1 =
A_2 z_1 z_2 + B_2 z_1 \ov{z}_2 + C_2 \ov{z}_1 z_2 + D_2 \ov{z}_1 \ov{z}_2, 
\end{equation}
with
\begin{equation}\label{eq:ABCD2}
  A_2 = \frac 12 \left(s^2+1\right)  {B_1}, 
  \qquad \quad
  B_2 =    s {A_1} \sqrt{s^2+1}, 
 \qquad \quad 
  C_2 =  {A_1} \left(s^2+1\right), 
  \qquad \quad D_2 =  \frac 12  s {B_1} \sqrt{s^2+1}.
\end{equation}
Up to higher order terms, we have that 
$$
  (v_2)^1 = (A_2 + B_2) z_1 + (C_2 + D_2) \ov{z}_1.
$$

Taylor expanding \eqref{eq:th-hat-1-}, up to higher-order terms $\hat{\th}^1_{(s)}$ 
transforms into 
%
$$
  \hat{\th}^1_{(s)} + \left[ v_2 \hat{\th}^1_{(s)} + 2 i (v_2)^1 \hat{\th}' \right]. 
$$ 
We now multiply by a complex unit factor $e^{i \psi}$, and impose a  
closeness condition on $e^{i \psi}$ multiplied by the latter form, up to higher-order terms, since by Proposition 
\ref{p:CRcoord} it should be approximately a constant multiple of $d \tilde{z}_{\text CR}$, up to h.o.t..  We then find 
\begin{eqnarray}\label{eq:unit-factor} \nonumber 
0 & = & d \left\{ e^{i \psi} \left( \hat{\th}^1_{(s)} + \left[ v_2 \hat{\th}^1_{(s)} + 2 i (v_2)^1 \hat{\th}' \right] \right) \right\} 
\\ & = & e^{i \psi} \left\{i d \psi \wedge \hat{\th}^1_{(s)} + i \, d \psi \wedge \left[ v_2 \hat{\th}^1_{(s)} + 2 i (v_2)^1 \hat{\th}' \right]  \right\} \\ 
& + & e^{i \psi} \left\{ d \hat{\th}^1_{(s)} + d v_2 \wedge \hat{\th}^1_{(s)} + v_2 d \hat{\th}^1_{(s)} 
+ 2 i d((v_2)^1) \wedge  \hat{\th}' + 2 i (v_2)^1 d \hat{\th}'  \right\} + h.o.t.. \nonumber 
\end{eqnarray}
We also have 
$$
  d \hat{\th}^1_{(s)} = 2 (1+s^2) dz^2 \wedge dz^1 - 2 s \sqrt{1+s^2} d \oz^2 \wedge d \oz^1. 
$$

We expand $d \psi$ and $\hat{\th}^1_{(s)}$ in homogeneous powers of the $(z,t)$ coordinates (w.r.t. parabolic 
scaling, including differentials) as follows 
$$
  d \psi = (d \psi)_0 + (d \psi)_1 + (d \psi)_2 + \cdots. 
$$
Taylor-expanding the above system up to order one we obtain the relations 
\begin{equation}\label{eq:second-line}
\begin{cases}
  (d \psi)_0 \wedge \hat{\th}^1_{(s)} = 0; \\ 
  i (d \psi)_1 \wedge \hat{\th}^1_{(s)} + d \hat{\th}^1_{(s)} + (d v_2) \wedge \hat{\th}^1_{(s)} 
  + 2 i (v_2)^1 d \hat \th' + 2 i (d (v_2)^1) \wedge \hat \th'  = 0. 
  \end{cases}
\end{equation}
The first component is easy to solve setting $(d \psi)_0 = \mu \, \hat{\th}^1_{(s)}$ for some $\mu \in \R$. 

For the second component, recall that we have 
$$
  \hat{\o}^1_{1(s)} = -i (1+ 2 s^2)  \hat \th' ; \qquad \qquad \hat{A}^1_{\ou(s)} = 2 i s \sqrt{1+s^2}. 
$$
It then follows 
$$
 d \hat{\th}^1_{(s)} = \hat{\th}^1_{(s)} \wedge \hat{\o}^1_{1(s)} + \hat{A}^1_{\ou(s)} \hat{\th}' \wedge \hat{\th}^{\ou}_{(s)} 
 =  -i (1+ 2 s^2)  \hat{\th}^1_{(s)} \wedge \hat{\th}' + 2 i s \sqrt{1+s^2} \, \hat{\th}' \wedge \hat{\th}^{\ou}_{(s)}. 
$$
Moreover we have  
$$
  d \hat{\th}' = 2 i h_{1 \ou (s)}  \hat{\th}^1_{(s)} \wedge \hat{\th}^{\ou}_{(s)} = \frac{2i}{1+s^2} \hat{\th}^1_{(s)} \wedge \hat{\th}^{\ou}_{(s)}, 
$$
and that  (up to $\hat{\th}'$)
$$
  d (v_2)^1 = \hat{Z}_{1(s)} (v_2)^1 \hat{\th}^1_{(s)} + \hat{Z}_{\ou(s)} (v_2)^1 \hat{\th}^{\ou}_{(s)}. 
$$
By the above expression of $(v_2)^1$ and \eqref{eq:new-33}, this becomes 
\begin{eqnarray*}
  d (v_2)^1  & = & \left[ (s^2 + 1/2) B_1 (|z_2|^2 - |z_1|^2) + B_2 \bar z_2^2 - C_2 \bar z_1^2 
  + s \sqrt{1+s^2} A_1 z_2^2 - s^2 A_1 z_1^2 \right] \hat{\th}^1_{(s)} \\ 
  & + & \left[ s \sqrt{1+s^2} B_1 (|z_2|^2 - |z_1|^2) - B_2  z_1^2 
  + s^2 A_1 \bar z_2^2 -  s \sqrt{1+s^2} A_1 \bar z_1^2 + C_2 z_2^2 \right] \hat{\th}^{\bar 1}_{(s)} 
  \quad \hbox{ mod } \hat{\th}'.  
\end{eqnarray*}
We next write 
$$
  (d \psi)_1 = (A_3 z_1 + B_3 \oz_1) \hat{\th}^1_{(s)} + (\ov{A}_3 \oz_1 + \ov{B}_3 z_1) \hat{\th}^{\ou}_{(s)}
  + C_3 \hat{\th}'; \qquad \quad \ov{A}_3 = - A_3, \quad \ov{B}_3 = - B_3.  
$$
The $\hat{\th}^{\ou}_{(s)} \wedge \hat{\th}^1_{(s)}$-component of the second equation in \eqref{eq:second-line}  is given by 
$$
   i \left( \ov{A}_3 \oz_1 + \ov{B}_3 z_1 \right) +  \frac{6}{1+s^2} (v_2)^1 = 0. 
$$ 
This determines $A_3$ and $B_3$ by
$$
  i \ov{A}_3 + \frac{6}{1+s^2} (C_2 + D_2) = 0; \qquad \qquad i \ov{B}_3 + \frac{6}{1+s^2} (A_2 + B_2) = 0,  
$$
giving 
\begin{equation}\label{eq:A3B3}
 A_3 = - \frac{3}{8} \frac{2 i s }{\sqrt{1+s^2}} (7 + 6 s^2); \qquad \qquad 
  B_3 = - \frac{3}{8} 2 i  (1 - 6 s^2). 
\end{equation}

Next, the $\hat{\th}' \wedge \hat{\th}^1_{(s)}$ component gives 
$$
  i C_3 + i (1 + 2 s^2) + \hat{T} v_2 - 2 i \hat{Z}_{1(s)} (v_2)^1 = 0. 
$$
Finally, the $\hat{\th}' \wedge \hat{\th}^{\ou}_s $-component of the second equation in \eqref{eq:second-line}   is given by 
$$
  2 i s \sqrt{1+s^2}  - 2 i \hat{Z}_{\ou(s)} (v_2)^1 = 0. 
$$
This is true because, as one can check  
$$
  \hat{Z}_{\ou(s)} (v_2)^1 =  (C_2 + D_2) + \frac{s}{\sqrt{1+s^2}} (A_2 + B_2);       \qquad \qquad 
  \hat{Z}_{1(s)} (v_2)^1= (A_2 + B_2) + \frac{s}{\sqrt{1+s^2}} (C_2 + D_2). 
$$
By a direct computation it follows that 
$$
   \hat{Z}_{\ou(s)} (v_2)^1  = s \left(s^2+1\right)^{1/2}; \qquad
    \qquad  \hat{Z}_{1(s)} (v_2)^1 = \frac{1}{8}  \left(2 s^2+1\right). 
$$
These also imply 
$$
  C_3 =  - \frac{3}{4}  (1 + 2 s^2) = - \frac{3}{2} s^2 - \frac 34. 
$$

Let us now try to integrate for the phase $\psi$. There holds 
$$
   \hat{\th}^1_{(s)} = (1+s^2) dz_1 - s \sqrt{1+s^2} d \ov{z}_1 + h.o.t.;  \qquad 
   \qquad \hat{\th}' = i \left[ (z_1 d \ov{z}_1 - \ov{z}_1 d z_1) + (d \ov{z}_2 - d z_2) \right] + h.o.t..    
$$
In this way, we have that $(d \psi)_1$ becomes 
\begin{eqnarray*}
     & & \left[ (A_3 z_1 + B_3 \ov{z}_1) (1+s^2) + (A_3 \ov{z}_1 + B_3 z_1) s \sqrt{1+s^2} - i  C_3 \oz_1 \right] dz_1   
     + conj. - i C_3 (dz_2 - d \oz_2) \\ 
     & = & \left\{ \left[ (1+s^2) A_3 + s \sqrt{1+s^2} B_3 \right] z_1 +  \left[ (1+s^2) B_3 + s \sqrt{1+s^2} A_3 - i  C_3 \right] 
     \oz_1  \right\} dz_1 + conj. \\ & - & i C_3 (dz_2 - d \oz_2). 
\end{eqnarray*}
Since $ (1+s^2) B_3 + s \sqrt{1+s^2} A_3 - i  C_3  = 0$,   we get 
$$
  (d \psi)_1 = A_4 z_1 dz_1 + B_4 dz_2 + conj., 
$$
with  
\begin{equation} \label{eq:AB4}
 A_4 = - 6 i  s \sqrt{1+s^2}; \qquad \qquad B_4 = - i C_3 =  \frac{3}{4} i  (1 + 2 s^2). 
\end{equation}
Integrating, we find 
$$
  (\psi)_2  = \frac{1}{2} A_4 z_1^2 + B_4 z_2 + conj. .  
$$
Taylor-expanding, we then get 
\begin{equation}\label{zcrs}
  dz'_{\text CR}  = \hat{\th}^1_{(s)} (1 + v_2 + i \psi_2) + 2 i (v_2)^1 \hat{\th}' + h.o.t.. 
\end{equation}

Writing  the 0-th and 2nd order terms of the right hand side, we obtain%
\begin{eqnarray*}
&&[1+A_{1}(z_{1}^{2}+\bar{z}_{1}^{2})+\frac{1}{2}iA_{4}(z_{1}^{2}-\bar{z}%
_{1}^{2})+B_{1}|z_{1}|^{2}+iB_{4}(z_{2}-\bar{z}_{2})] \\
&&\{(1+s^{2})z_{2}dz_{1}-s\sqrt{1+s^{2}}\bar{z}_{2}d\bar{z}%
_{1}-(1+s^{2})z_{1}dz_{2}+s\sqrt{1+s^{2}}\bar{z}_{1}d\bar{z}_{2}\} \\
&&-2[(A_{2}+B_{2})z_{1}+(C_{2}+D_{2})\bar{z}_{1}][z_{1}d\bar{z}_{1}-\bar{z}%
_{1}dz_{1}+d\bar{z}_{2}-dz_{2}].
\end{eqnarray*}
Further expanding this, gives%
\begin{eqnarray} \nonumber
&&(1+s^{2})z_{2}dz_{1}-s\sqrt{1+s^{2}}\bar{z}_{2}d\bar{z}_{1}  \label{SOT} 
-(1+s^{2})z_{1}dz_{2}+s\sqrt{1+s^{2}}\bar{z}_{1}d\bar{z}_{2}   \\
&&+[A_{1}(z_{1}^{2}+\bar{z}_{1}^{2})+\frac{1}{2}iA_{4}(z_{1}^{2}-\bar{z}%
_{1}^{2})+B_{1}|z_{1}|^{2}+iB_{4}(z_{2}-\bar{z}_{2})]   \\
&&\lbrack (1+s^{2})dz_{1}-s\sqrt{1+s^{2}}d\bar{z}%
_{1}]-2[(A_{2}+B_{2})z_{1}+(C_{2}+D_{2})\bar{z}_{1}]  \notag 
 \; \lbrack z_{1}d\bar{z}_{1}-\bar{z}_{1}dz_{1}+d\bar{z}_{2}-dz_{2}].  
\end{eqnarray}
We next set%
\begin{equation*}
w=z_{2}-\bar{z}_{2},
\end{equation*}

\noindent and rewrite the terms involving $z_{2}$ as%
\begin{eqnarray}\label{eq:dz2} \nonumber
z_{2} &=&1+\frac{1}{2}w-\frac{1}{2}|z_{1}|^{2};\text{ \ } \qquad \bar{z}_{2}=1-\frac{%
1}{2}w-\frac{1}{2}|z_{1}|^{2}; \\
dz_{2} &=&\frac{1}{2}dw-\frac{1}{2}(\bar{z}_{1}dz_{1}+z_{1}d\bar{z}_{1});  \qquad \quad %
\text{ }d\bar{z}_{2}=-\frac{1}{2}dw-\frac{1}{2}(\bar{z}_{1}dz_{1}+z_{1}d\bar{%
z}_{1}).
\end{eqnarray}
Write (\ref{SOT}) as%
\begin{eqnarray}
&&C_{6}(z_{2}dz_{1}-z_{1}dz_{2})+D_{6}(\bar{z}_{2}d\bar{z}_{1}-\bar{z}_{1}d%
\bar{z}_{2})  \label{SOT2} \\
&&+[A_{1}(z_{1}^{2}+\bar{z}_{1}^{2})+\frac{1}{2}iA_{4}(z_{1}^{2}-\bar{z}%
_{1}^{2})+B_{1}|z_{1}|^{2}+iB_{4}w][C_{8}dz_{1}+D_{8}d\bar{z}_{1}]  \notag \\
&&-2[(A_{2}+B_{2})z_{1}+(C_{2}+D_{2})\bar{z}_{1}][z_{1}d\bar{z}_{1}-\bar{z}%
_{1}dz_{1}-dw],  \notag
\end{eqnarray}
where 
\begin{equation} \label{eq:CD68}
C_{6}=C_{8}=1+s^{2}; \qquad \quad \text{ \ }D_{6}=D_{8}=-s\sqrt{1+s^{2}}.
\end{equation}
We now substitute $C_{6}(z_{2}dz_{1}-z_{1}dz_{2})$ $+$ $D_{6}(\bar{%
z}_{2}d\bar{z}_{1}-\bar{z}_{1}d\bar{z}_{2})$ $=$ $C_{6}$[$(1+\frac{1}{2}%
w)dz_{1}$ $+$ $\frac{1}{2}z_{1}^{2}d\bar{z}_{1}$ $-$ $\frac{1}{2}z_{1}dw]$ + 
$D_{6}[(1-\frac{1}{2}w)d\bar{z}_{1}$ $+$ $\frac{1}{2}\bar{z}_{1}^{2}dz_{1}$ $%
+$ $\frac{1}{2}\bar{z}_{1}dw]$ into (\ref{SOT2}) and collect terms involving 
$w$ as follows:%
\begin{eqnarray}
&&(iB_{4}C_{8}+\frac{1}{2}C_{6})wdz_{1}+[2(A_{2}+B_{2})-\frac{1}{2}%
C_{6}]z_{1}dw  \label{z1w} \\
&&+(iB_{4}D_{8}-\frac{1}{2}C_{6})wd\bar{z}_{1}+[2(C_{2}+D_{2})+\frac{1}{2}%
D_{6}]\bar{z}_{1}dw.  \notag
\end{eqnarray}
A direct computation shows that %
\begin{equation}
iB_{4}C_{8}+\frac{1}{2}C_{6} =2(A_{2}+B_{2})-\frac{1}{2}C_{6}
\label{z1w-1} 
=(1+s^{2})(-\frac{3}{2}s^{2}-\frac{1}{4}).  
\end{equation}
Similarly, we have%
\begin{equation}
iB_{4}D_{8}-\frac{1}{2}C_{6} =2(C_{2}+D_{2})+\frac{1}{2}D_{6}
\label{z1w-2} 
=s\sqrt{1+s}(\frac{3}{2}s^{2}+\frac{5}{4}).  
\end{equation}
In view of (\ref{z1w-1}) and (\ref{z1w-2}) we can write (\ref{z1w}%
) as%
\begin{equation}
d[(1+s^{2})(-\frac{3}{2}s^{2}-\frac{1}{4})z_{1}w  \label{z1we} 
+s\sqrt{1+s}(\frac{3}{2}s^{2}+\frac{5}{4})\bar{z}_{1}w].  \notag
\end{equation}

\noindent On the other hand, we can write terms only involving $z_{1}$ and $%
\bar{z}_{1}$ in (\ref{SOT2}) as%
\begin{equation}
(K_{11}z_{1}^{2}+K_{\bar{1}\bar{1}}\bar{z}_{1}^{2}+K_{1\bar{1}%
}|z_{1}|^{2})dz_{1}  \label{z11bar} 
+(N_{11}z_{1}^{2}+N_{\bar{1}\bar{1}}\bar{z}_{1}^{2}+N_{1\bar{1}%
}|z_{1}|^{2})d\bar{z}_{1},  \notag
\end{equation}
where 
$$
K_{11} =(A_{1}+\frac{1}{2}iA_{4})C_{8};\text{ \ }  \qquad 
K_{\bar{1}\bar{1}} =\frac{1}{2}D_{6}+(A_{1}-\frac{1}{2}%
iA_{4})C_{8}+2(C_{2}+D_{2}); 
$$
$$
K_{1\bar{1}} =B_{1}C_{8}+2(A_{2}+B_{2}); \qquad 
N_{11} =\frac{1}{2}C_{6}+(A_{1}+\frac{1}{2}iA_{4})D_{8}-2(A_{2}+B_{2}); 
$$
$$
N_{\bar{1}\bar{1}} =(A_{1}-\frac{1}{2}iA_{4})D_{8}; \qquad 
N_{1\bar{1}} =B_{1}D_{8}-2(C_{2}+D_{2}).
$$
Observe that%
\begin{equation}
K_{\bar{1}\bar{1}} =\frac{1}{2}N_{1\bar{1}}=s\sqrt{1+s^{2}}(\frac{3}{2}%
s^{4}+\frac{3}{4}s^{2}-1),  \label{z11bar-1} \qquad
N_{11} =\frac{1}{2}K_{1\bar{1}}=(1+s^{2})(-\frac{3}{2}s^{4}-\frac{9}{4}%
s^{2}+\frac{1}{4}),  
\end{equation}
and 
\begin{equation}
K_{11} =s(1+s^{2})^{3/2}(\frac{3}{2}s^{2}+\frac{15}{4}),  \label{z11bar-2}
\qquad
N_{\bar{1}\bar{1}} =-s^{2}(1+s^{2})(\frac{3}{2}s^{2}-\frac{9}{4}).  
\end{equation}
In view of (\ref{z11bar-1}) and (\ref{z11bar-2}), we can express (%
\ref{z11bar}) as%
\begin{equation}
d\{\frac{1}{3}K_{11}z_{1}^{3}+\frac{1}{3}N_{\bar{1}\bar{1}}\bar{z}%
_{1}^{3}+K_{\bar{1}\bar{1}}\bar{z}_{1}^{2}z_{1}+N_{11}z_{1}^{2}\bar{z}_{1}\}.
\label{z11bare}
\end{equation}
Altogether, from (\ref{z1we}) and (\ref{z11bare}) we obtain $%
\tilde{z}_{\text CR} $ (see (\ref{zcrs})) as follows:%
\begin{eqnarray} \label{eq:tzcr} \nonumber
&& \tilde{z}_{\text CR}  = (1+s^{2})z_{1}-s\sqrt{1+s^{2}}\bar{z}_{1} 
+(1+s^{2})(-\frac{3}{2}s^{2}-\frac{1}{4})z_{1}w+s\sqrt{1+s^2}(\frac{3}{2}%
s^{2}+\frac{5}{4})\bar{z}_{1}w \\
&&+s(1+s^{2})^{3/2}(\frac{1}{2}s^{2}+\frac{5}{4})z_{1}^{3}-s^{2}(1+s^{2})(%
\frac{1}{2}s^{2}-\frac{3}{4})\bar{z}_{1}^{3} 
+s\sqrt{1+s^{2}}(\frac{3}{2}s^{4}+\frac{3}{4}s^{2}-1)\bar{z}_{1}^{2}z_{1}
\\
&&+(1+s^{2})(-\frac{3}{2}s^{4}-\frac{9}{4}s^{2}+\frac{1}{4})z_{1}^{2}\bar{z}_{1} + h.o.t.. \nonumber
\end{eqnarray}

\noindent The CR normal coordinate $z'_{\text CR}$ w.r.t.  Jerison-Lee's frame 
reads 
\begin{equation*}
z'_{\text CR} =\frac{\tilde{z}_{\text CR} }{\sqrt{1+s^{2}}}.
\end{equation*}

We want next to determine the $t$-component of  CR normal coordinates. Recall the definition 
of $\hat{\th}$ and \eqref{eq:dz2}: 
after some cancellations one can check that 
$$
  \hat{\th}' = i \left\{ z_1 d \oz_1 - \oz_1 dz_1 - dw + \frac{1}{2} |z_1|^2 dw - 
  \frac 12 w (\oz_1 dz_1 + z_1 d\oz_1)  \right\}. 
$$
We now need to consider the conformal change of contact form 
$$
 \check{\theta}' = e^{2v} \hat{\th}' = (1 + 2 v_2 + \cdots ) \hat{\th}'. 
$$
Recalling that $v_2 = A_1 (z_1^2 + \oz_1^2) + B_1 |z_1|^2$ we obtain that 
$$
  \check{\theta}' = \left( 1 + 2 A_1 (z_1^2 + \oz_1^2) + 2 B_1 |z_1|^2 \right) i 
  \left\{ \left( z_1 - \frac 12 w z_1 \right) d \oz_1 - \left( \oz_1 + \frac 12 w \oz_1 \right) dz_1 
  - \left( 1 - \frac 12 |z_1|^2 \right) dw  \right\} + h.o.t.. 
$$
From straightforward computations one finds 
\begin{eqnarray*}
\check{\theta}' & = & i \left( z_1 + 2 A_1 (z_1^2 + \oz_1^2) z_1 + 2 B_1 |z_1|^2 z_1 
- \frac{1}{2} w z_1 \right) d \oz_1 \\ &  - &  i \left( \oz_1 + 2 A_1 (z_1^2 + \oz_1^2) \oz_1 + 2 B_1 |z_1|^2 \oz_1 
+ \frac{1}{2} w \oz_1 \right) d z_1 \\ 
& - & i \left( 1 + 2 A_1 (z_1^2 + \oz_1^2) + 2 B_1 |z_1|^2 - \frac{1}{2} |z_1|^2 \right) dw + h.o.t.. 
\end{eqnarray*}
Therefore, from \eqref{eq:tzcr} we deduce 
\begin{eqnarray*}
 d\tilde{z}_{\text CR}  & = & (1 + s^2) dz_1 - s \sqrt{1+s^2} d \oz_1 + (1+s^2) \left( - \frac{3}{2} s^2 - \frac 14  \right)
 (w dz_1 + z_1 dw) \\ & + & s \sqrt{1+s^2} \left( \frac 32 s^2 + \frac 54 \right) (w d \oz_1 + \oz_1 dw) \\ & + & 
 K_{11} z_1^2 dz_1 + N_{\ou \ou} \oz_1^2 d \oz_1 + K_{\ou \ou} (2 \oz_1 z_1 d\oz_1 + \oz_1^2 dz_1) 
 + N_{11} (2 z_1 \oz_1 dz_1 + z_1^2 d \oz_1) + h.o.t.. 
\end{eqnarray*}
One can then expand $\check{\theta}'  + i \bar z'_{\text CR}  dz'_{\text CR}  - i z'_{\text CR}  d \bar z'_{\text CR} $ to find that 
\begin{equation}\label{eq:tcr}
t'_{\text CR}  = - i w (1 + 1/2 |z_1|^2) + i s |z_1|^2 (z_1^2 - \bar z_1^2) + i s^2  (\bar z_1^4 - z_1^4)  +h.o.t..
\end{equation}
We can summarize the above discussion into the following result. 

\begin{pro}\label{l:cr-coord-expl}
The CR normal coordinates on Rossi spheres w.r.t. $\check{\theta}' = e^{2v} \hat{\theta}'$ are given by the formulas $z'_{\text CR} =\frac{\tilde{z}_{\text CR} }{\sqrt{1+s^{2}}}$, 
with $\tilde{z}_{\text CR} $ as in \eqref{eq:tzcr} and $t'_{\text CR} $ as in \eqref{eq:tcr}. 
\end{pro}

\medskip

\noindent We next collect some useful formulas derived from the latter proposition. 
Taylor-expanding $z'_{\text CR}$ one finds 
$$
   |z'_{\text CR} |^2 = |z_1|^2 \left( 1 + \frac{1}{2} |z_1|^2 \right) - s \left( z_1^2 + \bar z_1^2 
   + w(z_1^2 - \bar z_1^2) \right) + \frac{1}{2} s^2 \left( 4 |z_1|^2 - 4 |z_1|^4 - z_1^4 - \bar z_1^4 \right)
   + h.o.t., 
$$
while taking its square we obtain 
\begin{eqnarray} \nonumber \label{eq:zCR^4}
|z'_{\text CR} |^4 & = & |z_1|^4 \left( 1 + |z_1|^2 \right) - s |z_1|^2 \left( (z_1^2 + \bar z_1^2) (2 + |z_1|^2)  
   + 2 w(z_1^2 - \bar z_1^2) \right) \\ & + &  s^2 \left( (z_1^4 + \bar z_1^4) (1 - |z_1|^2) + 2 |z_1|^4 (3 |z_1|^2 - 1) 
   + 2 w (z_1^4 - \bar z_1^4 ) \right) + h.o.t.. 
\end{eqnarray}
The square of $t'_{\text CR} $ is given by 
$$
  (t'_{\text CR})^2 = -  w^2 (1 + |z_1|^2) + 2 s w |z_1|^2 (z_1^2 - \bar z_1^2) + 2 s^2 w (\bar z_1^4 - z_1^4) + h.o.t..
$$
Summing the latter formula and \eqref{eq:zCR^4} we obtain that, up to higher-order terms 
\begin{eqnarray*}
 (\rho'_{\text CR})^4 & = & (1 + |z_1|^2) (|z_1|^4 - w^2) - s |z_1|^2 (z_1^2 + \bar z_1^2) (2 + |z_1|^2) 
 \\ & + & s^2 
 \left[ (z_1^4 + \bar z_1^4) (1 - |z_1|^2) + 2 |z_1|^4 (3 - |z_1|^2) \right]. 
\end{eqnarray*}
It is also useful to expand the quantity $e^{v_2} (\rho'_{\text CR})^{-2}$, related to the conformal covariance 
for the Green's function, which up to higher-order terms is given by 
\begin{eqnarray} \label{eq:v2-rho'}
& & e^{v_2} (\rho'_{\text CR})^{-2}  =  \frac{|z_1|^2+4}{4 \sqrt{(|z_1|^2+1) \left(|z_1|^4-w^2\right)}} 
+ s \frac{\left(z_1^2+\bar z_1^2\right) \left(12 |z_1|^4 +8 |z_1|^2-6 w^2\right)}{8 \left((|z_1|^2+1) \left(|z_1|^4-w^2\right)\right)^{3/2}} 
 \\ & + & s^2 \left[ \frac{ (z_1^4 + \bar z_1^4) ( 20 |z_1|^6 + 8 |z_1|^4  + 4 w^2 - 5 |z_1|^2 w^2)  - 4 |z_1|^{10}+58 |z_1|^6 w^2 
+24 |z_1|^4 w^2  -24 |z_1|^2 w^4 
}{8  \left((|z_1|^2+1) \left(|z_1|^4-w^2\right)\right)^{5/2}} \right]. \nonumber
\end{eqnarray}
Note that w.r.t. the contact from $\hat{\theta} = \frac{1}{2} \hat{\th}'$ the CR normal coordinates and the 
{\em Heisenberg distance} would be $(z_{\text CR}, t_{\text CR}) = \left( \frac{z'_{\text CR}}{\sqrt{2}},  \frac{t'_{\text CR}}{{2}}\right)$ 
and $\rho_{\text CR} = \frac{\rho'_{\text CR}}{\sqrt{2}}$ respectively. 

\

%
%
%
%
%

\section{Proof of Theorem \ref{t:mass}}\label{s:Green}

In this section we determine the Green's function for the conformal sub-Laplacian 
on Rossi spheres, up to an error of order $s^3$. This allows to estimate the 
mass of Rossi spheres, which turns out to be negative for $s \neq 0$ small. This is done by comparing the 
expression of the Green's function in 
CR-normal coordinates, locally near the pole, and by deriving a formal expansion in 
$s$ \underline{globally} away from the pole with respect to the standard coordinates 
$(z_1, z_2)$ of $S^3$. 

\medskip

\subsection{Formal expansion of the Green's function in powers of $s$}

Let $L_s$ denote the conformal sub-Laplacian for the $J_{(s)}$-structure on $S^3$. 
For $s = 0$, the fundamental solution of $L_0 G_0 = 64 \pi^2 \d_p$ with pole at $p = (0,1)$ is given by  
\begin{equation}\label{eq:G0}
 G_0 = 2 ((1 - z_2)(1-\ov{z}_2))^{-\frac{1}{2}}. 
\end{equation}
We next solve formally, up to an error $O(s^3)$, $L_s G_s = 0$ away from $p$ in power series of $s$ in the 
 form 
\begin{equation}\label{eq:Gss}
G_s = G_0 + s G_1 + \frac{1}{2} s^2 (G_2 + \a G_0 - G_3),
\end{equation}
where $G_1, G_2$ are suitable explicit singular functions near $p$, $\alpha \in \R$ and 
$G_3$ is a H\"older continuous function near $p$ for which we would need to determine only $G_3(p)$. 
We chose to expand the second-order term including separately $\a G_0$: this will be useful 
later in order to fix the distributional component of the solution at the pole $p$. In principle this 
should be done also for the first-order term, but by our choice of $G_1$ this further 
correction will not be necessary.

For the above expansion, the following formulas will be used 
\begin{eqnarray} \label{eq:z1z1f} 
& & Z_{1} Z_1 z_1^a (1-z_2)^b (1-\bar z_2)^c = 
z_1^{a-2} (1-z_2)^{b-2} (1-\bar z_2)^c \left(b (|z_2|^2-1) (2 a \bar z_2 (z_2-1)-|z_2|^2+1) \right. \\ & + & \left. (a-1) a \bar z_2^2 (z_2-1)^2+b^2 (|z_2|^2-1)^2\right); \nonumber
\end{eqnarray}
\begin{equation}\label{eq:zouzouf}
  Z_{\bar 1} Z_{\bar 1} z_1^a (1-z_2)^b (1-\bar z_2)^c = (c-1) c z_1^{a+2} (1-z_2)^b (1-\bar z_2)^{c-2}; 
\end{equation}
\begin{eqnarray}\label{eq:zouz1f}
 & &  Z_{\bar 1} Z_{1} z_1^a (1-z_2)^b (1-\bar z_2)^c \\ & = & \nonumber
 z_1^a (1-z_2)^{b-1} \left(-(1-\bar z_2)^{c-1}\right) (a (z_2-1) ((c+1) \bar z_2-1)+b (c (|z_2|^2-1)+(\bar z_2-1) z_2)); 
\end{eqnarray}
\begin{equation}\label{eq:z1zouf}
  Z_1 Z_{\bar 1} z_1^a (1-z_2)^b (1-\bar z_2)^c = 
 -c z_1^a (1-z_2)^{b-1} (1-\bar z_2)^{c-1} ((a+1) \bar z_2 (z_2-1)+b (|z_2|^2-1)),
\end{equation}
with similar ones for $\bar z_1^a (1-z_2)^b (1-\bar z_2)^c$, passing to conjugates. 

To find the first-order correction $G_1$, we differentiate the relation $L_s G_s = 0$ with 
respect to $s$, evaluating it for $s = 0$.  Using \eqref{eq:dot-L-ddot-L} and \eqref{eq:dot-R-ddot-R}, this yields 
$$
  L_0 G_1 = - \dot{L} G_0 =  8 Z_1 Z_1 G_0 + 8 Z_{\ov{1}} Z_{\ov{1}} G_0 
  \qquad \quad \hbox{ on } S^3 \setminus \{p\}, 
$$
where $\dot{L} = \frac{d}{ds}|_{s=0} L_s$. The right-hand side is given by 
$$
 \frac{12 \left((\bar z_2-1)^2 \bar z_1^2+z_1^2 (z_2-1)^2\right)}{((\bar z_2-1) (z_2-1))^{5/2}}. 
$$
By  formulas \eqref{eq:z1z1f}-\eqref{eq:z1zouf}, the first-order 
correction $G_1$ to $G_s$ can be chosen as 
\begin{equation}\label{eq:G1}
G_1 =   \frac{1}{2}(z_1^2 + \ov{z}_1^2) \left[ \frac{1}{1-z_2} + \frac{1}{1-\ov{z}_2} + 2  \right]
  \frac{1}{\left( (1-z_2) (1-\ov{z}_2) \right)^{\frac 12}}.
\end{equation}

We pass next to the second order expansion for $G_s$: we will  find it up to a smooth function that can be 
determined at $p$, which is enough for our purposes. 
Differentiating the relation $L_s G_s = 0$ twice with 
respect to $s$ and evaluating at $s = 0$ we obtain (with analogous notation to above for the $s$-derivatives)
$$
  L_0 \ddot{G} = - 2 \dot{L} \dot{G} - \ddot{L} G_0. 
$$
Recalling from \eqref{eq:dot-L-ddot-L}, \eqref{eq:dot-R-ddot-R} that $\ddot{L} = 4 L_0$, we have  
$$
  L_0 (\ddot{G} + 4 G_0) = -  2 \dot{L} \dot{G} = 16 Z_1 Z_1 G_1 + 16 Z_{\ov{1}} Z_{\ov{1}} G_1. 
$$
It is possible to show by direct computation, again from \eqref{eq:z1z1f}-\eqref{eq:z1zouf}, that $16 Z_1 Z_1 G_1 + 16 Z_{\ov{1}} Z_{\ov{1}} G_1$ 
equals 
\begin{eqnarray} \nonumber 
 &  & \frac{-1}{\left((1-z_2) (1-\ov{z}_2)  \right)^{\frac 72}} \left[  
z_1^4 \left( 30 (z_2-1)^3 + 6 (z_2-1)^2 (\ov{z}_2-1) - 12 (z_2-1)^3 (\ov{z}_2-1) \right) \right. \\ \nonumber 
& + & \ov{z}_1^4  \left(  30 (\ov{z}_2-1)^3 + 6 (\ov{z}_2-1)^2 (z_2-1) - 12 (\ov{z}_2-1)^3 (z_2-1) \right) \\
& + & 30 (\ov{z}_2-1)^5 + 30 (z_2-1)^5 + 18 (\ov{z}_2-1)^4 (z_2-1) + 18 (z_2-1)^4 (\ov{z}_2-1) \\ 
& + & 6 (\ov{z}_2-1)^5 (z_2-1)^2 + 6 (z_2-1)^5 (\ov{z}_2-1)^2 - 18 (\ov{z}_2-1)^4 (z_2-1)^3 - 
18 (\ov{z}_2-1)^3 (z_2-1)^4 \nonumber \\ 
& - & \left. 12 (z_2-1)^3 (\ov{z}_2-1)^5 - 12 (\ov{z}_2-1)^3 (z_2-1)^5 \right],  \nonumber 
\end{eqnarray}
where we grouped the terms by homogeneity in $z_2-1$ and $\bar z_2-1$.

We can invert $L_0$ explicitly for the terms with factors $z_1^4$ and $\ov{z}_1^4$. The solution is given by 
$$
  G_{2,1} := \frac{(z_1^4+ \ov{z}_1^4) \, g_{2,1}}{\left((1-z_2)(1+\ov{z}_2)\right)^{\frac{5}{2}}}, 
$$
where 
\begin{eqnarray*}
 g_{2,1} & := & \frac{3}{8} (\ov{z}_2-1)^2 +\frac{3}{8} (z_2-1)^2 + \frac{1}{4}  (z_2-1)  (\ov{z}_2-1) 
  + \frac{3}{2}  (z_2-1)^2  (\ov{z}_2-1)^2 \\ & - & \frac{3}{4}  (z_2-1)^2  (\ov{z}_2-1) 
  - \frac{3}{4}  (z_2-1)  (\ov{z}_2-1)^2. 
\end{eqnarray*}
For the other terms, we can only find an explicit approximate solution. We set  
\begin{eqnarray*}
g_{2,2} & = & (z_2-1)^4 + (\bar z_2-1)^4 - \frac{4}{3} (z_2-1)^4 (\bar z_2-1) - \frac{4}{3} (\bar z_2-1)^4 (z_2-1) 
\\ & + & 4 (\bar z_2-1)^3 (z_2-1)^2 + 4 (z_2-1)^3 (\bar z_2-1)^2 \\ & + & \frac{11}{3} (z_2-1)^4 (\bar z_2-1)^2 
+   \frac{11}{3} (\bar z_2-1)^4 (z_2-1)^2 + 6 (z_2-1)^3 (\bar z_2-1)^3,   
\end{eqnarray*} 
and  
$$
  G_{2,2} := \frac{3}{4} \frac{g_{2,2}}{\left((1-z_2)(1-\ov{z}_2)\right)^{\frac{5}{2}}}. 
$$
Defining 
\begin{equation}\label{eq:G2}
 G_2 =  G_{2,1}+G_{2,2}, 
\end{equation}
still by \eqref{eq:z1z1f}-\eqref{eq:z1zouf} one finds that 
\begin{equation}\label{eq:G2'}
 L_0 \, G_2 - 16 Z_1 Z_1 G_1 - 16 Z_{\ov{1}} Z_{\ov{1}} G_1 = - 12 \frac{(z_2-1)^2 + (\bar z_2-1)^2 - 3 (z_2-1)(\bar z_2-1)}{\left((1-z_2)(1-\ov{z}_2)\right)^{\frac{1}{2}}} =: \Xi(z_2,\bar z_2),  
\end{equation}
with the right-hand side now bounded on  $S^3$. 

It will be now sufficient to add a more regular correction (which is H\"older continuous by standard regularity theory) to solve the equation for $G_2$ pointwise, away from $p$. 
%
%
%
From \eqref{eq:G2'}, setting $G_3 = L_0^{-1} \Xi(z,w)$  
we then find that 
$$
  L_0 (G_2 - G_3) - 16 Z_1 Z_1 G_1 - 16 Z_{\ov{1}} Z_{\ov{1}} G_1  = 0 \qquad 
  \quad \hbox{ on } S^3 \setminus \{p\}, 
$$
which corresponds to \eqref{eq:Gss} up to the term $s^2 \, \a \, G_0$, which will be determined later. 
To obtain $G_3(p)$, we use the Green's representation formula, convoluting  $\Xi(z_2,\bar z_2)$ with $G_0$: 
$$
   G_3(p) = \frac{1}{64 \pi^2} \int_{S^3}  - 24 \frac{(z_2-1)^2 + (\bar z_2-1)^2 - 3 (z_2-1)(\bar z_2-1)}{\left((1-z_2)(1+\ov{z}_2)\right)} 
   \hat{\theta} \wedge d \hat{\theta}. 
$$
%
%
%
%
%
%
The Taylor expansion of the integrand in $z_2, \ov{z}_2$ is 
\begin{eqnarray*}
\left(24 \bar z_2^5+24 \bar z_2^4+24 \bar z_2^3+24 \bar z_2^2-24\right) + 
\left(-24 \bar z_2^5-24 \bar z_2^4-24 \bar z_2^3-24 \bar z_2^2-48 \bar z_2\right) z_2 \\ 
+ (24-24 \bar z_2) z_2^2 +(24-24 \bar z_2) z_2^3 
+(24-24 \bar z_2) z_2^4
+(24-24 \bar z_2) z_2^5 + \cdots . 
\end{eqnarray*}
Integrated, this gives 
$$
 \int_{S^3}  (24 + 48 |z_2|^2) \, \hat{\theta} \wedge d \hat{\theta} = 48 \, 2 \pi^2 + 96 \, \pi^2 =  192 \pi^2, 
$$
which implies that 
\begin{equation}\label{eq:G3}
G_3(p) =  3. 
\end{equation}

In conclusion, we found that 
$$\ddot{G} = G_2 - G_3 + \a G_0, 
$$ 
i.e. \eqref{eq:Gss}, where $\a$ is a real number to be determined later. 
We proved therefore 
the following result. 

\begin{pro}\label{l:formal}
For every compact set $K$ in $S^3 \setminus \{p\}$, $p = (0,1)$, there exists a constant 
$C_K > 0$ such that 
the  function $G_s$ in \eqref{eq:Gss} satisfies 
$$
  | L_s G_s | \leq C_K s^3 \qquad \quad \hbox{ on }  \quad K. 
$$
\end{pro}

\medskip 

\subsection{Rigorous estimates}

We prove next that the function $G_s$ in Lemma \ref{l:formal} well matches with the 
expression of the Green's function of $L_s$ in CR normal coordinates. 
Recall from the end of Section \ref{s:CR-Rossi} that $\rho_{\text CR}^2 = \frac{1}{2} (\rho'_{\text CR})^2$: 
then from \eqref{eq:v2-rho'} we obtain that 
\begin{eqnarray}\label{eq:Green-CR-m-2} 
& & 2 e^{v_2} \rho_{\text CR} ^{-2}  =  \frac{|z_1|^2+4}{\sqrt{(|z_1|^2+1) \left(|z_1|^4-\mathtt{w}^2\right)}} 
+ s \frac{\left(z_1^2+\bar z_1^2\right) \left(12 |z_1|^4 +8 |z_1|^2-6 \mathtt{w}^2\right)}{2 \left((|z_1|^2+1) \left(|z_1|^4-\mathtt{w}^2\right)\right)^{3/2}} 
\\ & + & s^2 \left[ \frac{ (z_1^4 + \bar z_1^4) ( 20 |z_1|^6 + 8 |z_1|^4  + 4 \mathtt{w}^2 - 5 |z_1|^2 \mathtt{w}^2)  - 4 |z_1|^{10}+58 |z_1|^6 \mathtt{w}^2 
+24 |z_1|^4 \mathtt{w}^2  -24 |z_1|^2 \mathtt{w}^4 
}{2  \left((|z_1|^2+1) \left(|z_1|^4-\mathtt{w}^2\right)\right)^{5/2}} \right] \nonumber \\ 
& + & O(s^3 \rho^{-2}),  \nonumber
\end{eqnarray}
where $\mathtt{w} = z_2 - \bar z_2$. 
Given the covariance property of the Green's function ($G_{(\tilde{\th})} = e^{u} G_{(\th)}$ if $\tilde{\th} = 
e^{2u} \th$), 
we aim to compare  this expression to the  function $G_s$ 
in Lemma \ref{l:formal} on a suitable small annulus centered around $p$. 
We do it 
term by term for the Taylor series in $s$, and for this purpose the following  formulas will be useful.
Since $z_2 - \bar{z}_2$ is purely imaginary, we can write 
$$
  |z_1|^4 + |z_2 - \bar{z}_2|^2 = |z_1|^4 - (z_2 - \bar{z}_2)^2 = 
  \left( |z_1|^2 + (z_2 - \bar{z}_2) \right) \left( |z_1|^2 - (z_2 - \bar{z}_2) \right).
$$
As $|z_1|^2 + |z_2|^2 = 1$, we get
\begin{equation}\label{z4w2}
 |z_1|^4 - \mathtt{w}^2 = |z_1|^4 + |z_2 - \bar{z}_2|^2 = (1 + z_2) (1+\bar{z}_2)(1 - z_2) (1-\bar{z}_2). 
\end{equation}
 Setting $\mathtt{v} = z_2 + \bar{z}_2 - 2$ (which is real), we have that $z_2 = 1 + \frac{\mathtt{v}}{2} + \frac{\mathtt{w}}{2}$, 
 which implies 
 $$
   |z_2|^2 = 1 + \mathtt{v} + \frac{\mathtt{v}^2}{4} - \frac{\mathtt{w}^2}{4} + o(\rho^4). 
 $$
Squaring this relation, we obtain 
\begin{equation}\label{eq:square}
  |z_1|^4 = \mathtt{v}^2 + \frac{\mathtt{v}^3}{2} - \frac{\mathtt{v} \mathtt{w}^2}{2} + o(\rho^6). 
\end{equation} 
We also have that $|z_2|^2 = 1 + \mathtt{v}$ up to an error $O(\rho^4)$, so $\mathtt{v}^2 = - |z_1|^4 + o(\rho^4)$. These imply 
that 
\begin{equation}\label{eq:id-2}
 |z_2|^2 + 1 -(z_2 + \bar z_2) = \frac{1}{4} |z_1|^2 - \frac{1}{4} \mathtt{w}^2 = \frac{1}{4} |z_1|^2 + \frac{1}{4} |\mathtt{w}|^2
 + o(\rho^4).  
\end{equation}
Furthermore, there holds 
\begin{equation}\label{eq:id-3}
  1 + z_2 + \bar{z}_2 + |z_2|^2 \simeq 3 + \mathtt{v} + |z_2|^2 \simeq 2 + 2 |z_2|^2 = 4 - 2 |z_1|^2
   + o(\rho^2). 
\end{equation}
Recalling our notation from Section \ref{s:CRcoord}, we have then the following result.

\medskip

\begin{lem}\label{lem:comparison}
For $\a = - \frac 34$, the following estimate holds  
\begin{equation}\label{eq:eu2G}
2 e^{v_2} \rho_{\text CR} ^{-2}  = \left( G_0 + s G_1 + \frac{1}{2} s^2 (G_2 + \a G_0) \right) 
 +  o(s^2) O''(\rho^{-2}) + o_\rho (1), 
\end{equation}
where $ o_\rho (1) \to 0$ as $\rho \to 0$. 
\end{lem}

\begin{pf} We analyse separately different orders in $s$ for the left-hand side and the 
first term in the right-hand side of \eqref{eq:eu2G}. 

\
  
\noindent {\bf Zero-th order in $s$.} Recalling that $G_0 = 2 \left((1-z_2)(1-\bar{z}_2) \right)^{- \frac 12}$, we need to compare the two quantities 
\begin{equation}\label{eq:two-quant}
  \frac{|z_1|^2+4}{\sqrt{(|z_1|^2+1) \left(|z_1|^4-\mathtt{w}^2\right)}}   \qquad 
  \hbox{ and } \qquad \frac{2}{ \left((1-z_2)(1-\bar{z}_2) \right)^{\frac 12}}. 
\end{equation}
Taylor-expanding the terms involving $|z_1|^2$ in the left-hand side we are left with comparing 
$$
  \frac{4 \left( 1 - \frac{1}{4} |z_1|^2 \right)}{\sqrt{\left(|z_1|^4-\mathtt{w}^2\right)}}   \qquad 
  \hbox{ and } \qquad \frac{2}{ \left((1-z_2)(1-\bar{z}_2) \right)^{\frac 12}}. 
$$
Using \eqref{z4w2} and multiplying by  $\left((1-z_2)(1-\bar{z}_2) \right)^{\frac 12}$, we are left with the comparison of 
$$
  \frac{4 \left( 1 - \frac{1}{4} |z_1|^2 \right)}{\left((1+z_2)(1+\bar{z}_2) \right)^{\frac 12}}   \qquad 
  \hbox{ and } \qquad 2. 
$$
From \eqref{eq:id-3} we are left with comparing 
$$
  \frac{4 \left( 1 - \frac{1}{4} |z_1|^2 \right)}{\left( 4 - 2 |z_1|^2 \right)^{\frac 12}}   \qquad 
  \hbox{ and } \qquad 2, 
$$
which holds true up to an error of order $O(\rho^4)$. Therefore the two quantities in \eqref{eq:two-quant} 
coincide up to an error of order $O(\rho^2)$.

%
%
%
%
%
%
%
%

\ 

\noindent {\bf First order in $s$.}
Recalling \eqref{eq:G1}, we have that 
$$
  G_1 = \frac{1}{4} (z_1^2 + \bar z_1^2) \left[ \frac{4-3 z_2 - 3 \bar z_2  + 2 |z_2|^2}{(1-z_2)(1 - \bar z_2)} \right] G_0. 
$$
Considering the first-order term in $s$ of \eqref{eq:Green-CR-m-2},
we need to compare the two quantities 
$$
  \frac{\left(z_1^2+\bar z_1^2\right) \left(12 |z_1|^4 +8 |z_1|^2-6 \mathtt{w}^2\right)}{2\left((|z_1|^2+1) \left(|z_1|^4-\mathtt{w}^2\right)\right)^{3/2}}  
  \qquad 
    \hbox{ and } \qquad \frac{1}{4} (z_1^2 + \bar z_1^2) \left[ \frac{4-3 z_2 - 3 \bar z_2  + 2 |z_2|^2}{(1-z_2)(1 - \bar z_2)} \right] G_0. 
$$
Using the expression of $G_0$, dividing by $\left(z_1^2+\bar z_1^2\right) $ and multiplying by $2$ we need to compare 
$$
  \frac{\left(12 |z_1|^4 +8 |z_1|^2-6 \mathtt{w}^2\right)}{\left((|z_1|^2+1) \left(|z_1|^4-\mathtt{w}^2\right)\right)^{3/2}}  
  \qquad 
    \hbox{ and } \qquad \left[ \frac{4-3 z_2 - 3 \bar z_2  + 2 |z_2|^2}{\left( (1-z_2)(1 - \bar z_2) \right)^{\frac 32}} \right]. 
$$
Using \eqref{eq:zCR^4}, this is equivalent to the comparison of 
$$
  \frac{\left(12 |z_1|^4 +8 |z_1|^2-6 \mathtt{w}^2\right)}{\left((|z_1|^2+1) (1+z_2)(1 + \bar z_2)\right)^{3/2}}  
  \qquad 
    \hbox{ and } \qquad 4-3 z_2 - 3 \bar z_2  + 2 |z_2|^2.
$$
Using \eqref{eq:id-3} and Taylor-expanding the left-hand side  we arrive to comparing 
$$
  \frac{(1 - 3/4 |z_1|^2) \left(12 |z_1|^4 +8 |z_1|^2-6 \mathtt{w}^2\right)}{8}  
  \qquad 
    \hbox{ and } \qquad 4-3 z_2 - 3 \bar z_2  + 2 |z_2|^2.
$$
Using instead \eqref{eq:id-2} we transform the right-hand side, arriving to the comparison of 
$$
  \frac{(1 - 3/4 |z_1|^2) \left(12 |z_1|^4 +8 |z_1|^2-6 \mathtt{w}^2\right)}{8}  
  \qquad 
    \hbox{ and } \qquad |z_1|^2 + \frac 34 |z_1|^4 - \frac 34 \mathtt{w}^2, 
$$
which is again true up to an error of order $O(\rho^6)$. Therefore, we get matching 
of the first-order terms in $s$ in both sides of \eqref{eq:eu2G} up to an error $O(\rho^2)$.

\
  
\noindent {\bf Second order in $s$.}
Recalling again \eqref{eq:Green-CR-m-2} and the fact that $G_2$ comes with a factor $\frac 12$, let us first  compare
$$
  \frac{ (z_1^4 + \bar z_1^4) ( 20 |z_1|^6 + 8 |z_1|^4  + 4 \mathtt{w}^2 - 5 |z_1|^2 \mathtt{w}^2)  
  }{2  \left((|z_1|^2+1) \left(|z_1|^4-\mathtt{w}^2\right)\right)^{5/2}}   \quad 
      \hbox{ and } \quad \frac 12 G_{2,1} := \frac 12 \frac{(z_1^4+ \ov{z}_1^4) \, g_{2,1}}{\left((1-z_2)(1+\ov{z}_2)\right)^{\frac{5}{2}}}, 
$$
where, up to order $O(\rho^8)$
$$
 g_{2,1}:=  \frac{3}{8} (\ov{z}_2-1)^2 +\frac{3}{8} (z_2-1)^2 + \frac{1}{4}  (z_2-1)  (\ov{z}_2-1) 
  - \frac{3}{4}  (z_2-1)^2  (\ov{z}_2-1) 
  - \frac{3}{4}  (z_2-1)  (\ov{z}_2-1)^2. 
$$
Factoring out $ (z_1^4 + \bar z_1^4)$ and using \eqref{z4w2}, \eqref{eq:id-3} we need to compare 
$$
  \left( 1 - \frac 54 |z_1|^2 \right) \frac{  ( 20 |z_1|^6 + 8 |z_1|^4  + 4 \mathtt{w}^2 - 5 |z_1|^2 \mathtt{w}^2)  
  }{64}   \quad 
      \hbox{ and } \quad   \frac 12 g_{2,1}. 
$$
Expanding $g_{2,1}$ and using \eqref{eq:square}  we arrive to the comparison of 
$$
 \left( 1 - \frac 54 |z_1|^2 \right) \frac{  ( 20 |z_1|^6 + 8 |z_1|^4  + 4 \mathtt{w}^2 - 5 |z_1|^2 \mathtt{w}^2)  
   }{64}   \quad 
       \hbox{ and } \quad  \frac 12  \frac{1}{16} \left(-3 \mathtt{v}^3+4 \mathtt{v}^2+3 \mathtt{v} \mathtt{w}^2+
       2 \mathtt{w}^2\right), 
$$
which is correct, up to an error of order $O(\rho^{12})$.

We need next to compare 
$$
  \left[ \frac{  - 4 |z_1|^{10}+58 |z_1|^6 \mathtt{w}^2 
  +24 |z_1|^4 \mathtt{w}^2  -24 |z_1|^2 \mathtt{w}^4 
  }{2  \left((|z_1|^2+1) \left(|z_1|^4-\mathtt{w}^2\right)\right)^{5/2}} \right] \qquad \hbox{ and } 
  \qquad \frac 12 \frac{3}{4} \frac{g_{2,2}}{\left((1-z_2)(1-\ov{z}_2)\right)^{\frac{5}{2}}},
$$
where, up to higher order terms 
\begin{eqnarray*}
g_{2,2} & = & (z_2-1)^4 + (\bar z_2-1)^4 - \frac{4}{3} (z_2-1)^4 (\bar z_2-1) - \frac{4}{3} (\bar z_2-1)^4 (z_2-1) 
\\ & + & 4 (\bar z_2-1)^3 (z_2-1)^2 + 4 (z_2-1)^3 (\bar z_2-1)^2. 
\end{eqnarray*} 
Using again \eqref{z4w2}, we then need to compare 
$$
  \left[ \frac{  - 4 |z_1|^{10}+58 |z_1|^6 \mathtt{w}^2 
  +24 |z_1|^4 \mathtt{w}^2  -24 |z_1|^2 \mathtt{w}^4 
  }{2  \left((|z_1|^2+1) (1-z_2)(1-\ov{z}_2)\right)^{5/2}} \right] \qquad \hbox{ and } 
  \qquad \frac{3}{8}  g_{2,2}.
$$
As before, we are then comparing 
$$
  \left( 1 - \frac 54 |z_1|^2 \right) \frac{  - 4 |z_1|^{10}+58 |z_1|^6 \mathtt{w}^2 
    +24 |z_1|^4 \mathtt{w}^2  -24 |z_1|^2 \mathtt{w}^4 
     }{64}   \quad 
         \hbox{ and } \quad \frac{3}{8}  g_{2,2}.
$$
In fact, we can add to $G_2$ any multiple of $G_0$. In the latter formula, 
we can then replace $g_{2,2}$ with $\tilde g_{2,2}$, 
where 
$$
  \tilde g_{2,2} = g_{2,2} - 2 (z_2-1)^2 (\bar z_2-1)^2. 
$$
It turns out that
$$
   \frac 38 \tilde g_{2,2} = \frac{1}{16} \mathtt{v} \left(\mathtt{v}^4-4 \mathtt{v}^2 \mathtt{w}^2+6 \mathtt{v} 
   \mathtt{w}^2+3 \mathtt{w}^4\right). 
$$
Using \eqref{eq:square} and the previous formula to expand $|z_1|^2$ as 
$|z_1|^2 = - \mathtt{v} - \frac{1}{4} \mathtt{v}^2 + \frac 14 \mathtt{w}^2$, the left-hand side in the 
above formula becomes 
$$
  \mathtt{v}^5/16 + (3 \mathtt{v}^2 \mathtt{w}^2)/8 - (\mathtt{v}^3 \mathtt{w}^2)/4 + 
  (3 \mathtt{v} \mathtt{w}^4)/16 + O(\rho^{12}),
$$
so it coincides with the right-hand side, i.e. with  $\frac 38 \tilde g_{2,2}$ up to error terms 
of order $O(\rho^{12})$. Therefore, also the second-order terms in $s$ of both sides of 
\eqref{eq:eu2G} coincide up to an error of order $O(\rho^{-2})$. 

It is standard to check that the above matching also holds up to computing first- 
and second-order derivatives, which then implies the conclusion. 
\end{pf}

\

\begin{pfn} {\sc of Theorem \ref{t:mass}.}
Consider a small annulus of the form 
$$
  A_r := \left\{ r \leq \rho \leq 2r \right\}, 
$$
and a smooth cut-off function $\chi_r$ satisfying 
$$
  \begin{cases}
  \chi_r = 1 & \hbox{ on } \{ \rho \leq r \}; \\ 
   \chi_r = 1 & \hbox{ on } \{ \rho \geq 2 r \}; \\ 
   |\nabla_b \chi_r| \leq \frac{C}{r}; & 
   |\nabla_b^2 \chi_r| + |\nabla_T \chi_r| \leq  \frac{C}{r^2}.   
  \end{cases}
$$
If $v$ is the conformal factor as in Proposition \ref{p:CRnormcoord} then, 
with obvious notation, the Green's function conformally transforms as $G_{\theta} = e^{-v} G_{\hat{\th}}$. 
Consider then  the function 
$$
  \check{G}_s = \chi_r \left( 2 \rho_{\text CR} ^{-2} - \frac 12 G_3(p) s^2 \right) + (1-\chi_r) e^{-v} G_{\hat \th}. 
$$
From the conformal covariance of $L_s$,  Lemma \ref{l:formal} and Lemma \ref{lem:comparison} 
it follows that, applying the conformal sub-Laplacian with respect to the contact form $\th$: 
$$
|L_{s}^v \check{G}_s| \leq C_r \, o(s^2) \qquad \quad \hbox{ pointwise on } S^3.
$$
It then follows from standard regularity theory that the Green's function $G_\th$ of the 
conformal sub-Laplacian satisfies $\| G_\th - \check{G}_s\|_{L^\infty(S^3)} = o(s^2)$. 
Sending $s$ to zero and recalling that $G_3(p) = 3$, we deduce 
\begin{equation}\label{eq:As2}
  A = - \frac{3}{2} s^2 + o(s^2). 
\end{equation}
Therefore, given that  $m = 12 \pi A$ (see \eqref{eq:G3} and \eqref{eq:mass-A}), 
we obtain the conclusion. 
\end{pfn}

\

\section{Proof of Theorem \ref{t:SQ}}\label{s:min}

\noindent In this section we prove Theorem \ref{t:SQ} by an implicit function argument 
and some asymptotic expansions, which crucially use also Theorem \ref{t:mass}.

We start by analysing the relation of the CR Sobolev quotient on Rossi spheres with 
 the minimizers on standard spheres found in 
\cite{JLExtr}. Recall  that in \cite{JLCRYam} it was proved that for any three-dimensional CR manifold one has 
$\mathcal{Y}(M,J) \leq \mathcal{Y}(S^3, J_{S^3})$, which in particular implies 
\begin{equation}\label{eq:ineq-quotient-s}
\mathcal{Y}(S^3,J_{(s)}) \leq \mathcal{Y}(S^3, J_{S^3} = J_{(0)}). 
\end{equation}
In \cite{JLExtr} it was proven that $\mathcal{Y}(S^3, J_{S^3})$ is precisely 
attained by the following functions, up to composing $(z_1,z_2)$ with elements of $SU(2)$ 
\begin{equation}\label{eq:var-l}
 \var_\l =  \l \left(  \frac{\left(  |z_1|^2 + |z_2+1|^2 \right)^2 -  (z_2 - \ov{z}_2)^2}{\left( \l ^2 |z_1|^2 + |z_2+1|^2 \right)^2 
  - \l ^4 (z_2 - \ov{z}_2)^2} \right)^{\frac{1}{2}}; \qquad \quad \l > 0. 
\end{equation}
Recalling that $\hat{\th} \wedge d \hat{\th}$ is a volume form double w.r.t. the 
Euclidean one, the $\var_\l$'s satisfy the following normalization condition 
\begin{equation}\label{eq:norm-var-l}
  \int_{S^3} \var_\l^4 \, \hat{\th} \wedge d \hat{\th} = 4 \pi^2 \qquad \qquad \hbox{ for all } \l > 0. 
\end{equation}
On the standard $S^3$, see \cite{FS74}, the \emph{Folland-Stein
space} $\mathfrak{S}^{1,2}(S^3)$ is defined as the completion of the (complex-valued) $%
C^\infty$ functions on $S^3$ with respect to the norm 
\begin{equation*}
\|u\|_{\mathfrak{S}^{1,2}} := \left( \int_{S^3} (u_{, 1} \overline{u}_{, \overline{1}} + u_{, \overline{1}} 
\overline{u}_{, 1}) \, \hat{\th} \wedge d \hat{\th} \right)^{\frac 12} + \left(
\int_{S^3} |u|^2 \hat{\th} \wedge d \hat{\th}  \right)^{\frac 12}.
\end{equation*}
Notice that, for $|s|$ small, this defines an equivalent norm on Rossi 
spheres too: from now on, this will be assumed  understood.

We show next that, if a minimizer for the CR-Sobolev quotient on Rossi spheres exists for 
 $|s|$ small, it must be close in $\mathfrak{S}^{1,2}(S^3)$ to some function 
$\var_\l$ as in \eqref{eq:var-l}. We have indeed the following result. 

\begin{lem}\label{l:closeness-S12}
Fix $s \in \R$, $|s|$ small. Assume $u_s > 0$ attains $\inf Q_{(s)} = \mathcal{Y}(S^3,J_{(s)})$. Then, if $u_s$ is normalized so that  $\int_{S^3} u_s^4 \, \hat{\th} \wedge d \hat{\th} = 4 \pi^2$, 
up to a homogeneous action on $S^3$ there exists $\l > 0$ such that 
$$
  \| u_s  - \var_\l \|_{\mathfrak{S}^{1,2}(S^3)} = o_s(1), 
$$
where $o_s(1) \to 0$ as $s \to 0$. 
\end{lem}

\begin{pf}
It is sufficient to notice that, if $Z_{1(s)}$ is as in \eqref{eq:Js}, then for all 
smooth $u$'s one has  
$$
\int_{S^3} (Z_{1(s)} u \, Z_{\ou(s)} \overline{u} + Z_{\ou(s)} u \, 
Z_{1(s)} \overline{u}) \, \hat{\th} \wedge d \hat{\th}  = (1+o_s(1)) 
\int_{S^3} (Z_{1} u \, Z_{\ou} \overline{u} + Z_{\ou} u \, 
Z_{1} \overline{u}) \, \hat{\th} \wedge d \hat{\th}. 
$$
Since we are assuming $u_s$ to be normalized in $L^4(S^3)$ as in the statement, 
its $\mathfrak{S}^{1,2}(S^3)$-norm is uniformly bounded from above, and therefore 
$$
\int_{S^3} (Z_{1} u_s \, Z_{\ou} \overline{u}_s + Z_{\ou} u_s \, 
Z_{1} \overline{u}_s) \, \hat{\th} \wedge d \hat{\th}  = \int_{S^3} (Z_{1(s)} u_s \, Z_{\ou(s)} \overline{u}_s + Z_{\ou(s)} u_s \, 
Z_{1(s)} \overline{u}_s) \, \hat{\th} \wedge d \hat{\th}  
+ o_s(1). 
$$
This relation implies that $u_s$ is nearly a minimizer also for $\mathcal{Y}(S^3, J_{S^3} = J_{(0)})$ and 
therefore, since the minimizers of the latter quantity must be of the form \eqref{eq:var-l}, 
the conclusion follows. 
\end{pf}

\medskip

\subsection{Finite-dimensional reduction}

Let $\var_\l$ be as in \eqref{eq:var-l}, and define the following family of functions
\begin{equation}\label{eq:mathcalM}
 \mathcal{M} = \left\{ \var_\l(U (\cdot)) : \; | \; \l > 0, U \in SU(2)  \right\}. 
\end{equation}
Even though $SU(2)$ is a four-dimensional Lie group, since $\var_\l$ is invariant by a complex rotation in $z_1$, 
the result of these compositions is a set of  three dimensions. 
We previously saw  that the functions in $\mathcal{M}$ 
 are global minimizers of the CR-Sobolev quotient $Q_{(s)}$
 on the standard $S^3$ when $s = 0$, where 
\begin{equation}\label{eq:Qs}
   Q_{(s)} (u) = \frac{\int_{S^3} u L_s u \, \hat{\th} \wedge d \hat{\th}}{\left( \int_{S^3} u^4 \hat{\th} \wedge d \hat{\th} \right)^{\frac 12}}. 
\end{equation}
In \cite{MuUg}, Lemma 5, it was proved that the linearization of the Yamabe equation (with 
$s = 0$) at $\mathcal{M}$  is minimally degenerate, in the sense that its kernel coincides with the tangent space to 
$\mathcal{M}$.

 As a consequence, one has that 
the  CR-Sobolev 
quotient on the standard sphere is \emph{non-degenerate in the sense of Bott} on $\mathcal{M}$. Thanks to this fact and to Lemma \ref{l:closeness-S12}, 
for $s$ small we can characterize with particular precision all the solutions of the CR-Yamabe equation 
lying in a fixed neighborhood (in $\mathfrak{S}^{1,2}$) of the manifold $\mathcal{M}$, 
and in particular the (hypothetical) minimal ones. 
We first show that the CR-Yamabe equation is always solvable, in a fixed neighborhood of $\mathcal{M}$, 
up to a Lagrange multiplier: see \cite{AM06} for a general reference on this method.

\begin{pro} \label{p:ls-reduction}
For $\var_\l$ as in \eqref{eq:var-l} there exists a unique $w_\l \in \mathfrak{S}^{1,2}(S^3)$, depending 
smoothly on $\l$, such that $\|w_\l\|_{\mathfrak{S}^{1,2}(S^3)} \leq C \, s$ and which satisfies 
\begin{equation}\label{eq:lagr}
\int_{S^3} \var_{\l}^2 \frac{\pa \var_{\l}}{\pa \l} \, w_\l  \, 
   \hat{\th} \wedge d \hat{\th} = 0; \qquad \qquad  L_s (\var_{\l} + w_\l)  - 2 
   (\var_{\l} + w_\l)^3 
   = \ell \, \var_\l^2 \frac{\pa \var_{\l}}{\pa \l} 
\end{equation}
for some $\ell \in \R$.  Moreover, there exists $\delta > 0$ with the following property: if 
there exists a critical point of $Q_{(s)}$ in a $\delta$-neighborhood of $\mathcal{M}$ (in $\mathfrak{S}^{1,2}$ norm), 
then it must  be of the form $\var_{\l} + w_\l$ up to a homogeneous action on $S^3$ and up to 
a scalar multiple, with $w_\l$ as above. 
\end{pro}

\begin{pf} 
For $\l > 1$, $\var_\l$ has a global maximum at $(z_1, z_2) = (0,1)$. Locally near these functions, 
all other extremals can be obtained composing on the right with elements of $SU(2)$. 
When also $\l$ varies, the extremals can be described locally near the $\var_\l$'s by 
$$
  \varSigma_{\L ,\g} = \left\{ \var_{{\bf a},\l}(z_1, z_2) := \var_\l (U_{\bf a}(z_1,z_2))\; : \;  
   {\bf a} \in  (-\g,\g)^3, \l \in [1/2,2\L]   \right\} \subseteq \mathcal{M}, 
$$
where 
\begin{equation}\label{eq:Ua}
U_{\bf a}(z_1,z_2) = \left( \exp \left( \begin{matrix}
      0 & a_1 + i \, a_2 \\ - a_1 + i \, a_2 & i \, a_3 
      \end{matrix} \right) \right)  \left(\begin{array}{c}
  z_1 \\ z_2 \\ 
      \end{array}
      \right), \qquad \quad {\bf a} = (a_1, a_2, a_3). 
\end{equation}

Consider next the CR-Yamabe equation on the standard sphere 
$$
L_0 u = 2 u^3 \qquad \qquad \hbox{ on } S^3. 
$$
It was proved in \cite{MuUg} (see Lemma 5 there) that solutions of the linearized equation at $\var_\l$
$$
  L_0 v = 6 \var_\l^2  v \qquad \qquad \hbox{ on } S^3
$$
are of the form 
$$
  v = l_0 \frac{\pa \var_{{\bf a},\l}}{\pa \l} + \sum_{i=1}^3 l_i \frac{\pa \var_{{\bf a},\l}}{\pa a_i},  
$$
where $l_i \in \R$ and where the latter derivatives are evaluated at ${\bf a} = 0$. 

\medskip

Define $\widetilde{W} = \widetilde{W}_\l$  to be the space of  functions $\tilde{w}$ satisfying the four constraints
\begin{equation}\label{eq:constr-W}
 \int_{S^3} \varphi_{{\bf a},\l}^2 \frac{\pa \varphi_{{\bf a},\l}}{\pa \l} \, \tilde{w}  \, 
   \hat{\th} \wedge d \hat{\th} = 0; \qquad \qquad 
   \int_{S^3} \varphi_{{\bf a},\l}^2 \frac{\pa \varphi_{{\bf a},\l}}{\pa a_i} \, \tilde{w}  \, 
     \hat{\th} \wedge d \hat{\th} = 0;  \quad i = 1, 2, 3, 
\end{equation}
where, again, the derivatives are evaluated at ${\bf a} = 0$.

It follows from the  classification result in \cite{MuUg} 
and Fredholm's theory that the operator 
$$
 A_{\widetilde{W}} : \tilde{w} \mapsto P_{\widetilde{W}} \left[   L_s  \tilde{w} - 6 \var_{\l}^2 
 \tilde{w}\right], 
$$
where $P_{\widetilde{W}}$ denotes the projection onto $\widetilde{W}$, is invertible from $\widetilde{W}$ 
in itself. 

Setting $S_{{\bf a},\l}(\tilde{w}) := L_s (\var_{\l} + \tilde{w}) - 2 
   (\var_{{\bf a},\l} + \tilde{w})^3$, equation \eqref{eq:lagr} becomes $P_{\widetilde{W}} S_{{\bf a},\l}(\tilde{w})  =0$. 
Since $A_{\widetilde{W}}$ is invertible (with the norm of the inverse uniformly bounded), we have that 
$$
  P_{\widetilde{W}} S_{{\bf a},\l}(\tilde{w})  =0 \qquad \Longleftrightarrow \qquad \tilde{w} = T_{{\bf a},\l}(\tilde{w}), 
$$  
where 
$$
 T_{{\bf a},\l}(\tilde{w})  = 
  - (A_{\widetilde{W}})^{-1} \left\{S_{{\bf a},\l}(0) - 2 \left[  
     (\var_{{\bf a},\l} + \tilde{w})^3  - \var_{{\bf a},\l}^3 - 3  \var_{{\bf a},\l}^2 \tilde{w}  \right] \right\}. 
$$
From the smoothness in $s$ of the $J_{(s)}$ structures it follows that $\| T_{{\bf a},\l}(0) \| = O(s)$. 
Moreover, it is quite standard that for $s$ and $\d$ small 
$$
   \| T_{{\bf a},\l}(\tilde{w}_1)  - T_{{\bf a},\l}(\tilde{w}_2) \| 
  = o(1)\|  \tilde{w}_1 - \tilde{w}_2 \|, \quad \| \tilde{w}_1 \|, \| \tilde{w}_2 \| \leq \d. 
$$
It follows that for $s$ small $T_{{\bf a},\l}$ is a contraction in a normed ball of 
radius $C \, s$ for $C > 0$ large and fixed, so in such a ball there exists a unique fixed point $w_\lambda$ of $T_{{\bf a},\l}$. 

In this way we found a (unique) solution to the problem 
$$
L_s (\var_{\l} + w_\l)  - 2 
   (\var_{\l} + w_\l)^3 
   = \ell \, \var_\l^2 \frac{\pa \var_{\l}}{\pa \l} + \sum_{i=1}^3 \ell_i \var_\l^2 \frac{\pa \var_{{\bf a},\l}}{\pa a_i}|_{{\bf a} = 0} 
$$
for some Lagrange multipliers $\ell, \ell_i$. However the last three vanish by {\em Palais' criticality principle}. In fact, 
let us recall that, being $(S^3, J_{(s)})$ a homogeneous space, $Q_{(s)}$ is invariant under the maps 
$U_{\bf a}$ as in \eqref{eq:Ua}. Therefore, with obvious notation, we have with the same Lagrange multipliers that 
$$
L_s (\var_{\l,{\bf a}} + w_{\l,{\bf a}})  - 2 
   (\var_{\l,{\bf a}} + w_{\l,{\bf a}})^3 
   = \ell \, \var_{\l,{\bf a}}^2 \frac{\pa \var_{\l,{\bf a}}}{\pa \l} + \sum_{i=1}^3 \ell_i \var_{\l,{\bf a}}^2 \frac{\pa \var_{{\bf a},\l}}{\pa a_i}, 
$$
for ${\bf a}$ in a neighborhood of zero. Differentiating with respect to $a_i$ and then scalar-multiplying 
by $\frac{\pa \var_{{\bf a},\l}}{\pa a_j}$ one obtains an invertible system for $(\ell_i)_i$, yielding that 
$\ell_i = 0$ for $i = 1, 2, 3$, as desired.

Let now $u$ be a critical point of $Q_{(s)}$ in a $\delta$-neighborhood of $\mathcal{M}$ for $s$ small. Then it 
satisfies $L_s u = \mu \, u^3$ for some Lagrange multiplier $\mu$. Since $u$ is close to the family of 
$\var_{\l}$'s, satisfying $L_0 \var_{\l} = 2 \var_{\l}^3$, the multiplier $\mu$ must be $\delta$-close to $2$. 

Defining $\tilde{u} = \mu^{- \frac 12} u$, this is still close of order $\d$ to $\mathcal{M}$, and 
it satisfies $L_s \tilde{u} = 2 \tilde{u}^3$, i.e. the second equation in \eqref{eq:lagr} with $\ell = 0$. 
By uniqueness of the fixed point, we must then have $\tilde{u} = \var_{\l} + w_\l$, up to a homogeneous 
action on $S^3$. 
This concludes the proof.  
\end{pf}

\begin{rem}\label{r:contraction-better}
In Proposition \ref{p:ls-reduction} it is possible to replace the $\var_\l$'s with other approximate 
solutions to the CR-Yamabe equation on Rossi spheres. With a better approximate solution, for example, 
one would then require a correction as in \eqref{eq:lagr} of smaller norm, yielding a more precise 
expansion for the quotient $Q_{(s)}$. This observation will be crucially used in the next two sections. 
\end{rem}

\medskip

\subsection{Expansion of the CR Sobolev quotient}

Recalling the latter statement in Proposition \ref{p:ls-reduction}, we analyze the 
CR Sobolev quotient on functions of the form $\var_\l + w_\l$, showing that 
it is strictly higher than the standard spherical one. We first show that the latter expansion 
is always even in $s$.

\begin{lem}
Let $s > 0$ be small, and let $w_\l^{(s)}$ and $w_\l^{(-s)}$ denote the counterparts of $w_\l$ in Proposition \ref{p:ls-reduction} 
for $s$ and $-s$ respectively. Then one has that 
$$
Q_{(s)}(\var_\l + w_\l^{(s)}) = Q_{(-s)}(\var_\l + w_\l^{(-s)}). 
$$ 
\end{lem}

\begin{pf}
Let $\iota : S^3 \to S^3$  be the diffeomorphism given in \eqref{eq:iota}. We notice that $\var_\l$ is invariant under $\iota$ 
and that, due to \eqref{B8}, \eqref{S1} and \eqref{S2},  for any $u \in \mathfrak{S}^{1,2}(S^3)$ one has 
$$
  Q_{(s)} (\iota^* u) = Q_{(-s)} (u). 
$$ 
From this covariance property and the uniqueness in Proposition \ref{p:ls-reduction} it follows that 
$w_\l^{(s)} = \iota^* w_\l^{(-s)}$, and therefore we get  
$$
  Q_{(s)} (\var_\l + w_\l^{(s)}) = Q_{(s)} (\iota^* (\var_\l + w_\l^{(-s)})) = Q_{(-s)}(\var_\l + w_\l^{(-s)}), 
$$
which is the desired conclusion. 
\end{pf}

\

We analyse next two situations. The first is when the parameter $\l$ in the previous lemma 
tends to infinity or to zero, and the second when $\log \l$ remains bounded. In the latter case we will show that the 
CR Sobolev quotient would be strictly  higher than $\mathcal{Y}(S^3, J_{S^3})$, which would give 
a contradiction to \eqref{eq:ineq-quotient-s}. On the other hand, we can also rule out the former case
using  the estimates on the Green's function in Section \ref{s:Green}, and in particular the 
negativity of the mass of $(S^3, J_{(s)})$ for $s$ small and non zero. The proofs of the next two results, beginning from the latter case, 
are given in the next two appendices.

\begin{pro}\label{p:expansion}
Let $\L > 1$ be a fixed number. Then there exist $C_\L > 0$ 
such that, for $\l \in [1/\L,\L]$ and for $s$ 
small one has  $Q_{(s)}(\var_\l + w_\l) =  4 \pi  + 
   s^2 \mathcal{A}_\l + \mathcal{B}_{\l,s}$, where 
\begin{equation} \nonumber 
   \mathcal{A}_\l  = 
  \frac{16 \pi  \l^2 (3 + 12 \l^2 + 2 \l^4 + 12 \l^6 + 3 \l^8)}{(1 + \l^2)^6}, 
\end{equation}
and where $|\mathcal{B}_{\l,s}| \leq C_\L s^3$. 
\end{pro}

\medskip

\begin{pro}\label{p:expansion-mass}
The following expansion holds true, uniformly in $s$ (small) 
\begin{equation}\label{eq:eq-in-prop-exp-0} \nonumber 
   Q_{(s)}(\var_\l + w_\l) =  4 \pi  - \frac{8}{3 }  \frac{m_s}{\l^2} + O(\frac{s^2}{\l^3}) = 
   4 \pi  + 48 \pi \frac{s^2}{\l^2} (1 + o_s(1)) + O(\frac{s^2}{\l^3}), 
\end{equation}
for $\l$ large. 
\end{pro}

\medskip

\begin{rem}
The above function  $\l \mapsto \mathcal{A}_\l$ is positive and strictly decreasing 
for $\l > 1$, see the picture below. Notice that the matching of the first-order 
correction terms for $\l$ large in the above two propositions: the expansions are 
indeed obtained with \underline{two completely different approaches}. However, while the {\em mass} does not appear in 
the expansions of Section \ref{s:pf-prop}, it is somehow {\em hidden} in the fact that there we are 
using standard coordinates on $S^3$, and not CR normal coordinates. 
    
\begin{figure}[h]\label{fig:graph} \begin{center}
 \includegraphics[angle=0,width=9.0cm]{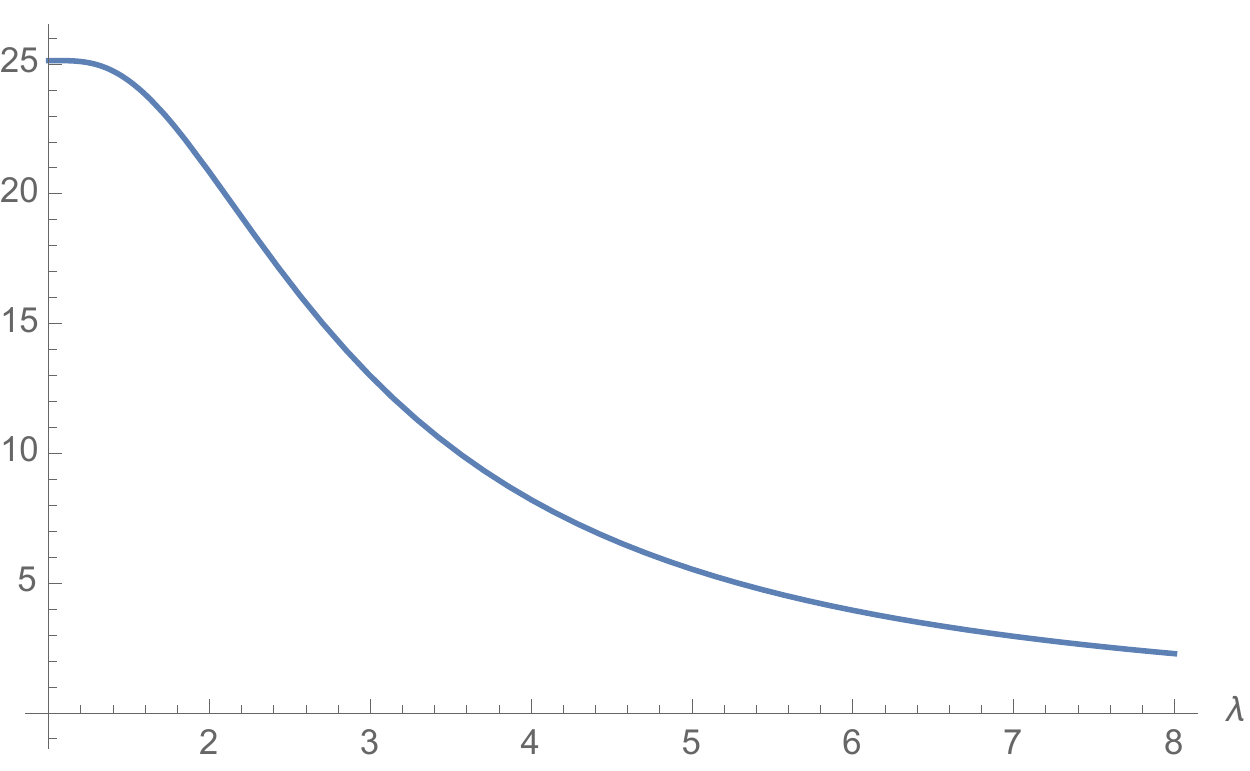} \end{center} 
\end{figure}
\end{rem}

We can finally prove our second main result. 

\

\begin{pfn} {\sc of Theorem \ref{t:SQ}.}
Assume by contradiction that $u$ is a minimizer of the CR-Sobolev quotient 
$Q_{(s)}$ for $s \neq 0$ small. By Lemma \ref{l:closeness-S12}, $u$ must then 
lie in a $\delta$-neighborhood of the manifold $\mathcal{M}$ defined in \eqref{eq:mathcalM}. 
From the 
second part of Proposition \ref{p:ls-reduction} we have also that 
$u = \var_\l + w_\l$ up to a homogeneous action on $S^3$, where $w_\l$ is as in the first part of the Proposition. The conclusion 
then follows from Proposition \ref{p:expansion} and Proposition \ref{p:expansion-mass}, 
which cover all ranges of $\lambda$ for $s$ small enough. 
\end{pfn}

\

\section{Appendix A: proof of Proposition \ref{p:expansion}}\label{s:pf-prop}

We consider the \emph{Cayley map} from $S^3$ into $\H^1$ given by 
\begin{equation}\label{eq:cayley}
  \mathcal{F}(z_1, z_2) = \left( \frac{z_1}{1+z_2}, {\text Re} \left( i \frac{1-z_2}{1+z_2}  \right)  \right), 
\end{equation}
with inverse 
$$
  \mathcal{F}^{-1}(z,t) = \left( \frac{2iz}{t + i \left(1 + |z|^2 \right)}, \frac{-t+i(1-|z|^2)}{t+i(1+|z|^2)} \right). 
$$
Using $\mathcal{F}$, we can derive explicit expressions for the CR maps on $S^3$. Letting 
$\mathfrak{d}_\l$ denote the natural dilation in the Heisenberg group 
$$
  \mathfrak{d}_\l(z,t) = (\l \, z, \l^2 \, t); \qquad \qquad \l > 0, 
$$
consider the map $\Phi_\l : S^3 \to S^3$ defined by 
$$
  \Phi_\l(p) = \left( \mathcal{F}^{-1} \circ \mathfrak{d}_\l \circ \mathcal{F} \right) (p). 
$$
By explicit computations one finds that the inverse  is given by 
\begin{equation}\label{eq:inverse-mobius}
  \Phi_\l^{-1}(z_1,z_2) = \left( \frac{2 \l (\ov{z}_2+1) z_1}{\l^2 |z_2+1|^2+\ov{z}_2+ |z_1|^2 -z_2} , 
  \frac{\l^2 |z_2+1|^2-\ov{z}_2- |z_1|^2 +z_2}{\l^2 |z_2+1|^2+\ov{z}_2+ |z_1|^2-z_2}  \right). 
\end{equation}
For later purposes the following formula will be useful 
\begin{equation}\label{eq:inverse-Mobius-bubble}
\var_\l(\Phi_\l^{-1}(z_1,z_2))^{-3} = \frac{1}{2} \l^{-1} \left( \frac{|1+z_2|^2}{\left( \l ^2 
|1+z_2|^2 + |z_1|^2 \right)^2 - (z_2 - \ov{z}_2)^2} \right)^{\frac{1}{2}}.
\end{equation}
Notice also that $\var_{\l=1} \equiv 1$ on $S^3$.

\medskip

\subsection{Approximate solutions}\label{ss:finding-w}

We construct next, on every compact interval in the range of $\lambda$, 
approximate solutions to the CR-Yamabe equation with $s \neq 0$  up 
an order $O(s^2)$, improving the accuracy of the $\var_{\l}$'s (approximate up to 
order $O(s)$) for $s \neq 0$.

\begin{lem}\label{l:expansion}
Let $\L > 1$ be a fixed number. Then there exist $C_\L > 0$ and regular functions 
$\hat{w}_\l$, depending smoothly on $\l$ such that for $\l \in [1/\L,\L]$ and for $s$ 
small one has  
$$
   L_s (\var_\l + s \hat{w}_\l) - 2  (\var_\l +s  \hat{w}_\l)^3 
   = f_\l, 
$$
with $\|f_\l\|_{L^\infty(S^3)} \leq C_\L s^2$.
\end{lem}

\begin{pf}
%
%
%
%
%
%
%
%
Recall that the extremals of the CR-Sobolev inequality (up to a homogeneous CR-action of $S^3$) have the expression 
in \eqref{eq:var-l}, namely 
$$
  \var_\l =  \l \left(  \frac{\left(  |z_1|^2 + |z_2+1|^2 \right)^2 -  (z_2 - \ov{z}_2)^2}{\left( \l ^2 |z_1|^2 + |z_2+1|^2 \right)^2 
  - \l ^4 (z_2 - \ov{z}_2)^2} \right)^{\frac{1}{2}} = 2 \l \left(  \frac{|1+z_2|^2}{\left( \l ^2 |z_1|^2 + |z_2+1|^2 \right)^2 
    - \l ^4 (z_2 - \ov{z}_2)^2} \right)^{\frac{1}{2}}, 
$$
and for all $\l > 0$ they satisfy the equation 
\begin{equation}\label{eq:eq-var-l}
  L_0 \var_\l = - 4 \Delta_b \var_\l + 2 \var_\l = 2 \var_\l^3 \qquad \qquad \hbox{ on } \quad S^3. 
\end{equation}
Our goal is to find a correction $s \hat{w}_\l$ such that $\var_\l + s \hat{w}_\l$ satisfies the 
CR-Yamabe equation on $(S^3,J_{(s)})$ up to an order $s^2$. 
Recalling \eqref{eq:dot-L-ddot-L} and \eqref{eq:dot-R-ddot-R}, it is sufficient  to solve for 
 $$
 - 4 \Delta_b \hat{w}_\l + 2 \hat{w}_\l - 6 \var_\l^2 \hat{w}_\l = \mathcal{G}_\l := 8 Z_1 Z_1 \var_\l +\hbox{conj.}.
$$
From a straightforward computation one has that 
\begin{eqnarray} 
\mathcal{G}_\l(z_1,z_2)  \nonumber  & = & 
 \frac{192 \l (\l ^2-1)^2 |1 + z_2|^8 {\text Re} \left[ z_1^2 \left(1 + z_2 + \l ^2 (z_2-1) \right)^2 \right]}{\left[ \left( \l ^2 |z_1|^2 + |z_2+1|^2 \right)^2 
  - \l ^4 (z_2 - \ov{z}_2)^2 \right] ^{\frac 52} \left[ \left(  |z_1|^2 + |z_2+1|^2 \right)^2 
    -  (z_2 - \ov{z}_2)^2 \right] ^{\frac 32}} \\ & = & \nonumber
     \frac 34 \l^{-4} \frac{(\l ^2-1)^2 {\text Re}  \left[ z_1^2 \left(1 + z_2 + \l ^2 (z_2-1) \right)^2 \right]}{\left[ \left( \l ^2 |z_1|^2 + |z_2+1|^2 \right)^2 
      - \l ^4 (z_2 - \ov{z}_2)^2 \right] ^{\frac 52}} \var_\l(z_1,z_2)^5. 
\end{eqnarray}
It is useful  to evaluate this expression after composing with the inverse CR 
map defined in \eqref{eq:inverse-mobius}: by direct computation, using also \eqref{eq:inverse-mobius}, one finds that 
\begin{eqnarray} \label{eq:inv-inv}\nonumber
\mathcal{G}_\l(\Phi_\l^{-1}(z_1,z_2)) & = &  \frac{3}{4} \l^{-3}(\l^2-1)^2 \left( \frac{|1+z_2|^2}{\left( \l ^2 |1+z_2|^2 + |z_1|^2 \right)^2 - (z_2 - \ov{z}_2)^2} 
  \right)^{\frac{3}{2}} \\ & \times & 
  \left[ z_1^2 \left( 1 - \ov{z}_2 + \l ^2 (1 + \ov{z}_2) \right)^4 + \ov{z}_1^2 \left( 1 - z_2 + \l ^2 (1 + z_2) \right)^4  \right].
\end{eqnarray}
Let us recall the  covariance of the conformal sub-Laplacian $L_\theta$: for a conformal contact form 
$\tilde{\theta} = u^{\frac{4}{Q-2}} \theta$ one has 
$$
  L_{\tilde{\theta}} \var = u^{-\frac{Q+2}{Q-2}} L_\theta (u \, \cdot). 
$$
Let $\mathfrak{L}_{\var_\l}$ be the 
linearized CR-Yamabe operator at $\var_\l$ on $(S^3, J_{(0)})$, i.e. 
\begin{equation}\label{eq:linearized-operator}
  \mathfrak{L}_{\var_\l} v = - 4 \D_b v + 2 v - 6 \var_\l^2 v, 
\end{equation}
and let $\mathfrak{w}_\l$ denote the pull-back of  $ \hat{w}_\l$ 
via $\Phi_\l$, namely 
\begin{equation}\label{eq:frak-w}
  \mathfrak{w}_\l(z) = \var_\l^{-1}(\Phi_\l^{-1}(z)) \,  \hat{w}_\l (\Phi_\l^{-1}(z)). 
\end{equation}
Then the  covariance of $L_\th$  implies  that 
\begin{equation}\label{eq:covariance-lin}
  (\mathfrak{L}_{\var_1 \equiv 1} \mathfrak{w}_\l)(x) =  \var_\l(\Phi_\l^{-1}(x))^{-3}   (\mathfrak{L}_{\var_\l} \hat{w}_\l) (\Phi_\l^{-1}(x)). 
\end{equation}
It follows from this formula and \eqref{eq:inv-inv} that the pull-back $\mathfrak{w}_\l$ satisfies the following equation on $S^3$, 
which has constant coefficients on the left-hand side   
\begin{equation}\label{eq:const-coeff}
- 4 \Delta_b \mathfrak{w}_\l - 4 \mathfrak{w}_\l = 
   12 (\l^2-1)^2 {\text Re}   \frac{(1 - \ov{z}_2 + \l^2 (1+ \ov{z}_2)) z_1^2}{(1 - z_2 + \l^2 (1+ z_2))^3}. 
\end{equation}
The latter equation can be solved explicitly in  $\mathfrak{w}_\l$ via Fourier decomposition: 
%
%
in fact, the right-hand side in \eqref{eq:const-coeff} is given by 
$$
  12 \frac{(\l^2-1)^2}{(\l^2+1)^2} {\text Re}  \frac{z_1^2 (1 - \Gamma \ov{z}_2)}{(1 - \Gamma z_2)^3}; 
  \qquad \qquad \hbox{ with } \quad \Gamma = \frac{1 - \l^2}{1 + \l^2}.  
$$
Since we have the expansion 
$$
  \frac{1}{(1 - \Gamma z_2)^3} = 1 + 3 \Gamma z_2 + 6 \Gamma^2 z_2^2 + 10 \Gamma^3 z_2^3 + 15 \Gamma^4 z_2^4 + \cdots,  
$$
we obtain that 
\begin{eqnarray} \nonumber 
12 \frac{(\l^2-1)^2}{(\l^2+1)^2} {\text Re}  \frac{z_1^2 (1 - \Gamma \ov{z}_2)}{(1 - \Gamma z_2)^3}
   & = & 
   12 \frac{(\l^2-1)^2}{(\l^2+1)^2} {\text Re}   \left\{ z_1^2 \left[ 1 + 3 \Gamma z_2 + 6 \Gamma^2 z_2^2 + 10 \Gamma^3 z_2^3 + 15 \Gamma^4 z_2^4 + \cdots \right. \right. \\  & - & \left. \left. \Gamma \ov{z}_2 \left( 1 + 3 \Gamma z_2 + 6 \Gamma^2 z_2^2 + 10 \Gamma^3 z_2^3 + 15 \Gamma^4 z_2^4 + \cdots \right) \right] \right\}.
\end{eqnarray}
While the set of functions of the right--hand side in the first line are spherical harmonics, 
i.e. satisfying 
\begin{equation}\label{eq:spher-simple}
 - \Delta_b (z_1^2 z_2^k) = (k+2) z_1^2 z_2^k, 
\end{equation}
the functions on the second line of the right-hand side are not. However, they can be easily modified in order to satisfy an eigenvalue equation. 
More precisely, one has that (see \cite{JLCR}) 
\begin{equation}\label{eq:spher-tricky}
 - \Delta_b \left( z_2^k \left( z_2 \ov{z}_2 - \frac{k+1}{k+4} \right) \right) = (10+3k) 
 \left( z_2^k \left( z_2 \ov{z}_2 - \frac{k+1}{k+4} \right) \right). 
\end{equation}
Hence we rewrite the  right-hand side in \eqref{eq:const-coeff} in the following way 
\begin{eqnarray*}
  & &   12 \frac{(\l^2-1)^2}{(\l^2+1)^2} {\text Re} \frac{z_1^2 (1 - \Gamma \ov{z}_2)}{(1 - \Gamma z_2)^3} \\ & = &  12 \frac{(\l^2-1)^2}{(\l^2+1)^2} {\text Re}  \left\{ z_1^2 \left[ 1 + 3 \Gamma z_2 + 6 \Gamma^2 z_2^2 + 10 \Gamma^3 z_2^3 + 15 \Gamma^4 z_2^4 + \cdots  
 \right. \right. \\ 
 & - & \Gamma \ov{z}_2 - 3 \Gamma^2 \left( z_2 \ov{z}_2 - \frac{1}{4} \right) - 6 \Gamma^3 z_2 \left( z_2 \ov{z}_2 - \frac{2}{5} \right) 
  - 10 \Gamma^4 z_2^2 \left( z_2 \ov{z}_2 - \frac{3}{6} \right) + \cdots \\ 
  & - & \left. \left. 3 \Gamma^2 \frac{1}{4} - 6 \Gamma^3 z_2 \frac{2}{5} - 10 \Gamma^4 z_2^2 \frac{3}{6} - \cdots \right] \right\}.  
\end{eqnarray*}
The latter expression can in turn be rewritten as 
\begin{eqnarray} \label{eq:rewritten} \nonumber
 & & 12 \Gamma^2 {\text Re}  \left\{ z_1^2 \left[ 
  \sum_{k=0}^\infty \frac{(k+1)(k+2)}{2} (\Gamma z_2)^k \left( 1 - \Gamma^2 \frac{k+3}{k+4}\right)  \right. \right. \\ 
  & - & \left. \left. \Gamma^2 \sum_{k=-1}^\infty \frac{(k+2)(k+3)}{2} (\Gamma z_2)^k  \left( z_2 \ov{z}_2 - \frac{k+1}{k+4} \right) \right] \right\}.
\end{eqnarray}
Recall that by \eqref{eq:const-coeff}, to obtain  $\mathfrak{w}_\l$, we need to invert the operator $- 4 \Delta_b - 4$ on the latter expression, so we have to 
divide the coefficients of the spherical harmonics respectively by (using \eqref{eq:spher-simple} and \eqref{eq:spher-tricky}) 
$4 (k+2) - 4 = 4(k+1)$ and by $4 (3k + 10) - 4 = 12(k+3)$. We then find 
\begin{equation}\label{eq:w-series}
\mathfrak{w}_\l =  \frac{3}{2} \Gamma^2 {\text Re}  \left\{ z_1^2 \left[ 
  \sum_{k=0}^\infty (k+2) (\Gamma z_2)^k \left( 1 - \Gamma^2 \frac{k+3}{k+4}\right)  
  - \frac 13 \Gamma^2 \sum_{k=-1}^\infty (k+2) (\Gamma z_2)^k  \left( z_2 \ov{z}_2 - \frac{k+1}{k+4} \right) \right] \right\},
\end{equation}
with $\Gamma = \frac{1-\l^2}{1+\l^2}$. Notice that since $|\Gamma| < 1$ all the above series are absolutely converging 
on $S^3$. Finally, the correction $\hat{w}_\l$ to $\var_\l$ for the CR-Yamabe equation can be obtained from \eqref{eq:frak-w}. 
\end{pf}

\medskip

\subsection{Second order expansion of the CR Sobolev quotient}

We want next to analyse  the order $s^2$ in the expansion of the CR Sobolev quotient.

\begin{lem}\label{p:expansion-2}
If $Q_{(s)}$ is as in \eqref{eq:Qs}, then we have that 
\begin{eqnarray}\label{eq:eq-in-prop-exp} \nonumber 
  Q_{(s)}(\var_\l + s \hat{w}_\l)
  =  4 \pi  + 
  \frac{16 \pi  \l^2 (3 + 12 \l^2 + 2 \l^4 + 12 \l^6 + 3 \l^8)}{(1 + \l^2)^6} s^2 + \mathcal{B}_{\l,s}, 
\end{eqnarray}
with $|\mathcal{B}_{\l,s}| \leq C_\L s^3$. 
\end{lem}

\begin{pf}
Recall that,  at $s = 0$, from \eqref{eq:dot-R-ddot-R} one has $\frac{d}{ds} R_s = 0 $ and $\frac{d^2}{d s^2} R_s = 8$. We use the 
choice of contact form 
$$
\hat{\theta} = \frac{1}{2} i \sum_{i=1}^2 \left( z_i d \ov{z}_i - \ov{z}_i d z_i \right); \qquad \qquad 
\hat{\theta} \wedge d \hat{\theta} = 2 \, d \s_{Eucl}. 
$$
From the expression of $- \ddot{\Delta}_b$ (in \eqref{eq:dot-L-ddot-L}) and of $\ddot{R}$ we have that  
the second derivative $\ddot{Q}(\var_\l)$ of $Q_{(s)}(\var_{\l})$ at $s = 0$ 
is given by 
$$
  \ddot{Q}(\var_\l) = \int_{S^3} \var_\l \left( - 4 \ddot{\Delta}_b  \var_\l  + \ddot{R}  \var_\l  \right) \hat{\theta} \wedge d \hat{\theta} 
  = \int_{S^3} \var_\l \left( - 16 \Delta_b \var_\l  + 8  \var_\l  \right) \hat{\theta} \wedge d \hat{\theta}. 
$$ 
Using \eqref{eq:eq-var-l}, this also becomes 
\begin{equation}\label{eq:ddot-Q}
 \ddot{Q}(\var_\l) =  8 \int_{S^3} \var_\l^4 \, \hat{\theta} \wedge d \hat{\theta} = 32 \pi^2, 
\end{equation}
since the integral is independent of $\l$ and since $\var_{\l=1} \equiv 1$.

\

Our next goal is to expand to second order in $s$ the quantity $Q_{(s)}(\var_\l+s \, \hat{w}_\l)$. We claim that 
\begin{eqnarray}\label{eq:claim-2} \nonumber
Q_{(s)}(\var_\l+s \, \hat{w}_\l) & = & 4 \pi + \frac{s^2}{2 \pi} \left( \frac 12 \ddot{Q}(\var_\l) -  \int_{S^3} \hat{w}_\l \mathfrak{L}_{\var_\l} \hat{w}_\l \, \hat{\theta} \wedge d \hat{\theta} \right) + o(s^2) \\ 
  & = & 4 \pi + \frac{s^2}{2 \pi}  \left( \frac 12 \ddot{Q}(\var_\l) -  \int_{S^3} \mathfrak{w}_\l \mathfrak{L}_{\var_1} \mathfrak{w}_\l \, \hat{\theta} \wedge d \hat{\theta} \right) + o(s^2). 
\end{eqnarray}
Here,  $\mathfrak{w}_\l$ is given in \eqref{eq:w-series} (see also \eqref{eq:frak-w}) 
and $\mathfrak{L}_{\var_\l}$ is given in \eqref{eq:linearized-operator}. The latter equality 
follows from the covariance property \eqref{eq:covariance-lin}. 
To check this claim, we want to expand $Q_{(s)}(\var_\l+s \, \hat{w}_\l)$, which we write as 
$$
  \frac{\int_{S^3} (\var_\l + s \hat{w}_\l) \left( L_0+ s \dot{L}  + \frac{1}{2} s^2 \ddot{L} \right) (\var_\l + s \hat{w}_\l)\hat{\theta} \wedge d \hat{\theta}}{\left(  \int_{S^3} (\var_\l+ s \hat{w}_\l)^4 \hat{\theta} \wedge d \hat{\theta} \right)^{\frac 12}}. 
$$
Expanding in $s$ we find that this quantity is equal to 
\begin{eqnarray*}
  & & Q_{(s)}(\var_\l+s \hat{w}_\l) \\ & = & \frac{\int_{S^3} \left[ \var_\l L_0 \var_\l + s \left( \var_\l \dot{L} \var_\l + 2 \hat{w}_\l L_0 \var_\l \right) + s^2 
  \left( \frac{1}{2} \var_\l \ddot{L} \var_\l + 2 \hat{w}_\l \dot{L} \var_\l + \hat{w}_\l L_0 \hat{w}_\l \right)  \right] \hat{\theta} \wedge d \hat{\theta}}{\left(  \int_{S^3} (\var_\l^4 
  + 4 s \var_\l^3 \hat{w}_\l + 6 s^2 \var_\l^2 \hat{w}_\l^2) \, \hat{\theta} \wedge d \hat{\theta} \right)^{\frac 12}} \\ & +& o(s^2). 
\end{eqnarray*}
The first-order term in $s$ vanishes, as one can see using the Euler equation for $\var_\l$, so we will just consider  the 
second-order term. Since $\mathfrak{w}_\l$ only consists of spherical harmonics of positive order, 
see \eqref{eq:w-series}, using \eqref{eq:frak-w} it also turns out that 
$$
  \int_{S^3} \var_\l^3 \hat{w}_\l \hat{\theta} \wedge d \hat{\theta}  = \int_{S^3} \mathfrak{w}_\l \hat{\theta} \wedge d \hat{\theta} = 0, 
$$ 
so there is no contribution to the expansion 
of the denominator from the first-order term (in $s$) in the denominator. 

Since $\int_{S^3} \var_\l L_0 \var_\l \, \hat{\theta} \wedge d \hat{\theta} = 8 \pi^2$ and  $\int_{S^3} \var_\l^4 \, \hat{\theta} \wedge d \hat{\theta} = 4 \pi^2$, we can collect these numbers in the numerator and denominator respectively  to get that 
$$
  Q_{(s)}(\var_\l+s \hat{w}_\l) = \frac{8 \pi^2}{(4 \pi^2)^{1/2}} \frac{1 + \frac{s^2}{8 \pi^2} \int_{S^3} ( \frac{1}{2} \var_\l \ddot{L} \var_\l 
  + 2 \hat{w}_\l \dot{L} \var_\l + \hat{w}_\l L_0 \hat{w}_\l) \, \hat{\theta} \wedge d \hat{\theta} }{\left( 1 + \frac{s^2}{4 \pi^2} \int_{S^3} 6 \var_\l^2 \hat{w}_\l^2 \, 
  \hat{\theta} \wedge d \hat{\theta} \right)^{\frac 12}} + o(s^2). 
$$
Taylor-expanding one finds 
\begin{eqnarray*}
& & Q_{(s)}(\var_\l+s \hat{w}_\l) \\ & = & \frac{8 \pi^2}{(4 \pi^2)^{1/2}} 
  \left[ 1 + \frac{s^2}{8 \pi^2}  \left( \int_{S^3} ( \frac{1}{2} \var_\l \ddot{L} \var_\l 
    + 2 \hat{w}_\l \dot{L} \var_\l + \hat{w}_\l L_0 \hat{w}_\l ) \, \hat{\theta} \wedge d \hat{\theta}  -   \int_{S^3} 6 \var_\l^2 \hat{w}_\l^2\,  \hat{\theta} \wedge d \hat{\theta} \right) \right] \\ & + & o(s^2). 
\end{eqnarray*}
We now use the fact that $\hat{w}_\l$ satisfies 
$$
 \mathfrak{L}_{\var_\l} \hat{w}_\l := L_0 \hat{w}_\l - 6 \var_\l^2 \hat{w}_\l = - \dot{L} \var_\l
$$
to deduce that
\begin{eqnarray}\label{eq:2ndorder0} \nonumber 
Q_{(s)}(\var_\l+s \hat{w}_\l) & = & \frac{8 \pi^2}{(4 \pi^2)^{1/2}} 
  \left[ 1 + \frac{s^2}{8 \pi^2}  \left( \int_{S^3} ( \frac{1}{2} \var_\l \ddot{L} \var_\l
    -  \hat{w}_\l \mathfrak{L}_{\var_\l} \hat{w}_\l ) \, \hat{\theta} \wedge d \hat{\theta} \right) \right] + o(s^2)
    \\ & = & \frac{8 \pi^2}{(4 \pi^2)^{1/2}} 
      \left[ 1 + \frac{s^2}{8 \pi^2}  \left( \int_{S^3} ( \frac{1}{2} \var_\l \ddot{L} \var_\l
        -  \mathfrak{w}_\l \mathfrak{L}_{\var_1} \mathfrak{w}_\l) \, \hat{\theta} \wedge d \hat{\theta} \right) \right]
        + o(s^2). 
\end{eqnarray}

\

\noindent We next compute the latter integral. To explicitly integrate spherical harmonics, we  need the following explicit formula (see Proposition 5.3 in \cite{JLCR}) 
\begin{equation}\label{eq:JL-formula}
  \int_{S^3} |z_1|^4 |z_2|^{2k} \hat{\theta} \wedge d \hat{\theta} = \frac{8 \pi^2}{(k+1)(k+2)(k+3)}. 
\end{equation}
Both $\mathfrak{w}_\l$ and $\mathfrak{L}_{\var_1} \mathfrak{w}_\l$ consist of two types of spherical harmonics, orthogonal to each-other. 
For the first series, taking real parts, we need to compute integrals of the form (notice that only products of conjugate terms contribute)
$$
  \frac{1}{4} \int_{S^3} \left( z_1^2 z_2^k + \ov{z}_1^2 \ov{z}_2^k \right)^2 \hat{\theta} \wedge d \hat{\theta} =  \frac{4 \pi^2}{(k+1)(k+2)(k+3)}. 
$$
For the second series, still taking real parts, we need to compute instead 
\begin{eqnarray*}
 & & \frac{1}{4} \int_{S^3}  \left[ z_1^2 z_2^k \left( z_2 \ov{z}_2 - \frac{k+1}{k+4} \right)  + 
    \ov{z}_1^2 \ov{z}_2^k \left( z_2 \ov{z}_2 - \frac{k+1}{k+4} \right)   \right]^2 \hat{\theta} \wedge d \hat{\theta} 
    \\ & = & \frac{1}{2} \int_{S^3} |z_1|^4 |z_2|^{2k}  
    \left(  |z_2|^4 - 2 |z_2|^2 \frac{k+1}{k+4} + \left( \frac{k+1}{k+4} \right)^2 \right) \hat{\theta} \wedge d \hat{\theta}. 
\end{eqnarray*}
Using \eqref{eq:JL-formula}, the expression becomes 
$$
 \frac{12 \pi^2}{(k+2)(k+3)(k+4)^2(k+5)}. 
$$
Therefore,  from \eqref{eq:rewritten} and \eqref{eq:w-series} we obtain 
\begin{eqnarray*}
 -  \int_{S^3} \mathfrak{w}_\l \mathfrak{L}_{\var_1} \mathfrak{w}_\l \,  \hat{\theta} \wedge d \hat{\theta} & = & -  24 \frac{3}{2} \Gamma^4 \left\{ \sum_{k=0}^{\infty} 
 \frac{(k+1)(k+2)^2}{2} \left( 1 - \Gamma^2 \frac{k+3}{k+4}\right)^2 \Gamma^{2k} \frac{2 \pi^2}{(k+1)(k+2)(k+3)} \right. \\ 
 & + & \left. \sum_{k=-1}^{\infty} \frac{\Gamma^4}{6} (k+2)^2 (k+3) \Gamma^{2k}\frac{6 \pi^2}{(k+2)(k+3)(k+4)^2(k+5)}   \right\}.
\end{eqnarray*}
After some simplification,  this gives 
$$
  - \int_{S^3} \mathfrak{w}_\l \mathfrak{L}_{\var_1} \mathfrak{w}_\l \, \hat{\theta} \wedge d \hat{\theta} = - 36 \Gamma^4 \pi^2 \left\{    \sum_{k=0}^{\infty} 
   \frac{k+2}{k+3} \left( 1 - \Gamma^2 \frac{k+3}{k+4}\right)^2 \Gamma^{2k} + \sum_{k=-1}^{\infty} \Gamma^4
   \frac{\Gamma^{2k}(k+2)}{(k+4)^2(k+5)}  \right\}.
$$
Notice that the last series starts from $k = -1$, so after relabelling we get 
$$
  -  \int_{S^3} \mathfrak{w}_\l \mathfrak{L}_{\var_1} \mathfrak{w}_\l  \, \hat{\theta} \wedge d \hat{\theta} = - 36 \Gamma^4 \pi^2 \left\{    \sum_{k=0}^{\infty} 
   \frac{k+2}{k+3} \left( 1 - \Gamma^2 \frac{k+3}{k+4}\right)^2 \Gamma^{2k} + \sum_{k=0}^{\infty} \Gamma^2
   \frac{\Gamma^{2k}(k+1)}{(k+3)^2(k+4)}  \right\}.
$$
After some manipulation, the series reduces to a finite one, and we find 
$$
  - \int_{S^3} \mathfrak{w}_\l \mathfrak{L}_{\var_1} \mathfrak{w}_\l  \, \hat{\theta} \wedge d \hat{\theta} =  8 \pi^2 \Gamma^4 \left(\Gamma^2-3\right). 
$$
Collecting this formula and \eqref{eq:2ndorder0}, from \eqref{eq:ddot-Q} and \eqref{eq:claim-2}  we obtain the second order expansion 
\begin{eqnarray}\label{eq:2ndorder} \nonumber 
Q_{(s)}(\var_\l+s \, \hat{w}_\l) & = & 4 \pi + \frac{s^2}{2\pi} \left( \frac 12 \ddot{Q}(\var_\l) - \int_{S^3} \mathfrak{w}_\l \mathfrak{L}_{\var_1} \mathfrak{w}_\l \, \hat{\theta} \wedge d \hat{\theta} \right) + o(s^2)\\ & = & 4 \pi + 4 \pi  s^2   \left(\Gamma^6-3 \Gamma^4+2\right) + o(s^2), 
\qquad \quad \Gamma = \frac{1-\l^2}{1+\l^2}. 
\end{eqnarray}
This concludes the proof.   \end{pf}

\medskip

We display next the  graph of the function  $4 \pi \left(\Gamma^6-3 \Gamma^4+2\right)$ in $\Gamma$. This shows that the second-order correction of the Sobolev 
quotient is always positive in $\lambda$, and tends to zero as $\l \to \infty$. 

\begin{figure}[h]\label{fig:graph2} \begin{center}
 \includegraphics[angle=0,width=9.0cm]{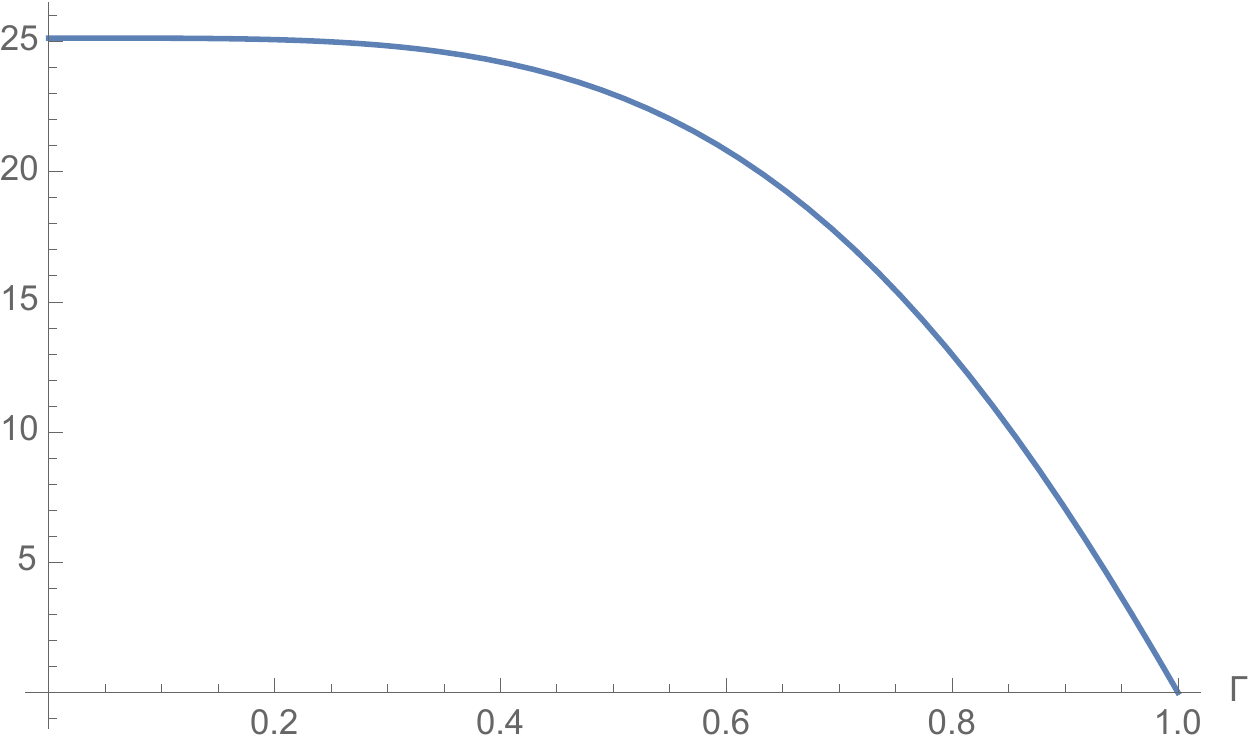} \end{center} 
\end{figure}

\medskip

\subsection{Conclusion} \label{ss:concl1}

We can use the observation in Remark 5.3, to work out the contraction argument in Proposition 
\ref{p:ls-reduction} starting from $\var_{\l} + s \hat{w}_\l$ instead of from $\var_{\l}$ only. 
Given the improved accuracy in Lemma \ref{l:expansion}, the contraction can be performed 
in a ball of radius $O(s^2)$ in $\mathfrak{S}^{1,2}(S^3)$, yielding a corresponding correction 
$\check{w}_\l$ of that order. By Lemma \ref{l:expansion} and the smoothness of $Q_{(s)}$, we then have that 
$$
  Q_{(s)} (\var_{\l} + s \hat{w}_\l + \check{w}_\l) = 
  Q_{(s)} (\var_{\l} + s  \hat{w}_\l)  + d Q_{(s)} (\var_{\l} + s \hat{w}_\l)  [\check{w}_\l] + O(\| \check{w}_\l \|^2) = 
  Q_{(s)} (\var_{\l} + s \hat{w}_\l) + O(s^4).  
$$
By uniqueness in the fixed point of the contraction, it must be $\var_{\l} + s \hat{w}_\l + \check{w}_\l = \var_{\l} + w_\l$, 
so the conclusion follows from Lemma \ref{p:expansion-2}.

\

\section{Appendix B: proof of Proposition \ref{p:expansion-mass}}\label{s:exp-mass}

The goal of this section is to expand $Q_{(s)}$ on the functions $\var_{\l} + w_\l$ given 
by Proposition \ref{p:ls-reduction} for large values of $\lambda$. Since the estimates 
of the previous section deteriorate for $\l$ in this range, we choose approximate solutions 
in terms of CR normal coordinates, better suited for highly-concentrated profiles.

Recall from the results in Section 5 of \cite{CMY17} that, given $p \in M$,  the Green's function of the conformal sub-Laplacian 
satisfies, in CR normal coordinates 
\begin{equation}\begin{split} \label{eq:exp-G}
G_{ p }=2 \, \rho^{-2} + A + O(\rho).
\end{split}\end{equation}

\medskip 

\subsection{Approximate solutions}

For $p \in S^3$, fix a small number $r > 0$ and define in CR normal 
coordinates a function $F$ such that 
$$
  \begin{cases}
  F(z,t) = |z|^2 & \hbox{ for } \rho \leq r; \\
  F \equiv 0 & \hbox{ for } \rho \geq 2 r. 
  \end{cases}
$$
In this way, $F$ can be extended via cut-offs to all of $S^3$ as the zero function away from $p$, 
so $F$ can be written as 
\begin{equation}\label{eq:F}
 F(z,t) = |z|^2 + O(\rho^5). 
\end{equation}
For  $\l > 0$ large,  let us consider a test function in 
CR normal coordinates as follows 
\begin{equation}\label{eq:breve-var}
  \breve{\var}_\l = \frac{\l}{\left( 1 + \l^2 F + \l^4 \tilde{G} \right)^{\frac 12}}, 
\end{equation}
where  $\tilde{G} = G_p^{-2}$.

\begin{lem}\label{l:brevevar}
In CR normal coordinates one has the expansion 
$$
    L_b \breve{\var}_\l = \breve{\var}_\l^3 \left( 2 + O(\rho^3) +  \l^{-2} O(\rho^2)  \right)
       +  \breve{\var}_\l^5 \left[ - \frac 32 |z|^2 \rho^2 (4 + \l^2 |z|^2) A + O(\rho^5) + O(\l^2 \rho^7) \right].
$$
\end{lem}

\begin{pf}
By direct computation we have that 
$$
 \breve{\var}_{\l,1} = - \frac 12 \frac{\l^3}{\left( 1  + \l^2 F(z) + \l^4 \tilde{G} \right)^{\frac 32}} 
 \left[  F_{,1} + \l^2  \tilde{G}_{,1} \right] = - \frac 12 \breve{\var}_\l^3  \left[  F_{,1} + \l^2  \tilde{G}_{,1} \right], 
$$
and similarly for its conjugate. As a consequence, we have that 
$$
  \breve{\var}_{\l,1\bar 1} = - \frac 12 \breve{\var}_\l^3
   \left[   F_{,1 \bar{1}} + \l^2  \tilde{G}_{,1 \bar 1} \right] + \frac 34 
   \breve{\var}_\l^5 \left|  F_{,1} + \l^2 \tilde{G}_{,1} \right|^2, 
$$
which implies  
$$
  \Delta_b \breve{\var}_\l  = - \frac 12 \breve{\var}_\l^3
     \left[   \D_b F + \l^2  \D_{b} \tilde{G} \right] + 
     \frac 32 \breve{\var}_\l^5 \left|  F_{,1} + \l^2 \tilde{G}_{,1} \right|^2. 
$$
By direct computation one finds (with $G = G_p$)  
$$
  \D_b \tilde{G} = - 2 G^{-3} \D_b G + 12 G^{-4} G_{,1} G_{, \bar 1}. 
$$
We then deduce 
$$
  L_b \breve{\var}_\l = 2 \, \breve{\var}_\l^3
       \left[   \D_b F + \l^2 (12 G^{-4} G_{,1} G_{, \bar 1} - 2 G^{-3} \D_b G ) \right] - 
       6 \, \breve{\var}_\l^5 \left|  F_{,1} + \l^2 \tilde{G}_{,1} \right|^2 + R \, \breve{\var}_\l. 
$$
We can next write  
$$
  R \, \breve{\var}_\l = R \, \breve{\var}_\l^3  \left(   \l^{-2}  +  F + \l^2 G^{-3} G \right).
$$
Since $G$ satisfies $L_b G = 0$, we get some cancellation and find that 
$$
  L_b \breve{\var}_\l = \breve{\var}_\l^3 \left( 2 \D_b F +   \l^{-2} R  +  R \, F \right)
  + 6 \, \breve{\var}_\l^5 \left[ 4 G^{-4} G_{,1} G_{, \bar 1}  \left( 1  + \l^2 F(z) + \l^4 \tilde{G} \right) - 
  \left|  F_{,1} + \l^2 \tilde{G}_{,1} \right|^2 \right]. 
$$
Using some further cancellation we then obtain 
$$
  L_b \breve{\var}_\l = \breve{\var}_\l^3 \left( 2 \D_b F +  \l^{-2} R +  R \, F \right)
  + 6 \, \breve{\var}_\l^5 \left[ 4 G^{-4} G_{,1} G_{, \bar 1}  \left( 1 + \l^2 F(z) \right) - 
  F_{,1} F_{, \bar 1} - \l^2 (F_{,1} \tilde{G}_{, \bar 1} + \tilde{G}_{,1} F_{, \bar 1}) \right]. 
$$
From Proposition A.5 in \cite{CMY17} (where a different but analogous notation is used) one has that, in CR normal coordinates 
$$
  Z_1 = (1 + O(\rho^4)) \overset{\circ}{Z}_1 + O(\rho^4) \overset{\circ}{Z}_{\bar 1} + O(\rho^5) \frac{\pa}{\pa t}; 
$$
$$
 \omega_{1}^1 = O(\rho^3) dz + O(\rho^3) d \bar z + O(\rho^2) \overset{\circ}{\theta}, 
$$
see \eqref{eq:circleformulas}. 
By direct computation,  one  then has  
$$
  G_{,1} = -\frac{i \sqrt{2} \bar z}{(t+i |z|^2) \rho^2} + O(1); \qquad  \qquad F_{,1} = \frac{\bar z}{\sqrt{2}} 
  + O(\rho^4); 
  $$
  $$
     \tilde{G}_{,1} = \frac{2 \sqrt{2} \bar z (|z|^2+i t)}{\left(A \rho^2 +2\right)^3} + O(\rho^6);  
  \qquad \qquad  \D_b F = 1 + O(\rho^3). 
$$
Using these expressions in the above formula for $L_b \breve{\var}_\l$ one finally finds 
$$
    L_b \breve{\var}_\l = \breve{\var}_\l^3 \left( 2 +  \l^{-2} O(\rho^2) +  O(\rho^3) \right)
    + 6 \, \breve{\var}_\l^5 \left[ - \frac 14 |z|^2 \rho^2 (4 + \l^2 |z|^2) A + O(\rho^5) + O(\l^2 \rho^7) \right], 
$$
which is the desired result. 
\end{pf}

\medskip

If the contact  form $\th$ involved in the definition of CR normal coordinates 
writes as $\th = e^{2v} \hat{\th}$, setting 
\begin{equation}\label{eq:bar-var-l}
  \bar \var_{\l} = e^{-v} \breve{\var}_\l, 
\end{equation}
by the covariance property of the conformal sub-Laplacian one has that 
\begin{equation}\label{eq:cov-laws}
\breve{\var}_\l^4 \, {\th} \wedge d {\th} = 
  \bar {\var}_\l^4 \, \hat{\th} \wedge d \hat{\th}; \qquad \qquad 
  \breve{\var}_\l L_b^{(\th)} \breve{\var}_\l \, {\th} \wedge d{\th}  = \bar{\var}_\l  L_s \bar{\var}_\l \, \hat{\th} \wedge d \hat{\th}, \quad L_s = L_b^{(\hat \th)}. 
\end{equation}
  These imply the invariance  
 $$
    Q_{(s)} (\bar{\var}_\l ) = \frac{\int_{S^3} \bar{\var}_\l  L_s \bar{\var}_\l  \, \hat{\th} \wedge d \hat{\th}}{\left( \int_{S^3} \bar{\var}_\l ^4 \hat{\th} \wedge d \hat{\th} \right)^{\frac 12}} = 
    \frac{\int_{S^3} \breve{\var}_\l  L_b^{(\th)} \breve{\var}_\l  \, {\th} \wedge d {\th}}{\left( \int_{S^3} \breve{\var}_\l^4 {\th} \wedge d {\th} \right)^{\frac 12}}.
 $$ 
We then get  the following consequence of Lemma \ref{l:brevevar}, concerning the 
differential of $Q_{(s)}$ at $\bar{\var}_\l$.

\begin{cor}\label{eq:error-norm}
There exists a constant $C > 0$ such that for all $s$ small and $\l$ large one has the inequality $|d Q_{(s)}(\bar \var_{\l}) [v]| \leq \frac{C}{\l^2} \|v\|_{\mathfrak{S}^{1,2}}$ for 
every $v \in \mathfrak{S}^{1,2}(S^3)$. 
\end{cor}

\begin{pf}
By direct computation, for $v \in \mathfrak{S}^{1,2}(S^3)$, one has 
\begin{equation}\label{eq:direct}
d Q_{(s)}(\bar \var_{\l}) [v] =  \frac{2}{\left( \int_{S^3} \bar{\var}_\l^4 \, \hat{\th} \wedge d \hat{\th} \right)^{\frac{3}{2}}} 
  \int_{S^3}  \left[ \left( \int_{S^3} \bar{\var}_\l^4 \, \hat{\th} \wedge d \hat{\th} \right) L_s \bar{\var}_\l 
  - \left( \int_{S^3} \bar{\var}_\l L_s \bar{\var}_\l  \hat{\th} \wedge d \hat{\th} \right) \bar{\var}_\l^3 \right]  v \, \hat{\th} \wedge d \hat{\th}. 
\end{equation}
From \eqref{eq:cov-laws} and Lemma \ref{l:brevevar} it follows that 
$$
  \int_{S^3} \bar{\var}_\l L_s \bar{\var}_\l  \hat{\th} \wedge d \hat{\th} = 2  \int_{S^3} \breve{\var}_\l^4 \, {\th} \wedge d {\th} 
  +  \int_{S^3} \left[ \breve{\var}_\l^4 (O(\rho^3) + \l^2 O(\rho^2)) + \breve{\var}_\l^6 (O(\rho^4) + \l^2 O(\rho^6)) \right] {\th} \wedge d {\th}.  
$$
Using a change of variable it is possible then to show  
$$
\int_{S^3} \breve{\var}_\l^4 \, {\th} \wedge d {\th} =
   -  \int_{S^3} \breve{\var}_\l L_b^\th \breve{\var}_\l  {\th} \wedge d {\th} + O(\l^{-2}). 
$$
Therefore, inserting the latter estimate and the result of Lemma \ref{l:brevevar} into 
\eqref{eq:direct} we find that 
$$
  |d Q_{(s)}(\bar \var_{\l}) [v]| \leq  \int_{S^3}  
  \breve{\var}_\l^3  \left[O(\rho^2) + O(\l^{-2}) + \breve{\var}_\l^5 
    \left( O(\rho^4) + \l^2 O(\rho^6) \right) 
   \right] \, |v| \, \th \wedge d \th. 
$$
Applying H\"older's inequality we get that 
\begin{eqnarray*}
 & & |d Q_{(s)}(\bar \var_{\l}) [v]| \\ & \leq & \left[ O(\l^{-2}) + \left( \int_{S^3} \breve{\var}_\l^4 O(\rho^{\frac 83}) \right)^{\frac{3}{4}} 
  +  \left( \int_{S^3} \breve{\var}_\l^{20/3} O(\rho^{\frac{16}{3}}) \right)^{\frac{3}{4}}  
  + \l^2 \left( \int_{S^3} \breve{\var}_\l^{20/3} O(\rho^{8}) \right)^{\frac{3}{4}} \right] \|v\|_{\mathfrak{S}^{1,2}}, 
\end{eqnarray*}
where all integrals are computed w.r.t. the volume form $\th \wedge d \th$. By the expression of $\breve{\var}_\l$, 
all terms are integrable and of order $\l^{-2}$, which concludes the proof. 
\end{pf}

\medskip

\subsection{Expansion of the CR Sobolev quotient} We expand next the CR Sobolev quotient 
$Q_{(s)}$ on the  approximate solutions $\bar{\var}_\l$ in \eqref{eq:bar-var-l}, obtaining the following result.

  \begin{lem}\label{prop_functional_at_infinity} 
  Let $\breve{\var}_\l$ be defined in \eqref{eq:breve-var}. Then for $\l$ 
  large one has the expansion  
 $$
   Q_{(s)} (\bar{\var}_\l) =  4 \pi  + 
     48 \pi \frac{s^2}{\l^2} (1+o_s(1)) + O \left( \frac{1}{\l^3} \right).
 $$ 
  
  \end{lem}

 \begin{pf}
 We use \eqref{eq:cov-laws},   Lemma \ref{l:brevevar} and integrate: 
expanding the numerator in $Q_{(s)}$ we find that 
\begin{eqnarray*}
 \int_{S^3} \breve{\var}_\l L_b^{(\th)} \breve{\var}_\l \th \wedge d \th & = & 2 \int_{S^3} \breve{\var}_\l^4 \th \wedge d \th 
  - \frac{3}{2} A \int_{\H^1} |z|^2 (4 + \l^2 |z|^2) \rho^2 \overset{\circ}{\var}_\l^6 \, \overset{\circ}{\theta} \wedge d \overset{\circ}{\theta}  \\ 
  & + & \int_{S^3} \var^4 (O (\l^2 \rho^2) + O(\rho^3)) \th \wedge d \th + \int_{S^3} \var^6 (O(\rho^5) + O(\l^2 \rho^7))
  \th \wedge d \th, 
\end{eqnarray*}
where 
$$
  \overset{\circ}{\var}_\l = \frac{\l}{\left( 1 + \l^2 |z|^2 + \frac 14 \l^4 (|z|^4 + t^2) \right)^{\frac 12}}; 
  \qquad \quad (z,t) \in \H^1.  
$$

For the first term, which also appears in the above expression, we Taylor-expand $\tilde{G}$ as
$$
  \tilde{G} = \left( \frac{2 + A \rho^2}{\rho^2} \right)^{-2} = \frac{1}{4} \rho^4 (1 - A \rho^2) + O(\rho^8). 
$$
Therefore, $\breve{\var}_\l$ expands as  
\begin{eqnarray*}
 \breve{\var}_\l & = & \frac{\l}{\left( 1 + \l^2 (|z|^2 + O(\rho^5)) + \frac{1}{4} \l^4 [\rho^4 (1-A \rho^2) + O(\rho^8)] \right)^{\frac 12}} 
  \\ & = & \left( 1 + \frac{1}{8} \frac{A \rho^6 \l^4}{1 + \l^2 |z|^2 + \frac{1}{4} \l^4 \rho^4} 
  + O \left( \frac{\rho^{12} \l^8}{(1 + \l^4 \rho^4)^2}  \right)\right)  \overset{\circ}{\var}_\l 
  \\ & = & \overset{\circ}{\var}_\l  + \frac 18 A \rho^6 \l^2 \overset{\circ}{\var}_\l^3 
  + O \left( \frac{\rho^{12} \l^8}{(1 + \l^4 \rho^4)^2}  \right) \overset{\circ}{\var}_\l. 
\end{eqnarray*}
Taylor-expanding the integral of the fourth power of $\breve{\var}_\l$ and using a 
change of variable we get that 
\begin{eqnarray*}
  \int_{S^3} \breve{\var}_\l^4 \th \wedge d \th  & = & \int_{\H^1} \overset{\circ}{\var}_\l^4 \overset{\circ}{\theta} \wedge d \overset{\circ}{\theta}  + \frac{1}{2} A \l^2 \int_{\H^1} \rho^6 \overset{\circ}{\var}_\l^6 \overset{\circ}{\theta} \wedge d \overset{\circ}{\theta} 
  + O(1/\l^3). 
\end{eqnarray*}
Hence, using the fact that $\int_{\H^1} \overset{\circ}{\var}_\l^4 \overset{\circ}{\theta} \wedge d \overset{\circ}{\theta}$ 
is independent of $\l$,  $Q_{(s)} (\bar{\var}_\l)$ becomes 
$$
  \frac{2 \left(  \int_{\H^1} \overset{\circ}{\var}_{\sqrt{2}}^4 \, \overset{\circ}{\theta} \wedge d \overset{\circ}{\theta}  + \frac{1}{2} A \l^2 \int_{\H^1} \rho^6 \overset{\circ}{\var}_\l^6 \overset{\circ}{\theta} \wedge d \overset{\circ}{\theta}  \right)  - 
  \frac{3 A}{2 } \int_{\H^1} |z|^2 \rho^2 (4 + \l^2 |z|^2) \overset{\circ}{\var}_1^6 \, \overset{\circ}{\theta} \wedge d \overset{\circ}{\theta} }{\left( \int_{\H^1} \overset{\circ}{\var}_{\sqrt{2}}^4 \, \overset{\circ}{\theta} \wedge d \overset{\circ}{\theta}  + \frac{1}{2} A \l^2 \int_{\H^1} \rho^6 \overset{\circ}{\var}_\l^6 \overset{\circ}{\theta} \wedge d \overset{\circ}{\theta}  \right)^{\frac 12}} + O(1/\l^3). 
$$
We can expand the denominator in the latter expression as 
\begin{eqnarray*}
  & & \left( \int_{\H^1} \overset{\circ}{\var}_{\sqrt{2}}^4 \, \overset{\circ}{\theta} \wedge d \overset{\circ}{\theta}  + \frac{1}{2} A \l^2 \int_{\H^1} \rho^6 \overset{\circ}{\var}_\l^6 \overset{\circ}{\theta} \wedge d \overset{\circ}{\theta}   \right)^{- \frac 12} 
  \\ & = & \left( \int_{\H^1} \overset{\circ}{\var}_{\sqrt{2}}^4 \, \overset{\circ}{\theta} \wedge d \overset{\circ}{\theta}  \right)^{- \frac 12} \left( 1 + \frac{\frac{1}{2} A \l^2 \int_{\H^1} \rho^6 \overset{\circ}{\var}_\l^6 \overset{\circ}{\theta} \wedge d \overset{\circ}{\theta}   }{\int_{\H^1} \overset{\circ}{\var}_1^4 \, \overset{\circ}{\theta} \wedge d \overset{\circ}{\theta} } \overset{\circ}{\theta} \wedge d \overset{\circ}{\theta}  \right)^{-\frac 12} \\ & = & \left( \int_{\H^1} \overset{\circ}{\var}_{\sqrt{2}}^4 \, \overset{\circ}{\theta} \wedge d \overset{\circ}{\theta}  \right)^{- \frac 12} \left( 1 - \frac{1}{4} A \l^2 
  \frac{ \int_{\H^1} \rho^6 \overset{\circ}{\var}_\l^6 \, \overset{\circ}{\theta} \wedge d \overset{\circ}{\theta}  }{\int_{\H^1} \overset{\circ}{\var}_1^4 \, \overset{\circ}{\theta} \wedge d \overset{\circ}{\theta} } \right) + O(1/\l^3), 
\end{eqnarray*}
which gives 
\begin{eqnarray*}
  & & Q_{(s)} (\bar{\var}_\l) \\ & = & \left( \int_{\H^1} \overset{\circ}{\var}_{\sqrt{2}}^4 \, \overset{\circ}{\theta} \wedge d \overset{\circ}{\theta}  \right)^{- \frac 12} \left[ 2 \int_{\H^1} \overset{\circ}{\var}_{\sqrt{2}}^4 \, \overset{\circ}{\theta} \wedge d \overset{\circ}{\theta} 
  + \frac{1}{2} {A} {\l^2} \int_{\H^1} \rho^6 \overset{\circ}{\var}_\l^6 \, \overset{\circ}{\theta} \wedge d \overset{\circ}{\theta}  -  \frac{3}{2} {A}\int_{\H^1} |z|^2 \rho^2 
  (4 + \l^2 |z|^2) \overset{\circ}{\var}_1^6 \, \overset{\circ}{\theta} \wedge d \overset{\circ}{\theta}  \right] \\ & + & O \left( \frac{1}{\l^3} \right), 
\end{eqnarray*}
equivalent to 
\begin{eqnarray}\label{eq:llllll} \nonumber
  Q_{(s)} (\bar{\var}_\l) & = & \left( \int_{\H^1} \overset{\circ}{\var}_{\sqrt{2}}^4 \, \overset{\circ}{\theta} \wedge d \overset{\circ}{\theta}  \right)^{- \frac 12} \left[ 2 \int_{\H^1} \overset{\circ}{\var}_{\sqrt{2}}^4 \, \overset{\circ}{\theta} \wedge d \overset{\circ}{\theta} 
  - \frac{1}{2} A \int_{\H^1} \left( 3 |z|^2 
    (4 + \l^2 |z|^2) - \l^2 \rho^4 \right) \rho^2 \overset{\circ}{\var}_\l^6 \, \overset{\circ}{\theta} \wedge d \overset{\circ}{\theta}  \right] \\ & + & O \left( \frac{1}{\l^3} \right).
\end{eqnarray}
The computation on page 177 in \cite{JLCRYam} (where $\theta_1$ in their notation equals  $2 \hat \theta$)
shows that  $\overset{\circ}{\var}_{\sqrt{2}}^4$ is the scaling factor for the volume of the Cayley map. 
Recalling that $\hat \theta \wedge d \hat \theta$ is twice the (induced) Euclidean volume on $S^3$, 
this implies 
\begin{equation}\label{eq:vol-S3}
\int_{\H^1} \overset{\circ}{\var}_{\sqrt{2}}^4 \, \overset{\circ}{\theta} \wedge d \overset{\circ}{\theta}  
= \int_{S^3} \hat \theta \wedge d \hat \theta = 4 \pi^2. 
\end{equation}
We now make the following change of variables $\l z \mapsto \sqrt{2} z'$, $\l^2 t \mapsto 2 t'$, 
and notice that 
$$
  \overset{\circ}{\var}_\l (z,t) = \frac{\l}{\sqrt{2}} \frac{\sqrt{2}}{\left( (1 + |z'|^2)^2 + (t')^2  \right)^{\frac 12}} 
  = \frac{\l}{\sqrt{2}} \overset{\circ}{\var}_{\sqrt{2}} (z',t'). 
$$
In this way we have  
$$
  \int_{\H^1} \left( 3 |z|^2 
      (4 + \l^2 |z|^2) - \l^2 \rho^4 \right) \rho^2 \overset{\circ}{\var}_\l^6 \, \overset{\circ}{\theta} \wedge d \overset{\circ}{\theta}  (z,t) = \frac{4}{\l^2} \int_{\H^1} \left( 3 |z'|^2 
            (2 + |z'|^2) -  (\rho')^4 \right) (\rho')^2 \overset{\circ}{\var}_{\sqrt{2}}^6 \, \overset{\circ}{\theta} \wedge d \overset{\circ}{\theta}  (z',t'). 
$$
As one can check by direct computations, the primitive w.r.t. $t$ of the integrand is 
\begin{eqnarray*}
  -\frac{3 \left(|z'|^6+8 |z'|^4+19 |z'|^2+8\right) |z'|^6 \log \left(\frac{t^2 \left(|z'|^4+4 |z'|^2+2\right)-2 t \left(|z'|^2+1\right) \sqrt{2 |z'|^2+1} \sqrt{t^2+|z'|^4}+\left(|z'|^2+1\right)^2 |z'|^4}{|z'|^4 \left(t^2+\left(|z'|^2+1\right)^2\right)}\right)}{ 
  2 \left(|z'|^2+1\right)^5\left(2 |z'|^2+1\right)^{3/2}} \\
  -\frac{t \sqrt{t^2+|z'|^4} \left(t^2 \left(3 |z'|^8+4 |z'|^6-17 |z'|^4-4 |z'|^2+2\right)+\left(|z'|^3+|z'|\right)^2 \left(3 |z'|^6-8 |z'|^4-55 |z'|^2-24\right)\right)}{\left(|z'|^2+1\right)^4\left(2 |z'|^2+1\right) \left(t^2+\left(|z'|^2+1\right)^2\right)^2}. 
\end{eqnarray*} 
As a consequence, we deduce that 
\begin{eqnarray*}
& & \int_\R \left( 3 |z'|^2 
            (2 + |z'|^2) -  (\rho')^4 \right) (\rho')^2 \overset{\circ}{\var}_{\sqrt{2}}^6 \, dt \\ 
 & = & \frac{-3 \left(|z'|^6+8 |z'|^4+19 |z'|^2+8\right) |z'|^6 \log \left(\frac{|z'|^4-2 \left(\sqrt{2 |z'|^2+1}-2\right) |z'|^2-2 \sqrt{2 |z'|^2+1}+2}{|z'|^4+2 \left(\sqrt{2 |z'|^2+1}+2\right) |z'|^2+2 \left(\sqrt{2 |z'|^2+1}+1\right)}\right)}{2 \left(|z'|^2+1\right)^5 \left(2 |z'|^2+1\right)^{3/2}} \\ &  - & \frac{4 \left(|z'|^2+1\right) \sqrt{2 |z'|^2+1} \left(3 |z'|^8+4 |z'|^6-17 |z'|^4-4 |z'|^2+2\right)}{2 \left(|z'|^2+1\right)^5 \left(2 |z'|^2+1\right)^{3/2}}.
\end{eqnarray*}
Multiplying this quantity by $2 \pi |z'|$, its primitive w.r.t. $|z'|$ is 
$$
-\frac{\pi  \left(3 \left(|z'|^2+2\right) |z'|^8 \log \left(\frac{|z'|^4-2 \left(\sqrt{2 |z'|^2+1}-2\right) |z'|^2-2 \sqrt{2 |z'|^2+1}+2}{|z'|^4+2 \left(\sqrt{2 |z'|^2+1}+2\right) |z'|^2+2 \left(\sqrt{2 |z'|^2+1}+1\right)}\right)+4 \sqrt{2 |z'|^2+1} \left(|z'|^6+6 |z'|^4+6 |z'|^2+1\right)\right)}{2 \left(|z'|^2+1\right)^4 \sqrt{2 |z'|^2+1}},
$$
whose difference between the values $|z'| \to + \infty$ and $|z'| = 0$ is $8 \pi$. 
Therefore, recalling that the volume form $\overset{\circ}{\theta} \wedge d \overset{\circ}{\theta}$ 
is four times the Euclidean one, we obtain that
$$
 \int_{\H^1} \left( 3 |z|^2 
     (4 + \l^2 |z|^2) - \l^2 \rho^4 \right) \rho^2 \overset{\circ}{\var}_\l^6 \, \overset{\circ}{\theta} \wedge d \overset{\circ}{\theta}   = 32 \pi. 
$$ 
Recalling \eqref{eq:vol-S3} and the fact that $A = - \frac{3}{2} s^2 (1+o_s(1))$,  from \eqref{eq:As2} 
and \eqref{eq:llllll}  
we deduce that 
$$
   Q_{(s)} (\breve{\var}_\l)  = 4 \pi - \frac{32\pi A}{\l^2} s^2 (1+o_s(1)) + O\left( \frac{1}{\l^3} \right) = 
  4 \pi + \frac{48\pi}{\l^2}  s^2 (1+o_s(1)) + O\left( \frac{1}{\l^3}  \right). 
$$
This concludes the proof.
\end{pf}

\medskip 

\subsection{Conclusion}

We can use the observation in Remark 5.3, to perform the contraction argument in Proposition 
\ref{p:ls-reduction} starting from $\breve{\var}_\l$ instead of from $\var_{\l}$ only. 
Given the improved accuracy in Lemma \ref{p:expansion}, the contraction can be performed 
in a ball of radius $O(\frac{1}{\l^2})$ in $\mathfrak{S}^{1,2}(S^3)$, yielding a corresponding correction 
$\breve{w}_\l$ of that order. By Lemma \ref{p:expansion} and the smoothness of $Q_{(s)}$, we then have 
similarly to Subsection \ref{ss:concl1}  
$$
  Q_{(s)} (\breve{\var}_\l+ \breve{w}_\l) = Q_{(s)} (\breve{\var}_\l) +  
  O(\| \breve{w}_\l \|^2)  = Q_{(s)} (\breve{\var}_\l) + O \left( \frac{1}{\l^4} \right).  
$$
By uniqueness of the fixed point, it must be $\breve{\var}_\l+ \breve{w}_\l = \var_{\l} + w_\l$, 
so from Lemma \ref{prop_functional_at_infinity} we get that 
\begin{equation}\label{eq:fin-fin-fin}
  Q_{(s)}( \varphi_{\l} + w_\l) =  4 \pi + 48 \pi \frac{s^2}{\l^2} (1 + o_s(1)) + O \left( \frac{1}{\l^3} \right).
\end{equation}
Notice that 
$$
  Q_{(0)}( \varphi_{\l} + w_\l) = Q_{(0)}( \varphi_{\l}) \equiv 4 \pi, 
$$
 and therefore
the term $O(\frac{1}{\lambda^{3}})$ appearing in \eqref{eq:fin-fin-fin}
is identically zero for $s = 0$, even and smooth in $s$. It therefore must be of the form 
$
  O \left(  \frac{s^2}{\l^3} \right). 
$
Hence the statement of the proposition holds true.

\bibliography{CR-Sobolev_references}

\bibliographystyle{plain}

\end{document}